\documentclass{article}
\usepackage{amsfonts}

\newenvironment{proof}{{\bf Proof}:}{\vskip 5mm }

\newtheorem{lemma}{Lemma}
\newtheorem{proposition}[lemma]{Proposition}

\newtheorem{theorem}[lemma]{Theorem}
\newtheorem{cor}[lemma]{Corollary}
\newtheorem{conjecture}{Conjecture}

\begin{document}
\title{Hodge structure on the fundamental group and its application to p-adic integration}
\author{Vadim Vologodsky}
\date{}
\maketitle
\begin{abstract}
We study the unipotent completion $\Pi^{DR}_{un}(x_0, x_1, X_K)$ of the de Rham fundamental groupoid [De]
of a smooth algebraic variety over a local non-archimedean field $K$ of characteristic $0$.
  We show that the vector space  
 $\Pi^{DR}_{un}(x_0, x_1, X_K)$ carries a certain additional structure. That is a $\mathbb{Q}^{ur}_p$-space 
 $\Pi_{un}(x_0, x_1, X_K)$ equipped  with a $\sigma$-semi-linear operator $\phi $, a linear operator $N$ satisfying the relation
   $N \phi = p \phi N $ and a weight filtration $W_{\bullet}$ together with a canonical isomorphism 
  $\Pi^{DR}_{un}(x_0, x_1, X_K)\otimes _K \overline K \simeq \Pi_{un}(x_0, x_1, X_K)\otimes _{\mathbb{Q}^{ur}_p} \overline K$.
We prove that an analog of the Monodromy Conjecture holds for  $\Pi_{un}(x_0, x_1, X_K)$. 

As an application, we show that the vector space   $\Pi^{DR}_{un}(x_0, x_1, X_K)$ possesses a distinguished element. In the other words,
given a vector bundle $E$ on $X_K$ together with a unipotent integrable connection, we have {\sf a canonical} isomorphism
$E_{x_0}\simeq E_{x_1}$ between the fibers. The latter construction is a generalization 
of  Colmez's $p$-adic integration ($rk\, E=2$) and
Coleman's $p$-adic iterated integrals ($X_K$ is a curve with good reduction). 

In the second part we prove that, if $X_{K_0}$ is a smooth  variety over an unramified extension of $\mathbb{Q}_p$   with
 good reduction  and $r \leq \frac{p-1}{2}$ then there is a canonical isomorphism $\Pi^{DR}_{r}(x_0, x_1, X_{K_0})\otimes B_{DR} \simeq
\Pi^{et}_{r}(x_0, x_1,  X_{\overline K_0})  \otimes B_{DR}$ compatible with the action of Galois group. 
( $\Pi^{DR}_{r}(x_0, x_1, X_{K_0})$ stands for the level $r$ quotient of  $\Pi^{DR}_{un}(x_0, x_1, X_K)$).
In particularly, it implies the Crystalline Conjecture for the fundamental group [Shiho] (for $r \leq \frac{p-1}{2}$) .
\end{abstract}

\tableofcontents

\section{Introduction} 

{\bf 1.1. Notation} Throughout this paper  $k$ stands for 
a finite field of characteristic $p>2$,
$W(k)$ is the ring of Witt vectors, $K_0$ is its field of fractions, 
$K$ is a finite extension of $K_0$ with the ring of integers $R\subset K$
and $e$ is its ramification index.

{\bf 1.2.} For a  smooth scheme $X$ over a field of characteristic $0$ we denote by  ${\cal D}{\cal M} (X)$ 
 the category of  vector bundles  on $X$ together with an integrable connection.

An object $E$ of the latter category is called unipotent of level $ r$  
if there exists  a filtration 
$$ 0=E_0 \subset E_1 \subset E_2 \subset \cdot \cdot \cdot \subset E_r=E $$
where  $E_i$ are subobjects,  such that  the  quotients $E_i/E_{i-1}$ are trivial.
The latter means that the vector bundle  $E_i/E_{i-1}$ is generated by its global parallel 
sections.

We denote the full subcategory of  ${\cal D}{\cal M} (X) $ 
consisting of unipotent bundles 
of level $r$ by  ${\cal D}{\cal M}_{r} (X) $ and put
${\cal D}{\cal M}_{un} (X):= \bigcup _r {\cal D}{\cal M}_{r} (X)$
The latter category carries a natural tensor structure.

{\bf 1.3.}  Let $X_{K_0}$ be a smooth scheme  finite type over $K_0$.
 For a  pair of points $x_0,x_1\in X_{K_0}(\overline K) $
we denote by $\pi^{et}_1(X_{\overline K}, x_0, x_1)$  the
set of isomorphisms between the corresponding functors from the category 
of finite etale schemes over $X_{\overline K}=X\times spec \overline K$
to the category of finite sets.

{\bf 1.4.} To state the first main result, we assume that the scheme $X_{K_0}$ 
 has {\it a good compactification over} $spec \, W(k)$. That is 
 a smooth proper scheme $\overline X$ over $spec\, W(k)$ with 
  a divisor $Z\subset \overline X$  with  normal crossings relative to 
  $spec \, W(k)$ and an isomorphism   
$$X_{K_0} \simeq (\overline X - Z)\times_{spec \, W(k)} spec \, K_0 $$ 
  
{\bf 1.5.} Assume that $r \leq \frac{p-1}{2} $.  For any object
$E$ of the category  ${\cal D}{\cal M}_{ r} (X_{\overline K}) $, a pair of points 
$x_0,x_1\in X(\overline K) $ and an etale path $\lambda \in \pi^{et}_1(X_{\overline K}, x_0, x_1)$
we construct a parallel translation:
\begin{equation}
T_{\lambda}: E_{x_0} \otimes _{\overline K} B_{DR} \simeq  E_{x_1} \otimes_{\overline K} B_{DR}
\label{pt}
\end{equation}
Here $B_{DR}$ is the certain universal ring of periods introduced by Fontaine in [Fo]. 
 
{\bf 1.6.}   To formulate  the result more precisely it is convenient to introduce 
a certain category, which is, in fact,  a modification of the fundamental groupoid [De].

A  point $x \in X_{K_0} (\overline K) $ defines a fiber functor 
$$ F_x : \; {\cal D}{\cal M}_{ r} (X_{\overline K}) \longrightarrow Vect_{\overline K} $$ 
to the  category of vector spaces over $ \overline K $.

Consider the category 
${\cal P}^{DR}_{r}(X_{\overline K})$,
whose objects $\underline {x} $  are points of
 $ X_{K_0} (\overline K) $ and  whose set morphisms between two points $\underline {x_0}, \underline {x_1} $ is the vector 
space of  morphisms between the corresponding fiber functors:
               $$Mor_{{\cal P}^{DR}_{r}(X_{\overline K})}  (\underline {x_0}, \underline {x_1})= Mor(F_{x_0},F_{x_1})$$         

In the same way we define the category ${\cal P}^{et}_{r}(X_{\overline K})$,
 simply by replacing vector bundles with connections by  unipotent 
etale locale systems of ${\mathbf Q}_p$-vector spaces.

The categories ${\cal P}^{DR}_{r}(X_{\overline K})$ and ${\cal P}^{et}_{\leq r}(X_{\overline K})$ 
carry  a natural action of the Galois group $Gal(\overline K /K_0)$.  
 
{\bf 1.7. }. Given a homomorphism of rings $C\to C^{\prime} $ and a 
category ${\cal A}$ with a structure  of a $C$-module on $Mor_{{\cal A}}(*;*)$
compatible with the composition, we denote by ${\cal A}\otimes _C C^{\prime}$
the category whose set of  objects is the same but whose set of
 morphisms   between objects $V$ and $W$ is 
$$ Mor_{{\cal A}}(V;W)
\otimes _C  C^{\prime} $$       
 
{\bf 1.9. Theorem A}.
{\it Assume that $r \leq \frac{p-1}{2}$. Then there exists a functor:
$$ I_{ X_{\overline K}}
:\,{\cal P}^{et}_{ r}(X_{\overline K}) \otimes_{Q_p} B_{DR} \longrightarrow
 {\cal P}^{DR}_{r}(X_{\overline K}) \otimes_{\overline K} B_{DR} $$
identical on objects and satisfying the following properties:

a) $I_{ X_{\overline K}}$ establishes an equivalence between the categories.

b) The functor $I_{ X_{\overline K}}$ is equivariant with respect to the action of the Galois group 
$Gal(\overline K/K_0)$.

c)Let  $f:\overline X \to \overline Y $ be a morphism. 
Then we have:
     $$ f_*^{DR}\circ I_{ X_{\overline K}}= I_{ Y_{\overline K}} \circ  f_*^{et} $$  
Here $ f_*^{DR}$ and $ f_*^{et}$ stand for the functors:
$$ f_*^{DR}:{\cal P}^{DR}_{r}(X_{\overline K})\longrightarrow
 {\cal P}^{DR}_{r}(Y_{\overline K})$$
$$ f_*^{et}:{\cal P}^{et}_{r}(X_{\overline K})\longrightarrow
 {\cal P}^{et}_{r}(Y_{\overline K})$$ }           

{\bf 1.10. Remark.} The property c) holds for any morphism $f: X_{\overline K} \to  Y_{\overline K} $
of smooth schemes over $spec \,\overline K $, which have a good compactification (as defined in 1.4).

In particularly,  the functor
$ I_{ X_{\overline K}}$ depends only on  the scheme  $X_{\overline K}$ itself and not on the choice
of a good model.

This is a true fact but it is not proven in this paper.

{\bf 1.11.} Unfortunately we were not able to  find a direct construction of this functor; 
our proof is very implicit. It is based on the theory of variations 
of p-adic Hodge structures  developed by Faltings [Fa]. ( We use this term for the category
${\cal M}{\cal F}^{\nabla}(\overline X)\otimes \mathbb{Q}$ of certain filtered $F$-crystals 
introduced in {\it loc. cit.})  

{\bf 1.12.} The proof of Theorem 1 consists of the following steps.

First, we construct the parallel translation (\ref{pt}) for a certain universal unipotent local system.

Let $\Pi ^{DR}_r(X_{K_0})$ be the vector bundle on $X_{K_0} \times X_{K_0} $ 
together  with a unipotent connection of level $r$ 
characterized by the following property:

given another object $E$ of the category ${\cal D}{\cal M}_{r} (X_{K_0} \times X_{K_0})$, we have a canonical isomorphism:
\begin{equation} 
Hom _{{\cal D}{\cal M} (X_{K_0} \times X_{K_0}) }(\Pi ^{DR}_r(X_{K_0}); E) \simeq  H^0_{DR} (\Delta ; \Delta ^* E )
\label{podlost}
\end{equation}
where $\Delta : X_{K_0} \to X_{K_0} \times X_{K_0}$ is the diagonal embedding.

The fiber of $\Pi ^{DR}_r(X_{K_0})$ over a point $(x_0,x_1)$ is canonically 
identified with  $Mor(F_{x_0},F_{x_1})$.  

The identity morphism $Id:\, \Pi ^{DR}_r(X_{K_0}) \to \Pi ^{DR}_r(X_{K_0}) $ gives rise to a parallel section 
\begin{equation}
 \mathbf{1}: \,
  {\cal O}_{\Delta} \hookrightarrow \Delta ^* \Pi ^{DR}_r(X_{K_0})
\label{int} 		
\end{equation}
 
{\bf 1.13. } The main step is a construction of a variation of $p$-adic Hodge structure on  $\Pi ^{DR}_r(X_{K_0})$.
The analogus construction  over the field of complex numbers is well known. It is based on an interpretation of the fibers
in terms of homology groups of a certain simplicial scheme. In our case, we found a different construction
 based on a simple linear algebra argument. It works over $\mathbb{C}$ as well.     

We proof by induction on $r$, that there exists a unique variation of Hodge structure on $\Pi ^{DR}_r(X_{K_0})$ such that  
the map (\ref{podlost}) is a morphism of Hodge structures, for any unipotent level $r$ variation $E$. For this we consider the following exact sequence:
$$0 \longrightarrow I \longrightarrow \Pi ^{DR}_r(X_{K_0})\longrightarrow \Pi ^{DR}_{r-1}(X_{K_0})\longrightarrow 0$$
By induction hypothesis the boundary terms are endowed with Hodge structures. We show that there is a unique Hodge structure on
$\Pi ^{DR}_r(X_{K_0})$ satisfying the following property: the map (\ref{int}) is a morphism of Hodge structures.

 Given  a point $x_0 \in X_{K_0}(K_0)$,  we denote by $\Pi ^{DR}_r(X_{K_0}, x_0)$ the restriction 
of the variation of Hodge structure $\Pi ^{DR}_r(X_{K_0})$ to the fiber $X_{K_0}\times x_0$. 
 
{\bf 1.14.} Next, we make use of the functor
\begin{equation}
\mathbf{D}: \; {\cal M}{\cal F}^{\mathbb{Q}}(\overline X) \longrightarrow   Sh^{et}_{\mathbb{Q}_p}(X_{K_0})
\label{is}
\end{equation}
from the category of variations of Hodge structures to the category of etale locally constant
sheaves on $X_{K_0}$. The latter was constructed by Faltings.
 
We proof that for any point $x \in X_{K_0}(\overline K)$ and a variation $E$ there is a canonical isomorphism of fibers
\begin{equation}
 (\mathbf{D}(E))_x \otimes B_{DR} \simeq E_x \otimes B_{DR}
\label{isom}
\end{equation}
(The latter is not entirely obvious if the point reduces to infinity in the special fiber).  

Applying the (\ref{is}) and (\ref{isom}) to $\Pi _r^{DR}(X_{K_0},x)$ we obtain a  
parallel translation:
$$T_{\lambda}: \Pi _r^{DR}(X_{K_0},x_0)_{x_0} 
\otimes  B_{DR} \simeq  \Pi _r^{DR}(X_{K_0},x_0)_{x_1} \otimes  B_{DR}$$ 

In particularly, we constructed an element 
$$T_{\lambda}(\mathbf{1}) \in \Pi _r^{DR}(X_{K_0},x_0)_{x_1} \otimes  B_{DR}=
Mor_{{\cal P}^{DR}_{r}(X_{\overline K})}(\underline x_0;\underline x_1)\otimes  B_{DR}$$

It gives rise to a morphism:
$$ Mor_{{\cal P}^{et}_{r}(X_{\overline K})}(\underline x_0;\underline x_1)\otimes  B_{DR} \to 
Mor_{{\cal P}^{DR}_{r}(X_{\overline K})}(\underline x_0;\underline x_1)\otimes  B_{DR}$$
Finally, we show that the latter is an isomorphism.

{ \bf 1.15. } To formulate our second main result we let $X_K$ be a smooth geometrically
connected scheme  finite over $K$.

 Denote by $K_{st} $  the ring of polynomials $K[l(p)]$ in a formal variable $l(p)$.
Let  $N:K_{st}\to K_{st} $ be the derivation such that $N(x)= 0$ for any $x \in K$ and
$N(l(p))=1$.
     
We define a map
          $$Log: K^* \to K_{st}  $$  
to be the unique homomorphism given by the series
$$ Log\, x =\sum_{i=1}^{\infty}\frac{(1-x)^i}{i}$$
on a neighborhood of $0\in K$ and satisfying $Log \, p = l(p) $.   
 
 Given a vector bundle $E$ on $X_K$ together with a unipotent integrable connection and a pair of 
points $x_0, x_1 \in X_K(K)$, we construct a {\it canonical} isomorphism
\begin{equation}
   C_{x_0;x_1;X_{K}}(E): E_{x_0} \otimes _K K_{st} \simeq E_{x_1} \otimes _K K_{st} 
\label{canpar}
\end{equation}
We would like to stress, that the latter does not depend on any additional choices like that of an etale path.
         
{\bf 1.16.}  To describe some properties of the canonical parallel translation  we need a few  more notations :  denote by
$$ F_x \otimes K_{st}:{\cal D}{\cal M}_{un} (X_K) \longrightarrow  Mod_ {K_{st}}$$
the fiber functor to the category of $ K_{st} $-modules. Let 
$$Mor^{\otimes}( F_{x_0} \otimes K_{st}; F_{x_1} \otimes K_{st}) \subset 
Mor( F_{x_0} \otimes K_{st}; F_{x_1} \otimes K_{st}) $$   
be the set of all morphisms $C$ satisfying with following properties:

i) $ C({\cal{O}}_{X_K} )= Id $

ii)$ C(E^*) = (C(E))^*:  
E^*_{x_0} \otimes _K K_{st} \simeq E^*_{x_1} \otimes _K K_{st} $.

( Here $E^*$ is the dual object).

iii)$C(E \otimes E^{\prime})=C(E)
\otimes C(E^{\prime})$

for any objects $E$ and $E^{\prime}$.
    
Given a morphism of smooth connected schemes 
\begin{equation}
 f:X_K \to Y_K
\label{durak}
\end{equation} 
we denote by $f_*$
the canonical map  
$$ Mor^{\otimes}( F_{x_0} \otimes K_{st}; F_{x_1} \otimes K_{st}) \longrightarrow
   Mor^{\otimes}( F_{f(x_0)} \otimes K_{st}; F_{f(x_1)} \otimes K_{st})$$

The map $N:\, K_{st} \to   K_{st} $ induces a $K$-linear endomorphism (denoted by the same latter
$N$) of the space $Mor( F_{x_0} \otimes K_{st}; F_{x_1} \otimes K_{st}) $.
 
{\bf 1.17. Theorem B.} 
{\it There exists a canonical parallel translation i.e. an element 
$$C_{x_0;x_1;X_{K}} \in Mor^{\otimes}(F_{x_0} \otimes K_{st};F_{x_1} \otimes K_{st})  $$  
satisfying the following properties:

1)For any triple $x_0,x_1,x_2 \in X(K) $          
 $$C_{x_1;x_2;X_{K}} \circ C_{x_0;x_1;X_{K}}=  C_{x_0;x_2;X_{K}}$$

2)For any object  $E$ of the category ${\cal D}{\cal M}_{un} (X_K)$ the parallel translation
(\ref{canpar}) is locally analytic in the variables $x_0$ and $x_1$.

3) For any morphism (\ref{durak}) $f_*(C_{x_0;x_1;X_{K}})=C_{f(x_0);f(x_1);Y_{K}}$

4)Let $K\subset L$ be a finite field extension. Then  $C_{x_1;x_2;X_{K}}=C_{x_1;x_2;X_{L}}$.

5)Let $f$ be a non-vanishing function on $X_K$,
 $E$ be the trivial vector bundle  with a basis $e_{-1}, e_{0}$ 
We define a connection $\nabla : E \to E\otimes \Omega^1_X $    
 by: 
 $$\nabla(s)= A\, s + ds $$
where      
 $$
A=
\left(\begin{array}{cc}
0 & \omega  \\
0 & 0 
\end{array}\right)
$$
Here  $\omega $ is the differential form  $ \frac{df}{f}$. Then
the parallel translation 
$$
C_{x_0;x_1;X_{K}}(E): \;  E_{x_0} \otimes _K K_{st}  
\longrightarrow E_{x_1} \otimes _K K_{st}
$$
is given by the matrix

$$
\left(\begin{array}{cc}
1 & c \\
0 & 1 
\end{array}\right)
$$
where $c=Log \, f(x_0) - Log \, f(x_1)$.
6)Suppose that $X_K$ is proper. Then
$$ (N^r(C_{x_0;x_1;X_{K}}))(E)=0 $$
for any object $E$ of the category  ${\cal D}{\cal M}_{ r+1} (X_K) $ }

{\bf 1.18.} We also propose the following conjecture.
\begin{conjecture}\label{con}
A collection 
$$ C_{x_0;x_1;X_{K}} \in Mor^{\otimes}(F_{x_0} \otimes K_{st};F_{x_1} \otimes K_{st}) $$
satisfying the above properties is unique.
\end{conjecture}
{\bf 1.19.} The above result generalizes the construction of iterated integrals on a curve with good reduction
by R.Coleman [Co], on the one hand, and, on the other hand, Colmez's theory of $p$-adic integration
[Colmez]. The former is essentially equivalent to our construction of the canonical parallel translation
in the case when $X_K$ is a curve with a good reduction, 
while the second author treated the case of 
arbitrary smooth scheme $X_K$, but when $E$ is of rank $2$. 
Conjecture \ref{con} is also proven by Colmez in rank $2$ case.
   
{\bf 1.20.} Let us explain the idea of our construction of  $C_{x_0;x_1;X_{K}}$.

First, we show that the  space  $\Pi^{DR} _r (x_0;x_1;X_{K})= Mor(F_{x_0},F_{x_1})$ possesses a certain additional structure.
 That is a  vector space $\Pi^{ur}_r (x_0;x_1;X_{K})$ over the maximal unramified extension $K^{ur}_0$,  together with an action of
a nilpotent linear operator $N$, an invertible $Fr$-linear operator $\phi$, 
satisfying the following relation:
$$N\cdot \phi= p \phi \cdot N$$
and an isomorphism 
$$\Pi^{ur} _r (x_0;x_1;X_{K}) \otimes _{K^{ur}_0} \overline K_{st} \simeq \Pi^{DR} _r (x_0;x_1;X_{K})
\otimes _{K} \overline K_{st}$$
commuting with the action of $N$.
(Compare with  a conjecture by Fontaine:  according to this conjecture  the de Rham cohomology of any smooth variety over $K$ possesses such a structure.) 

{\bf 1.21.}It turns out that if $X_K$ has a {\it good compactification}
 over $spec \; R$ (see 1.4 for the definition), the  space $(\Pi^{ur} _r (x_0;x_1;X_{K})^{\phi}$ of invariants   
is one-dimensional over $\mathbb{Q}_p$ and, moreover,
it has a unique element $C_{x_0;x_1;X_{K}}$, characterised by the property:
\begin{equation}
 C_{x_0;x_1;X_{K}}({\cal O}_{X_K})=id 
\label{pora}
\end{equation}

{\bf 1.22.} In general, 
we show that the vector space $\Pi^{ur} _r (x_0;x_1;X_{K})$ possesses a canonical (``weight'') filtration
      $$W_i\Pi^{ur} _r (x_0;x_1;X_{K}) \subset \Pi^{ur} _r (x_0;x_1;X_{K})$$
compatible with the one on $\Pi^{DR} _r (x_0;x_1;X_{K})$.

Let ${\cal V}_r\subset \Pi^{ur}_r (x_0;x_1;X_K)$ be the subspace  which consists of elements 
$v\in  \Pi^{ur}_r (x_0;x_1;X_K)$ satisfying the following properties:

i)$\phi v= v$

ii)$ N^av \in W_{-a-1} \, \Pi^{ur}_r (x_0;x_1;X_K)$ for any $0<a<r$.

We prove that $dim {\cal V}_r = 1$. Moreover, the canonical morphism:
    $$ {\cal V}_r \longrightarrow   \mathbb{Q}_p$$
is an isomorphism. 

The last assertion is derived from a variant of the Monodromy Conjecture [Illusie 2] for \\ 
 $\Pi ^{ur}_r (x_0;x_1;X_{K})$. In turn, the latter follows from a result of [Mokrane].

Hence, the space ${\cal V}_r$ contains a unique element $C_r (X_K;x_0;x_1)$ satisfying (\ref{pora}).

{\bf 1.23.}In the case, when $X_{K}$ is proper and has a smooth model over $spec \, R $, 
a similar construction has been independently found
 by Besser\footnote{Although unpublished, the text is available electronically, cf. our list of references}.

The weight filtration on $\Pi^{ur}_r (x_0;x_1;X_K)$ has been independently (and using  a different method)
constructed by Kim and Hain ( [Kim]               
 
We would like to stress that Theorem A, at least in the case when $X_{K_0}$ is proper and $x_0 = x_1$,
 is equivalent to  the  Crystalline Conjecture for the fundamental group invented by Shiho.
 He also  announced  a proof of this conjecture (see  [Shiho]).

We were informed by Beilinson that he was aware that Theorem 1 is true.

{\bf Acknowledgments.}
I owe an enormous  debt of gratitude...

... first and foremost, to my advisor, David Kazhdan, for his guidance and support.

... to Sasha Beilinson.  All I know about $p$-adic Hodge Theory I have learned  from him.

... to Joseph Bernstein, Sasha Goncharov, Johan de Jong, Maxim Kontsevich and Barry Mazur for their support and many interesting discussions.

... to my dear friends - Dima Arinkin, Sasha Braverman, Misha Finkelberg, Dennis Gaitsgory and Dima Tamarkin. Their explanation helped me in the moments of despair.  
 
This paper is an extended version of author's Ph.D. Thesis.
It was finished when the author was staying at the Institut des Hautes Etudes Scientifiques in summer 2001.

{\sf Dear Reader}: may we suggest that you skip the section 2: it contains definitions 
and results, which are assumed to be known, and proceed directly to section 3.

\section{Logarithmic crystals}
In this section we remind some aspects of the theory of logarithmic crystals developed by 
Kato [Kato]. It will be used throughout the text.

{\bf 2.1. Logaritmic structure}.  
A {\it prelogarithmic structure } on a scheme $X$ is a pair $(M_X; \alpha _X)$, where 
$M_X$ is a sheaf (in the etale topology) of commutative monoids on $X$ and $\alpha _X$ is a morphism from $(M_X)$
into the multiplicative monoid ${\cal O} _X$ .

A prelogarithmic structure is called {\it logarithmic} if  the restriction
     $$\alpha _{|\alpha ^{-1}({\cal O}^*_X)}: \alpha^{-1}({\cal O}^*_X) \to {\cal O}^*_X$$
is an isomorphism.

A morphism between logarithmic schemes $(X,M_X; \alpha _X)$ and $(Y,M_Y; \alpha _Y)$
is a morphism of underlying schemes $f: X\to Y $ together with a morphism 
$f^*_M: f^{-1}M_Y \to M_X $ such that the obvious  square commutes.

It is easy to see that any prelogarithmic structure maps to a corresponding universal logarithmic structure,
called associated. If  $f: X\to Y $ is morphism of schemes and $Y$ is equipped with a logarithmic structure
$(M_Y; \alpha _Y)$ we define a logarithmic structure $f^*(M_Y)$ on $X$ to be the one associated associated to 
$f^{-1}M_Y$.

A logarithmic structure $(M_X; \alpha _X)$ is called fine if locally for etale topology it is 
associated to a prelogarithmic structure of the form $(P_X; \beta _X)$,
where $P_X$ is a constant sheaf of monoids of value $P$, with $P$ finitely generated and integral.

A morphism $i: U\hookrightarrow T$  of fine logarithmic schemes is called {\it a thickening}
 (resp. {\it a nilpotent thickening }) if 
the underlying map of schemes is a closed immersion defined by a nilideal 
(resp. a nilpotent ideal ) and the map $ i^*M_T \to M_U $
is an isomorphism.

Logaritmic smootheness is defined using Grothendieck's infinitisimal lifting property for logarithmic thickening.
 
{\bf 2.2.} What follows is a list of basic examples which will appear in this text.

a) {\bf The trivial logarithmic structure} - $M_X={\cal O}^*_X$

b){\bf  Normal crossings divisors}.   
Let $X$ be a regular scheme and $D\subset  X$ be a normal crossings divisor.  relative  to $specS$.
 The latter  means that etale locally $D$ is given by the equation $t_1t_2\cdot \cdot \cdot t_i=0$ 
where $t_1,t_2, \cdot \cdot \cdot , t_i$   is a part of regular system of local parameters at a point. 
We denote by $U = X\backslash D \hookrightarrow X$ the obvious inclusion and define 
    $$M_X={\cal O}_X \cap j_*{\cal O}^*_U \to {\cal O}_X $$
If in addition we are given a smooth morphism $X \to S$ and  $D\subset  X$ is a normal crossings divisor relative  to $specS$,
the morphism $X \to S$ of log schemes, where $S$ is equipped with the trivial logarithmic structure is smooth.

c){\bf Semi-stable reduction}. Let $S$ stand for the spectrum of a discrete valuation ring, $s$ be its closed point and
$X/S$ be a scheme with semi-stable reduction. The special fiber $X_s$ is a divisor with normal crossing, hence
by virtue of b) it defines a logarithmic structure on $X$. If we introduce a log stucture on $S$ applying b)  to  
to the pair $s \hookrightarrow S $, the  morphism of log schemes $X\to S $ is smooth.
   
d){\bf Logarithmic points}. Let $X$ be the spectrum of a field $k$.
Define $M_X = \mathbb{N} \times k^*  \to k  $ such that $(0,a) \mapsto a$ and $(n,a) \mapsto 0$ if 
$n>0$. The canonical morphism from the logarithmic point to the point with the trivial log structure is not smooth.

{\bf 2.3.} Let $f: X\to Y $ be a morphism of log schemes and  $E$ be a sheaf of ${\cal O}_X$-modules. A logarithmic derivation
with values in $E$ is a pair $(D, \delta log)$, where $D: {\cal O}_X \to E$ is a  derivation relative  to
$Y$ and $\delta log: M_X \to E$ is a monoid homomorphism such that $D(\alpha _X (m)) = \alpha (m) \delta log(m)$ and
$\delta log f^*(m) = 0$ for any section $m$ of $M_Y$. 

There is a universal derivation - "Kahler differentials" $\Omega^1_{X/Y}(\d log\infty )$. Define 
$\Omega ^q_{X/Y}(\d log\infty )= \bigwedge _{{\cal O}_X}^q \Omega^1_{X/Y}(\d log\infty )$. Then with the natural map
$d:\Omega ^q_{X/Y}(\d log\infty ) \to \Omega ^{q+1}_{X/Y}(\d log\infty )$  $(\Omega ^*_{X/Y}(\d log\infty ); d)$ becomes a 
complex.

Note that in the setting of our example b) the sheaf  $\Omega^1_{X/D}(\d log\infty )$ coincides with the sheaf of  
 differential 1-forms with logarithmic singularities at $D$. 

{\bf 2.4. Log cristalline site}. Let   $X \to S $ be a morphism of fine log schemes,
 $p^N{\cal O}_S = 0$ for some integer $N$,
$I$ be a sheaf of ideals on $S$ endowed with a PD-structure $\gamma $. We assume that  $\gamma $ extends to $X$.
 The objects cristalline site  are triples $(U,T, \delta )$,
where $U$ is a scheme etale over $X$ and 
$(T, \delta )$ is a logarithmic PD-thickening over $(S, \gamma)$.
 
Similarly, the nilpotent crystalline site consists of all PD-nilpotent thickenings.
 
As usual  one can define  logarithmic crystals.

 Let $i: X \to Z $ be a exact closed immersion into a log smooth scheme over 
$S$ (i.e. a map such that the underlying map of schemes is a closed immersion and
$i^* M_Z \to M_X $ is an isomorphism  ) and $D$ be the PD-envelope of $i$.

 One can check that the category of crystals on the  crystalline site is 
equivalent to the category of $ {\cal O}_D $ modules $E$ on $D$ (for the etale topology)
 with an integrable connection
           $$ \nabla : E \to E \otimes _{ {\cal O}_Z} \Omega^1_{Z/S}(\d log\infty ) $$
having the following property:

Let $x \in X$ be a point, $\overline x $ be its separable closure of $x$ and $t_i$
$\; (1\leq i \leq r) $ be elements of $(M_Z)_{\overline x}$ such that $d \,
log(t_i)_{1\leq i \leq r}$ 
is a basis of $\omega ^1_{Z/S,\overline x}$. Then, for any section $f$ of $E$ on a neighborhood 
of $\overline x$ and any $i$ , there exist $m_1, \cdot \cdot \cdot ,m_k$ ,$ n_1, \cdot \cdot \cdot
, n_k \in \mathbb{N}$ such that 
     $$(\prod _{1\leq j \leq k} (\nabla ^{log}_{t_i} - m_j)^{n_j})(f)=0$$
Here $\nabla ^{log}_{t_i}$ is defined by : if $\nabla (f)= \sum _{1\leq i \leq r} f_i \otimes
d\, log(t_i)$, then $\nabla ^{log}_{t_i}(f)=f_i$. (It is easy to prove that if the latter 
condition holds for one choice of $(t_i)_{1\leq i \leq r}$ then it holds for any choice
of  of $(t_i)_{1\leq i \leq r}$).
  
Similarly, one can describe crystals on the nilpotent crystalline site.
Suppose the immersion $i: X\to D$ is defined by a sheaf of ideals $J_X$
on $D$. A crystal on the nilpotent site corresponds  to  a sheaf $E$
of modules over $ {\cal O}_{\widehat{D}} = \lim {\cal O}_D/J_X^{[d]}$ complete with respect
 to the topology defined by $J_X^{[d]}$, together with an integrable connection
 $$ \nabla : E \to E \otimes _{ {\cal O}_Z} \Omega^1_{Z/S}(\d log\infty ) $$
  
 Moreover, the de Rham complex  computes the crystaline cohomology of $E$.
More precisely, if we denote by $\pi: (X/S)_{cris} \to X_{et} $ the natural map of topoi, we have a canonical isomorphism
$$   R\pi(E) \simeq E_D \otimes \Omega _{Z/S}^*(\d log\infty ) $$

If the converse is not said, we work with  the nilpotent site. 

{\bf 2.5. Filtered  logarithmic crystals}.

As it was  stated above,  the  category of coherent sheaves  equipped with an integral  connection 
 on smooth proper scheme  over ${\mathbf Z}_p $ is equivalent to the category 
of coherent crystals on the special fiber. In particularly it can be constructed just from the special fiber.
On the other hand the category of coherent sheaves together with an integrable connection and a filtration
satisfying Griffiths-transversality is not one of "the crystalline nature".
For example, endomorphisms of the special fiber do no act on it. In this section we define 
(following [Fa2])a  certain modification of this category, which admits a "crystalline" description.

A filtered  log crystal on $X/S$ is a log crystal $E$ together with a decreasing sequence of subsheaves 
$F^n(E) \subset E $, such that on any  logarithmic PD-thickening $(U,T, \delta )$ $F^n(E)_T$ is a subsheaf of 
${\cal O}_T$-modules of $E_T$ with ${\cal J}_T^{[n]} \cdot F^n(E)_T \subset  F^{n+m}(E)_T $ 
(here ${\cal J}_T^{[n]}$ is the PD-power of the ideal of $U \hookrightarrow T )$ and for any morphism
of PD-thickenings  $ f:(U,T, \delta ) \to (U^{\prime},T^{\prime}, \delta  ^{\prime}) $ the subsheaf 
$F^n(E)_T$ is equal to 
$$F^n(E)_T= \sum_{k+m\geq n } {\cal J}_T^{[m]} \cdot f^*_T (F^k(E)_{T^{\prime}}) $$     

{\bf 2.6. Example}. Assume that $X$ is log smooth over $S$.   Then any log crystal $E$ (as it was explaned above the latter 
datum is equivalent to giving a ${\cal O}_X$-module $E_X $ with an integrable
 log connection  together with a filtration  $F^n(E_X) \subset E_X $ by ${\cal O}_X$-submodules 
 satisfying Griffiths-transversality: 
$$\nabla (F^n(E_X)) \subset F^{n-1}(E_X)\otimes \Omega _{X/S}^1(\d log\infty ) $$
 gives rise to a filtered logarithmic crystal.

As usual the theory of log crystals can be extended to the case of a $p$-adic formal base $S$.

\section{Log $F$-crystals on the formal disk}
The results of this section are mostly known ([Co], [Fa2], [Ogus]) or at least dwell on well known ideas.
Nevertheless, we could not find them in the literature in a form that we can make use of them directly. 

{\bf 3.1. Notations:} $K$ is a finite extension of $\mathbf {Q}_p$, with ring of integers $R$ and 
residue field $k$, $W(k) \subset R $ is the ring of Witt vectors, $K_0$ is its field of fractions.

We start with some auxiliary  results on logarithmic connections.
 
Let   $\overline S $ be a  formal scheme  over $Spf \, R$,  isomorphic to
$$ Spf \, R [[x_1,x_2,\cdot \cdot \cdot ,x_n ]]$$
  $D=\cup D_i  \hookrightarrow \overline S $
 be a  divisor with normal crossing given by the equations $x_i=0$.
We would like to stress that we do not fix coordinates on   $\overline S $.

  Denote by $s$ the point  $ Spf \, R  =  \bigcap _i D_i \hookrightarrow \overline S $ 
The divisor $D$ gives rise to a logarithmic structure on  $\overline S $ (see {\bf 2.2.}). 
We would like to stress that we do not fix coordinates on   $\overline S $.

Let  $\overline S_{an}$ be the corresponding  rigid analytic space over $K$.
For any $0< \alpha \leq 1 $  we denote by $\overline S_{an}(\alpha ) \subset \overline S_{an} $
the open disk of radius $ \alpha $. That is an open subspace of $\overline S_{an} $ whose 
$\overline K $-points  are given by inequalities $|x_i| <  \alpha \; (i=1,2 \cdot \cdot \cdot , n) $
 The subspace  $\overline S_{an}(\alpha )$ does not depend on the choice
of coordinates. Let $ S_{an}(\alpha )$ stand for the complement $\overline S_{an}(\alpha ) - 
D_{an}$.

{\bf 3.2.} Consider   a vector bundle $E$ on $\overline S_{an}(\alpha )$ with a logarithmic integrable 
connection $\nabla E \to E \otimes   \Omega^1_{\overline S}(log \infty )$.  
Here $\Omega^1_{\overline S}(log \infty )$ stands for the sheaf of differential forms with at most 
logarithmic singularities at $D_{an} \hookrightarrow \overline S_{an}(\alpha ) $.

 The connection gives rise to a family of commuting  linear operators ("residues")
 acting on the fiber  $E_{s}$ .
 To construct it, we choose local coordinates $x_i$ as above and consider the map
$(1\otimes cont(\frac{d}{dx_i})) \nabla: E \to E $. 
Here $cont(\frac{d}{dx_i})$ is the contraction
operator. 
It is easy to see that the latter map
 descends to the fiber  $E_{s}$.      
                 Moreover,  one can check that the  resulting map
 does not depend on the choice of local coordinates. We denote it by $N_i$.
                 Clearly, that $N_iN_k=N_kN_i$.

Conversly, given a vector space $E_{s}$  over $K$ and a family of linear operators
$N_i:E_{s} \to E_{s}$  we consider  the trivial vector bundle  
$\widetilde E_{s} =E_{s} \otimes  A(\overline S_{an}(\alpha )) $ 
(here  $ A(\overline S_{an}(\alpha ))$ is the ring of analytic functions 
on  $\overline S_{an}(\alpha )$) with a 
logarithmic connection given by $\nabla (e \otimes f)= 
e \otimes df + \sum_i N_i(e)\otimes f \frac{dx_i}{x_i}$. The connection is integral 
provided that     
$N_i$ commute. As before, one can check (although it is slightly less trivial), 
that  $(\widetilde E_{s}; \nabla)$ does not depend on the choice of coordinates.
 
  \begin{lemma}\label{ma}
 Let $E$ be  vector bundle  on  $\overline S_{an}(\alpha )$ with a logarithmic integrable 
connection. Assume that operators the $N_i$ are nilpotent.
 Then for  a sufficiently  small $\beta \leq \alpha $ there exists a unique isomorphism 
$$(E; \nabla) \simeq (\widetilde E_{s},\nabla)$$
of the bundles with connection on  $\overline S_{an}(\beta ) $ identical on $ E_{s}$.
\end{lemma} 
For a proof one can just repeat the well known argument over the field of complex numbers ([De]).

If we define  {\it the category of vector bundles with a logarithmic integrable 
connection on $\overline S_{an}(0)$} to be the injective limit of the corresponding
categories on $\overline S_{an}(\alpha )$,
 the previous lemma can be reformulated as follows:
\begin{lemma}
The  category of bundles with connection with nilpotent residues $N_i$
on  $\overline S_{an}(0)$ is equivalent to the category of vector spaces over $K_0$
together with an actions of commuting operators $N_i$. 
\end{lemma}
Proof is omitted.

In particularly, we have an action of $N_i$ on any object of the latter category.
 
We are going to show that in certain cases the isomorphism in Lemma \ref{ma} extends to a larger disk.

{\bf 3.3. Logarithmic extension}. 
Following Coleman [Co] 
we define a canonical  extension $A_{Log}(\overline S_{an}(\alpha ))$ of the ring of 
rigid analytic functions $A(\overline S_{an}(\alpha ))$ on $S_{an}(\alpha )$:
$$ A_{Log}(\overline S_{an}(\alpha ))=A(\overline S_{an}(\alpha ))
[\{l(f):f \in {\cal O}_{\overline S}\cap  (A(S_{an}(\alpha ))) ^*\} ]/I $$
where $I$ is the ideal 
$$\langle \{l(fg)-l(f)-l(g): \; f,\, g \in {\cal O}_{\overline S} \cap 
(A(S_{an}(\alpha )))^* \}, $$
$$\{ l(f)=\sum_{i=1}^{\infty}\frac{(1-f)^i}{i}
\; : |f(x)-1|<1 \textstyle{for \; all \;} x \in S_{an}(\alpha ) \} \rangle $$
Let $A^f_{Log}(\overline S_{an}(\alpha )) \subset A_{Log}(\overline S_{an}(\alpha ))$ be the subring generated by
$l(f)$ where $f\in {\cal O}_{\overline S}$ is a function whose divisor is supported on 
$D \subset \overline S $.

  It is easy to check (see [Co]), that the algebra $ A_{Log}(\overline S_{an}(\alpha ))$ 
is isomorphic to the polynomial 
algebra $A(\overline S_{an}(\alpha ))[l(x_i), l(p)]$ 
and that  $A^f_{Log}(\overline S_{an}(\alpha )) =A(\overline S_{an}(\alpha ))[l(x_i)]$.
 Of course these isomorphisms depend on the choice of coordinates. 

The module $ A_{Log}(\overline S_{an}(\alpha ))$ carries a canonical logarithmic connection: 
$$ d: \; A_{Log}(\overline S_{an}(\alpha )) 
\to A_{Log}(\overline S_{an}(\alpha ))\otimes _{A(\overline S_{an}(\alpha ))}
 \Omega ^1 _{\overline S_{an}(\alpha )} $$

In addition,  $ A_{Log}(\overline S_{an}(\alpha ))$ and  $ A^f_{Log}(\overline S_{an}(\alpha ))$ 
is endowed with a family of commuting parallel endomorphisms   
$$N_i: A_{Log}(\overline S_{an}(\alpha )) \to A_{Log}(\overline S_{an}(\alpha ))$$
 We define $N_i$ 
to be the unique derivation
such that 
$$N_i (A(\overline S_{an}(\alpha ))[l(p)] = 0 \;, \; N_i l(x_i)=1 \; , \; N_i l(x_j)=0  $$

{\bf 3.4. The unipotent nearby cycles functor $\Psi ^{un} $}.
The following is inspired by Beilinson's construction of the unipotent  nearby cycles functor [Be].
Given a vector bundle $E$  on $\overline S_{an}(\alpha )$ 
 together with a logarithmic integrable connection
 we define the space   of   unipotent nearby cycles to
be the space  of parallel sections $E \otimes   A^f_{Log}(\overline S_{an}(\alpha ))$:
$$\Psi ^{un}(E)= ker\, (\nabla :\; E\otimes A^f_{Log}(S_{an}(\alpha )) 
\to   E \otimes  A^f_{Log}(\overline S_{an}(\alpha ))\otimes 
 \Omega ^1 _{\overline S_{an}(\alpha )})$$
The action of $N_i$ on  logarithmic extension $A_{Log}(S_{an}(\alpha ))$ induces one on the vector 
space  $\Psi ^{un}(E)$.
\begin{lemma}\label{dima}
 $dim_{K_0}\Psi ^{un}(E) \leq rk\,E $
\end{lemma}
{\bf Proof:} Let $\overline S_{K} = Spf \, K[[x_1, \cdot \cdot \cdot , x_n]]$. Given a vector bundle $G$ with an integrable log
connection on  $\overline S_{an}(\alpha )$ we define one on $\overline S_{K}$ : 
      $${\cal E}= E\otimes _{A(\overline S_{an}(\alpha ))}K[[x_1, \cdot \cdot \cdot , x_n]]$$
We have $dim \, E^{\nabla} \leq dim \, {\cal E}^{\nabla}$. Hence, it suffices to check the statement for vector bundles on 
 $\overline S_{K}$. In this case the result is known [De]. 
 
\begin{proposition}\label{idi}
Assume that the residues $N_i: E_s \to E_s $ are nilpotent. Then the following are equivalent:

i)  $dim_{K_0}\Psi ^{un}(E)= rk\,E $

ii) A canonical morphism 
$$(\Psi ^{un}(E) 
\otimes _{K} A^f_{Log}(\overline S_{an}(\alpha )))^{N_{\cdot}=0} 
 \longrightarrow  E $$
is an isomorphism.

Moreover, for sufficiently small $\beta \leq \alpha $  the restriction 
 $E|_{\overline S_{an}(\beta )}$ satisfies the equivalent conditions above.
\end{proposition}
Proof is obvious.

{\bf 3.5. Unipotent connections.}
Denote the category of bundles together with a unipotent logarithmic connection  on
 by ${\cal D}{\cal M}^{Log}_{un}( \overline S_{an}(\alpha ))$.
 By definition, a logarithmic connection on $E$ is unipotent if it 
can be obtained by successive extensions bundles with trivial logarithmic connection.
Logarithmic connection is called trivial if the bundle is generated by parallel sections.

\begin{theorem}
a) A connection on $E$ is unipotent if and only if  $dim_{K_0}\Psi ^{un}(E)= rk\,E $.\\
b)The functor $\Psi ^{un}$ establishes an equivalence between the category \\
 $ {\cal D}{\cal M}^{Log}_{un}( \overline S_{an}(\alpha ))$  and the category of the of 
finite-dimensional  vector spaces over $K_0$ endowed with commuting nilpotent endomorphisms $N_i$
\end{theorem}   
\begin{proof}
 We show by induction on level $r$ that the functor 
$$ {\cal D}{\cal M}^{Log}_{r}( \overline S_{an}(\alpha )) \longrightarrow  {\cal D}{\cal M}^{Log}_{r}( \overline S_K)$$ 
is an equivalence of categories (we use notations introduced in the proof of (\ref{dima}).
Assume that it is  already known  for  levels  $< r$.
To proof that the latter holds for level $r$ , it is sufficient to check 
that the evident map
$$Ext^1_{{\cal D}{\cal M}^{Log}_{r}( \overline S_{an}(\alpha ))}(A_{Log}(\overline S_{an}(\alpha )), E) \to
Ext^1_{{\cal D}{\cal M}^{Log}_{r}( \overline S_K)}({\cal O}_{\overline S_K}; E \otimes {\cal O}_{\overline S_K})$$
is an isomorphism for any $E$ of level $< r$. First, it is clear that the map is surjective.
 
The $Ext^1$-groups can be computed by the logarithmic de Rham complex. Hence, the statement follows from the following simple fact:

if a section $f \in E \otimes {\cal O}_{\overline S_K}$ satisfies the property that $\nabla(f) \in 
 \Omega^1_{\overline S}(log \infty ) \otimes E$ then $f$ itself is a section of $E$.
It completes the proof.
\end{proof}

{\bf 3.6.} Let $L$ be a finite extension of $K$.
 
We define a locally analytic function $Log \, :  L^* \to L_{st}= L[l(p)]$
to  be the unique homomorphism given by the series
$$ Log\, z =\sum_{i=1}^{\infty}\frac{(1-z)^i}{i}$$
on a neighborhood of $0\in L$ and satisfying $Log \, p = l(p) $.

Given a point $x \in S_{an}(\alpha )(L)$ it gives rise to a homomorphism 
        $$ A_{Log}(\overline S_{an}(\alpha ))\longrightarrow L_{st} $$
 
The latter induces a map
$$ \Psi ^{un}(E)  \longrightarrow  E_x \otimes _L   L_{st} $$

\begin{theorem}\label{glaza}
Let $E$ be a vector bundle  with a unipotent logarithmic connection. Then there is a canonical
isomorphism:
$$ \Psi ^{un}(E) \otimes _{K}   L_{st}  \longrightarrow  E_x \otimes _L   L_{st} $$ 
\end{theorem}
\begin{proof}
Since the functor $\Psi ^{un}$  from the category 
 $ {\cal D}{\cal M}^{Log}_{un}( \overline S_{an}(\alpha ))$ is exact,
 it suffices to check the statement in the case of trivial connection, when it is a tautology.
\end{proof}   
\begin{cor}
Under the assumption of Theorem, 
there is a canonical isomorphism
 of the fibers:
$$ E_x \otimes _L   L_{st} \to E_y \otimes _L   L_{st} $$
where $x,\, y  \in   S_{an}(\alpha )(L)$
\end{cor}
The latter corollary is a local version of Theorem B.

{\bf 3.7. F-crystals}. For the rest of this section we assume that $\overline S$ is a formal disk over $spec \, W(k)$
(i.e. $\overline S \simeq Spf \, R [[x_1,x_2,\cdot \cdot \cdot ,x_n ]]$).

Choose a lifting $\widetilde{Fr}: \overline S \to \overline S $ of the 
absolute Frobenius
$Fr:\overline S \times spec \, k \to \overline S \times spec \, k $ 
{\it compatible with the logarithmic structure on $ \overline S $ }.
It is easy to see that 
the endomorphism $\widetilde{Fr}$ acts the space  $\overline S_{an}(\alpha ) $.

An $F$-crystal  on  $\overline S_{an}(\alpha ) $ is a vector bundle $E$ together with a 
logarithmic connection and an isomorphism  
 $$\phi:\widetilde{Fr}^* E \simeq E$$
compatible with the connection.

We claim that the categories of  $F$-crystals for different liftings of the Frobenius
are equivalent to one another. It follows from the fact that the connection defines a  canonical
isomorphism $\widetilde{Fr}^* E \simeq \widehat{Fr}^* E $, where $\widehat{Fr}$ is another lifting
of the  Frobenius.
 
\begin{lemma}\label{chudo}
Let $E$ be a $F$-crystal on  $\overline S_{an}(\alpha ) $.
The vector space  $E^{\nabla}$  of parallel sections carries a canonical action of $\phi $
(i.e. the action is independent of the choice of lifting of the  Frobenius).    
\end{lemma}
Proof is obvious.                      

The logarithmic extension 
 $ A_{Log}(\overline S_{an}(\alpha ))$ is an indobject in the category of $F$-crystals:
a log lifting of the Frobenius $\widetilde{Fr} ^*: A(\overline S_{an}(\alpha )) \to A(\overline S_{an}(\alpha ))$
extends to 
$$\widetilde{Fr} ^*: A_{Log}(\overline S_{an}(\alpha ))\to  A_{Log}(\overline S_{an}(\alpha )) $$
The induced structure of (ind)$F$-crystal does not depend on the choice of lifting.
Clearly, the submodule  $ A^f_{Log}(S_{an}(\alpha ))$ also inherits a stucture of (ind)$F$-crystal.
  
Let ${\cal H}_n$ be  the ring $K_0[\phi , N_i ]$  generated by  
$\phi,\, N_i \; (i=1,2, \cdot \cdot
\cdot n)$  with the following relations
\begin{equation} 
 N_i \phi = p \phi N_i  \; , \phi \cdot a = Fr(a) \cdot \phi \; , 
         N_i \cdot a = a  \cdot  N_i
\label{hekke} 
\end{equation}
where $a \in K_0 $.

We introduce  a structure of Hopf algebra ( over $\mathbb{Q}_p$) on ${\cal H}_n$: 
the  comultiplication is given  by $N_i \to 1\otimes N_i + N_i \otimes 1 $ 
and  $\phi \to \phi \otimes \phi $. 
 
Next, we let $E$ to be a $F$-crystal on   $\overline S_{an}(\alpha )$.  
By lemma (\ref{chudo})  the vector space $\Psi ^{un}(E)$ carries a canonical action of $\phi $.
It is easy to check that the relations (\ref{hekke}) are satisfied i.e. $\Psi ^{un}(E)$ is a     
${\cal H}_n$-module.

{\bf Remark}. It is instructive to compare the  ${\cal H}_n$-module  $\Psi ^{un}(E)$
with the naive one : $E_s$. Quite surprisingly {\it there is no canonical isomorphism} between
these vector spaces. In fact , the action of $\phi $ on the fiber   
$E_s$ {\it does depend } on the choice of lifting of    the Frobenius.

\begin{theorem}
 Let  $0< \alpha \leq 1$ and $E$ be a fiber bundle on  $\overline S_{an}(\alpha ) $ 
together with a logarithmic integrable connection. Assume that $E$ has a structure of 
$F$-crystal. Then  the connection is unipotent.
In particularly, we have  a canonical isomorphism
\begin{equation}
(\Psi ^{un}(E) 
\otimes _{K_0} A^f_{Log}(S_{an}(\alpha )))^{N_{\cdot}=0} 
 \simeq E 
\label{tesis}
\end{equation}
\end{theorem}  
\begin{cor}\label{poraa}
The tensor category of  $F$-crystals  on  $\overline S_{an}(\alpha ) $ is
 equivalent to the category of finite dimensional over $K_0$ ${\cal H}_n$-modules.
\end{cor} 
\begin{proof}
Choose a $F$-crystal structure on $E$.
\begin{lemma} 
(Grothendieck) Let $V$ finite dimensional  ${\cal H}_n$-module.
 Suppose $\phi$ is an invertible operator. Then $N_i$ are nilpotent.  
\end{lemma}
Proof is omitted.

It follows that the residues $N_i$ of the connection are nilpotent. Hence by proposition (\ref{idi})
 we have an isomorphism of $F$-crystals
\begin{equation}
(\Psi ^{un}(E|_{\overline S_{an}(\beta )}) 
\otimes _{K_0} A^f_{Log}(S_{an}(\beta )))^{N_{\cdot}=0} 
 \simeq E|_{\overline S_{an}(\beta )}
\label{konec}
\end{equation}
for small $\beta $. To extend this isomorphism to  the larger disk, we make use of the following
trick invented by Dwork. Pick a lifting $\widetilde{Fr}$ of the 
Frobenius and any $\gamma < \alpha $. There exists an integer $N$ such that
   $$ \widetilde{Fr}^N(\overline S_{an}(\gamma )) \subset \overline S_{an}(\beta )$$    
Then (\ref{konec}) induces an isomorphism
$$p: \; (\widetilde{Fr}^N)^*((\Psi ^{un}(E) 
\otimes _{K_0} A^f_{Log}(S_{an}(\alpha )))^{N_{\cdot}=0}) 
 \simeq (\widetilde{Fr}^N)^* E $$
We claim that $\phi ^N \circ p \circ (\phi)^{-N} $ gives the desired isomorphism (\ref{tesis}).
\end{proof}

{\bf 3.8. Log $F$-crystals on schemes}.
Let $X$ be a log scheme in characteristic $p$. A log $F$-crystal on $X$ is a crystal $E$ on the log
crystalline site $(X/spec \, \mathbb {Z} _p )_{cris}$ ($spec \, \mathbb {Z} _p $ is endowed
with the trivial log structure) together with    
 a morphism $\phi: Fr^* E \to E$, which is an isomorphism 
in the category $\textstyle{Crystals}\otimes \mathbb {Q} $.
\begin{theorem} 
Let $i: X\hookrightarrow Y $ be a nilpotent log thickining of schemes over $spec \, \mathbb {F}_p$.
Assume that $Y$ can be embedded in a smooth log scheme over $spec \, \mathbb {Z} _p $ . 
The restriction functor $i^*$ induces an equivalence  of  categories:
 $$i^*: \{ \textstyle{coherent \, F-crystals \; on \; Y} \} \otimes \mathbb {Q}\simeq
\{ \textstyle{coherent \, F-crystals \; on \; X} \} \otimes \mathbb {Q}$$
\end{theorem}
\begin{proof}
This is another application of Dwork's trick. 
Choose a closed emmersion $s: Y\hookrightarrow Z $ in a log smooth scheme over 
 $spec \, \mathbb {Z} _p $. Let $D_s$ (resp. $D_{s\cdot i}$) be the PD-completion
of the PD-envelope of $s$ (resp. $s\cdot i$ ). There is a canonical map $\widetilde {i}:
D_{s\cdot i}\to  D_s$.  
Since the question is local, we can assume that
the Frobenius extends to $\widetilde{Fr}: D_s \to D_s $. 

There exists an integer $N$ and a morphism $p: D_s \to D_{s\cdot i}$ such that
$\widetilde {i} \cdot p =\widetilde{Fr}^N \; , \; p \cdot \widetilde {i}=\widetilde{Fr}^N$.
 We claim that  $p^*$ defines the inverse functor:
  $$  \{ \textstyle{F-crystals \; on \; X} \} \otimes \mathbb {Q}\to  
 \{ \textstyle{F-crystals \; on \; Y} \} \otimes \mathbb {Q}$$
\end{proof}
{\bf Convention}.  A  $F$-crystal on an arbitrary  log  scheme $X$ is, by the definition,   a log crystal on 
$(X\times spec \, \mathbb{F}_p\, /spec \, \mathbb {Z} _p )_{cris}$. 
 
\begin{cor}
Let $X$ be a semi-stable scheme over $spec \, R$  and $i: \,  X_k := X \times _{spec\, R} spec \, k 
\hookrightarrow X$ be its special fibre. 
The restriction functor $i^*$ induces an equivalence  of  categories:
 $$i^*: \{ \textstyle{coherent \, F-crystals \; on \; X} \} \otimes \mathbb {Q}\simeq
\{ \textstyle{coherent \, F-crystals \; on \; X_k} \} \otimes \mathbb {Q}$$
\end{cor}

{\bf 3.9. Log $F$-crystals on the log point $spec \; R $}. 
What follows is borrowed from [Fa2].

We endow the scheme $spec \, R $ with a log structure given by the closed point.

As an application of the above results we  describe  
the category of log $F$-crystals on $spec \, R \; /spec \,\mathbb {Z} _p $.
  
Let $E$ be such a crystal. We are going to show that the vector space $K \otimes _R E_{spec \, R}$
carries certain additional structure.
 
Choose a uniformizer $\pi \in R$ and let $f(T) \in W(k)[T]$ be its minimal polynomial over $W(k)\subset R $. That is an Eisenstein polynomial of degree
$e$, where  $e$ is  ramification index  of $R/\mathbf {Z}_p$. We have 
       $$R=  W(k)[T]/f(T) $$
Denote by  $V$ the PD-completion of the ring obtained
 by adjoining to $ W(k)[T]$ divided powers $\frac{f^i}{i!}$
(or, equivalently, $\frac{T^{i\cdot e}}{i!}$) . The Frobenius extends to 
$\widetilde{Fr} : spec \, V  \to spec \, V $, $T \to T^p$.
We  endow the scheme $Spf \, V $  with the log structure $M_{Spf \, V }$ 
corresponding to the divisor $T=0$.
  There is a natural immersion
     $$i_{\pi}: Spf \, R \hookrightarrow Spf \, V $$
with $i^*_{\pi}(T)=\pi$ 
By the very definition a coherent  crystal $E$ on 
the log scheme $(spec \, R \, , i^* M_{spec \, V }) $ is a coherent module 
$E_{i_{\pi}, \,  spec \, V }$ on $spec \, V $ with a logaritmic integrable connection 
     $$E_{i_{\pi}, \, spec \, V } \to  E_{i_{\pi}, \, spec \, V } \otimes _{W(k)[T]} \Omega ^1 _{W(k)[T]}(\d log\infty ) $$
Suppose now that $E$ is a $F$-crystal: 
$$ \phi:\widetilde {Fr} ^*( E _{i_{\pi}, \, spec \, V } )\otimes \mathbf {Q} \simeq   E _{i_{\pi}, \, spec \, V }\otimes \mathbf {Q} $$ 
It gives rise to a log $F$-crystal $E _{i_{\pi}}$ on the disk $\overline S_{an}(q^{- \frac{1}{(p-1)e}})$.
(Here $q$ is the cardinality of residue field $k$. So $|\pi | = q^{- \frac{1}{e}} $.)

In turn, the latter gives rise to a ${\cal H}_1$-module  $\Psi^{un}(E _{i_{\pi}}) $.   
The above construction  leads to a functor:
$$\Psi_{\pi}^{un}: \{\textstyle{log \; F-crystals \; on \; spec \, R}\} \to
 \{\textstyle{modules \; over}\; {\cal H}_1 \}$$
A priory it depends on the choice of a uniformizer $\pi$. In fact, it does not. More precisely, we have the following result.
\begin{lemma}
For any uniformizers $\pi, \pi^{\prime} \in R$, there exists a canonical isomorphism of functors:
 $$\Psi_{\pi}^{un}\simeq  \Psi_{\pi^{\prime}}^{un}$$
\end{lemma}
\begin{proof}
Choose an isomorphism
$$   
\def\normalbaselines{\baselineskip20pt
\lineskip3pt  \lineskiplimit3pt}
\def\mapright#1{\smash{
\mathop{\longrightarrow}\limits^{#1}}}
\def\mapdown#1{\Big\downarrow\rlap
{$\vcenter{\hbox{$\scriptstyle#1$}}$}}
\matrix{Spf \, R & \stackrel{i_{\pi}}{\hookrightarrow}    &Spf \, V  \cr
 \mapdown{Id} &  & \mapdown{\theta } \cr
 Spf \, R & \stackrel{i_{\pi^{\prime}}}{\hookrightarrow}    &Spf \, V  \cr}
$$
By the definition, for any crystal $E$,  we have a canonical isomorphism $r_{\theta}: \theta^* E _{i_{\pi^{\prime}}, \, spec \, V } \simeq E _{i_{\pi}, \, spec \, V }$
identical on $E_{spec \, R}$. 

On the other hand, there is a canonical isomorphism of functors: 
$$\Psi^{un} \theta^* \simeq \Psi^{un} \; : 
\{F- \textstyle{crystals \; on}\; \overline S_{an}(q^{- \frac{1}{(p-1)e}})\} \longrightarrow  \{\textstyle{modules \; over}\; {\cal H}_1 \}$$
It gives rise to $\theta ^*:\Psi_{\pi}^{un}\simeq  \Psi_{\pi^{\prime}}^{un}$. Moreover , the latter commutes with the isomorphism
$\Psi_{\pi}^{un}(E)\otimes K_{st}\simeq E_{spec \, R}\otimes K_{st}$ (see Theorem \ref{glaza}). Hence, it does not depend on the choice of $\theta$.
 It completes the proof. 
\end{proof}

Moreover, Corollary (\ref{poraa}) implies the following result.  
\begin{theorem}\label{psi} 
There is an equivalence of tensor categories: 
$$\Psi^{un}:\; \{\textstyle{coherent \; log\; F-crystals \; on \; spec \, R}\}
 \otimes  \mathbf {Q}$$ 
$$\simeq \{\textstyle{finite-dimensional \; over}
\; K_0\; {\cal H}_1-\textstyle{modules}\} $$
\end{theorem}    
 
By the construction  there is a canonical isomorphism:
\begin{equation}
\Psi^{un}(E)\otimes_{K_0} K_{st} \simeq  E_{T=\pi}\otimes_{K} K_{st}\simeq  E_{spec \, R}
\otimes_{R} K_{st} 
\label{vseuge}
\end{equation}
 
We define a $K$-linear operator  $N: K_{st} \to K_{st}$ to be the derivation of the ring $ K_{st}$ such that
$N(l(p))= e $. (Remind that $e$ stands for the ramification index $K$ over $K_0$). 

It endows both sides of (\ref{vseuge}) with an action of $N$. The action  on the left-hand side 
comes from the action on both factors,  while  $E_{spec \, R}$ is equipped with the trivial action of $N$.

\begin{lemma}\label{I}  
The isomorphism (\ref{vseuge}) commutes with the action of $N$.
\end{lemma}
Proof is a simple direct computation.

In particularly, we obtain a canonical isomorphism:
$$E_{spec \, R}\simeq (\Psi^{un}(E)\otimes_{K_0} K_{st})^{N=0}$$

{\bf 3.10.} It is convenient to have slightly more general formulation of the latter theorem.
Consider the log scheme ${\cal S}_r = (spec \, R;\, M_r \, \alpha) $, where 
$ \alpha : \, M_r \to R $
is the log structure associated to the prelogarithmic structure $\beta : \, \mathbb{N}^r \to R $, 
$ \, \beta ((0,\cdot \cdot \cdot ,1,\cdot \cdot \cdot ,0))= \pi $. Clearly, the log structure
does not depend on the choice of a uniformizer. If $r=1$, the latter  coincides
with the log structure given by the closed point $spec \, k \hookrightarrow spec \, R$  
\begin{theorem}
a)There is a tensor functor 
$$\Psi^{un}:\; \{\textstyle{unipotent \; log\; crystals \; on}\;{\cal S}_r \}
 \otimes  \mathbb {Q} \to Vect_{K_0} $$
to the category of finite-dimensional vector spaces over $K_0$ endowed with 
commuting nilpotent
 operators $N_i$, such that for any unipotent log crystal $E$ there is a canonical
isomorphism 
 $$\Psi^{un}(E)\otimes_{K_0} K_{st} \simeq  E_{spec \, R}
\otimes_{R} K_{st}$$
b)The tensor category 
 $$ \{\textstyle{ log\; F-crystals \; on}\;{\cal S}_r \}
 \otimes  \mathbb {Q}$$ 
is equivalent to the tensor category of ${\cal H}_r$-modules finite-dimensional over $K_0$. 
\end{theorem}    
  
{\bf 3.11. Hyodo-Kato isomorphism}.
 Let $  \overline X $ be a proper  scheme finite type over  $spec \, R$, 
$Z \hookrightarrow \overline X$  be a divisor, such that
etale locally $\overline X$ is isomorphic to 
\begin{equation}
spec \, R[t_1, \cdot \cdot \cdot ,t_i, s_1,
\cdot \cdot \cdot , s_j]/(t_1\cdot \cdot \cdot t_i-\pi ^m) 
\end{equation}
(here $\pi \in R$ is
a uniformizer and $m$ is an integer) , with  $Z$  given by 
the equations $\prod _{k\in J} s_k=0$ for some subset $J\subset \{1, \cdot \cdot \cdot , s_j\}.$
 
{\sf Example:}  $p: \,  \overline X \to spec \, R  $ is a proper scheme with a {\sf semi-stable reduction}, $Z=0$.

We endow  the scheme $\overline X$ with the logarithmic  structure given by the  divisor 
 $D=\overline X_k \bigcup Z $ i. e.
 $$M_{\overline X}={\cal O}_{\overline X}\cap {\cal O}^*_{\overline X -D}$$
and  $spec \, R$ with the log structure corresponding to the closed point.
Put $X_K: = \overline X - D$ 
It is easy to see that  the canonical morphism $p: \,  \overline X \to spec \, R  $
 is log smooth.

One can prove  (see [Hyodo]) that the direct image $R^*p_* {\cal O}_{cris}$
is a coherent crystal on  $  spec\, R $.  Moreover, the Base Change Theorem implies that there is 
acanonical isomorphism 
       $$(R^*p_* {\cal O}_{cris})_{ spec\, R} \otimes K \simeq H^*_{DR}(X_K)$$

Hence, by the previous result, it gives rise to a  ${\cal H}_1$-module 
$$ H^*_{st} (\overline X) = \Psi^{un}(R^*p_* {\cal O}_{cris})$$ 
and a canonical isomorphism:
$$  H^*_{DR}(X_K) = (H^*_{st} (\overline X) \otimes _{K_0} K_{st})^{N=0} $$

This construction is due to Hyodo and Kato {\it loc.cit.}  

We denote the category of  coherent log crystals on a scheme $X$ by ${\cal C}(X)$.
 
One can make use of the Leray spectral sequence and the fact that the category ${\cal C}(spec \, R) \otimes
\mathbb{Q} $ has homological dimension $1$ to prove the following result.
\begin{lemma}\label{obman}
 For any coherent log crystals $E$ and $G$ the  sequence
$$0\longrightarrow Ext^1_{{\cal C}(spec \, R) \otimes
\mathbb{Q} }
({\cal O}_{cris};p_*( E^*\otimes G))\longrightarrow 
 Ext^1_{{\cal C}(\overline X) \otimes
\mathbb{Q} }(E;G) \longrightarrow
Hom ({\cal O}_{cris};R^1p_*( E^*\otimes G))\longrightarrow 0 $$
is exact.
\end{lemma}

\section{Construction of the canonical parallel translation $C_{x_0;x_1;X_{K_0}}$ in the case when $X_{K_0}$ has a good compactification}
{\bf 4.1.} Let  $X_{K_0}$ be a smooth geometrically connected  scheme finite type over $K_0$. 
Remind  that a good compactification of $X_{K_0}$ is a smooth proper scheme $\overline X$ over $spec\, W(k)$ with 
  a divisor $Z\subset \overline X$  with  normal crossings relative to 
  $spec \, W(k)$ and an isomorphism   
$X_{K_0} \simeq (\overline X - Z)\times_{spec \, W(k)} spec \, K_0 $. 

In this section we proof a special case of Theorem B when  $X_{K_0}$ has such a  compactification.
Athough it will not be used in the general construction, the proof of this special case is significantly simpler and  
a the construction has a clear geometrical interpretation. 
 
{\bf 4.2.} Given an embedding $X_{K_0} \hookrightarrow  \overline X $ as above, we can interpret 
objects of the category  ${\cal D}{\cal M}_{un} (X_{K_0})$ as  
 unipotent  crystals on the log scheme $\overline X_k= \overline X \times spec\, k $ .
More precisely we have the following result.
\begin{lemma}
There is an equivalence of the categories:
$${\cal D}{\cal M}_{un} (X_{K_0}) \simeq \{\textstyle{unipotent\; log \; crystals\; on}\;   
 \overline X_k \}\otimes \mathbb{Q}$$
\end{lemma}
\begin{proof} By {\bf 2.4.} we have a  fully faithful functor 
$$\{\textstyle{unipotent\; log \; crystals\; on}\; \overline X_k \}\otimes \mathbb{Q} 
\longrightarrow {\cal D}{\cal M}_{un} (X_{K_0}) $$
It remains to show that its image contains all bundles with a unipotent connection. 
The latter follows from the fact that for any unipotent log crystal $E$ on $\overline X_k$
we have
$$Ext^1({\cal O}_{\overline X_k} ;E) \simeq H^1_{cris}(\overline X_k; E)\otimes \mathbb{Q}_p 
 \simeq H^1_{DR}(X_{K_0}; E_{X_{K_0}}) 
\simeq Ext^1_{{\cal D}{\cal M}_{un} (X_{K_0})}({\cal O}_{X_{K_0}}; E)$$
Here the first $Ext$ group is computed in the category of unipotent crystals $\otimes \mathbb{Q}$.
In turn, the existence of the latter isomorphisms follows form the well known result saying
that de Rham cohomology of a 
smooth scheme $X_{K_0}$ over a field of characteristic $0$ can be computed   
by virtue of the logarithmic de Rham complex on a good compactification $\overline X_{K_0}$.
It completes  the proof.
\end{proof}
{\bf Remark}. The Lemma implies that $\{\textstyle{unipotent\; log \; crystals\; on}\;   
 \overline X_k \}\otimes \mathbb{Q}$ {\it is an abelian category}. In fact, one can easily check 
that the category of unipotent log crystals on   
 $\overline X_k $ itself is abelian. The latter property might be  false for the category of
 unipotent crystals on a nonproper scheme.  
  
As corollary, we can see  that ${\cal D}{\cal M}_{un} (X_{K_0})$ possesses
 an additional symmetry:
\begin{equation}
     Fr^*: {\cal D}{\cal M}_{un} (X_{K_0}) \longrightarrow {\cal D}{\cal M}_{un} (X_{K_0})
\label{frob}
\end{equation}
induced by the Frobenius $Fr: \; \overline X_k \to \overline X_k$.
 
{\bf 4.3.} For  a point $y$  given by a morphism $i: spec\,  k \to \overline X_k$ we define a fiber finctor
              $${\cal F}_y: {\cal D}{\cal M}_{un} (X_{K_0})\longrightarrow Vect_{K_0} $$
If the point does not lie  on the divisor  $Z_k \hookrightarrow X_k$ the vector space
  ${\cal F}_y E$ is  just the stalk 
$ K_0 \otimes _{W(k)} i_{cris}^*E $.
 In general, we make use of the theory developed in the previous section and define
$${\cal F}_y E = \Psi^{un} (i^*_{cris}E)$$
where $spec \, k $ is now {\it the log point} ({\bf 2.2}). 
\begin{theorem}\label{proper}
For any pair $y_0;\, y_1 \in \overline X_k(spec\,  k)$
there exists a unique element 
$$ C_{y_0;y_1;\overline X_k} \in Mor^{\otimes}({\cal F}_{y_0}; {\cal F}_{y_1} )$$ 
satisfying the following the following property:
    $$Fr_*(C_{y_0;y_1;\overline X_k})=C_{Fr(y_0);Fr(y_1);\overline X_k} $$
\end{theorem}
\begin{proof}
For an integer $r \geq 1 $ we consider the subcategory   
${\cal D}{\cal M}_{r} (X_{K_0})\subset  {\cal D}{\cal M}_{un} (X_{K_0}) $ of 
unipotent crystals of length $r$. We denote by $\Pi _r (y_0;y_1)$ the space of morphisms
between the fiber functors ${\cal F}_{y_0}$ and ${\cal F}_{y_1}$ restricted to the latter 
subcategory.
There are canonical surjections : $$\Pi _{r+1} (y_0;y_1)\to \Pi _r (y_0;y_1)$$
and  
 $$ Mor({\cal F}_{y_0}; {\cal F}_{y_1} )= \lim_{\leftarrow}\Pi _r (y_0;y_1)$$
Note that $\Pi _1(y_0;y_1)=K_0$.

Choose an integer $n$ such that $Fr^n$ acts trivially on $k$.
The Frobenius (\ref{frob}) induces a $K_0$-linear endomorphism: 
 $$ Fr^n :\;  \Pi _r (y_0;y_1) \longrightarrow  \Pi _r (y_0;y_1) $$                
 
Theorem follows from the following Lemma.

{\bf Key Lemma}. The map $ \Pi _r (y_0;y_1) \to K_0$ induces an isomorphism:
$$(\Pi _r (y_0;y_1))^{Fr^n} \simeq K_0 $$ 

{\sf Proof of Key Lemma}. Define a filtration on $\Pi _r (y_0;y_1)$ to be
 $$  F^i=ker(\Pi _r (y_0;y_1) \to \Pi _i (y_0;y_1)) $$
where $1 \leq i \leq r $. Similarly, we can define a filtration on the algebra 
$\Pi _r (y_0;y_0)$.

Next, the  vector space $\Pi _r (y_0;y_1)$ carries a structure of
 filtered module over the filtered algebra $\Pi _r (y_0;y_0)$. Moreover we have a canonical
isomorphism:
$$ Gr^*_F(\Pi _r (y_0;y_0)) \simeq  Gr^*_F(\Pi _r (y_0;y_1)) $$
The latter is induced by the action on $1\in K_0=\Pi _1(y_0;y_1)\subset 
Gr^*_F(\Pi _r (y_0;y_1)) $.
 
Hence, it suffices to show that  $(Gr^i_F(\Pi _r (y_0;y_0))^{ Fr^n}=0$ for $i>1$.

It is easy to see that $Gr^2_F(\Pi _r (y_0;y_0))=
 (H^1_{cris}(\overline X_k)\otimes _{W(k)}K_0)^* $
It is known that the eigenvalues of $Fr^n$ acting on the cohomology group are algebraic numbers
and 
       $$|v(\alpha)| \geq p^{\frac{n}{2}}$$          
for any eigenvalue $\alpha $  and an embedding $v:\, \overline{\mathbb{Q}}\hookrightarrow  
\mathbb{C}$.
On the other hand the algebra 
$ Gr^*_F(\Pi _r (y_0;y_0))$ is generated by $Gr^2_F(\Pi _r (y_0;y_0))$.
It implies that for any eigenvalue $\alpha $ of $Fr^n$ on   $Gr^i_F(\Pi _r (y_0;y_0))$
  we have
$$ |v(\alpha)| \leq p^{\frac{(1-i)n}{2}}$$ 
It completes the proof of the Key Lemma along with the Theorem. 
\end{proof}
{\bf 4.4.} In what follows we let $x_0, x_1$ be points of $X_{K_0}(K)$.
Such a point gives rise to a 
 morphism $spec\, R \to \overline X$ and we denote by
$y_i $ ($ \; i=0;\,1$)   the corresponding points of the special fiber $\overline X_k $. 
Let  $F_{x_i}$  be the usual fiber functors:
${\cal D}{\cal M}_{un} (X_{K_0}) \to Vect_K $.

By Theorem \ref{psi} from the previous section we have a canonical isomorphism of functors
\begin{equation}
    F_{x_i} \otimes _K K_{st} \simeq {\cal F}_{y_i} \otimes _{K_0} K_{st}
\label{podl} 
\end{equation}
Combining the latter with Theorem \ref{proper} we arrive to a canonical isomorphism
$$ C_{x_0;x_1; X_K} : \, F_{x_0} \otimes _K K_{st} \simeq F_{x_1} \otimes _K K_{st}$$
whose existence is proclaimed in Theorem B.

By the construction $ C_{x_0;x_1; X_K}$ is compatible with the tensor structure.

{\bf 4.5. Remark}. A priory an element  
$ C_{x_0;x_1; X_K} \in  Mor^{\otimes}(F_{x_0} \otimes K_{st};F_{x_1} \otimes K_{st})$
depends on the choice of a good model $\overline X $. It will be proven later that in fact 
$ C_{x_0;x_1; X_K}$ does not depend on this auxiliary choice  and, moreover, the construction 
is functorial with respect to any morphism $X_{K_0}\to X^{\prime}_{K_0}$. 

More precisely, the isomorphism (\ref{podl}) gives rise to
$$\Pi _r (y_0;y_1)\otimes _{K_0} K_{st} \simeq \Pi^{DR} _r (x_0;x_1) \otimes _K K_{st}$$
where $\Pi^{DR} _r (x_0;x_1) $ stands for the space of morphisms between the functors 
restricted to ${\cal D}{\cal M}_{r} (X_{K_0})$. 

Remind that for any log crystal $E$ on $\overline X_k$ the vector space ${\cal F}_{y_i}$ is endowed
with a canonical nilpotent operator $N$. It gives rise to an endomorphism $N$ 
of $\Pi _r (y_0;y_1)$.
  
In the other words
$\Pi^{DR} _r (x_0;x_1) $ possesses a structure of log $F$-crystal
($\otimes \mathbb{Q}$)  on the logarithmic scheme $spec \, R$ 
(with the log structure given by the closed point).
  It turns out that the latter structure 
depends on $X_{K_0}$ only and, moreover, it exists for  any smooth scheme $X_K$ which has a semi-stable model over $spec \, R$.

A proof of the above statements occupies the next two sections.

\section{Construction of the fundamental crystal $\Pi_{r} $}
{\bf 5.1. Motivation: the fundamental $D$-module}. 
Let $f: \; X \to S $ be a smooth proper morphism of smooth schemes
over a field $K$ of characteristic $0$.
 For a pair of points $x_0; x_1 \in X$ with $f(x_0)=f(x_1)$
we denote by $\Pi^{DR} _r (X/S;x_0;x_1)$ the space of morphisms between the corresponding fiber
functors from the category  ${\cal D}{\cal M}_{r} (X_{f(x_0)})$
of bundles together with a unipotent connection {\it on the fiber} $X_{f(x_0)}$.
It is easy to see that the vector spaces $\Pi^{DR} _r (X/S;x_0;x_1)$ form a vector bundle 
 $\Pi^{DR} _r(X/S)$ on $X\otimes _S X$. It  immediately follows form the construction that
 $\Pi^{DR} _r(X/S)$ possesses a natural integrable connection along the fibers of the
 map $X\otimes _S X \to S $.

We claim that, in fact, $\Pi^{DR} _r(X/S)$   carries a canonical {\it total connection} on 
 $X\otimes _S X$.  
The latter is an analog of the Gauss-Manin connection.

What follows is a variant of the above construction when 
$S$ is replaced by the log scheme $spec \, R$.
  
{\bf 5.2. Notations}. Let $\overline X$ be a proper  geometrically connected scheme finite type over  $spec \, R$, 
$Z \hookrightarrow \overline X$  be a divisor, such that

etale locally $\overline X$ is isomorphic to 
\begin{equation}
spec \, R[t_1, \cdot \cdot \cdot ,t_i, s_1,
\cdot \cdot \cdot , s_j]/(t_1\cdot \cdot \cdot t_i-\pi ^m) 
\label{C}
\end{equation}
(here $\pi \in R$ is
a uniformizer and $m$ is an integer) , with  $Z$  given by 
the equations $\prod _{k\in J} s_k=0$ for some subset $J\subset \{1, \cdot \cdot \cdot , s_j\}.$ 

We endow  $\overline X$ with the logarithmic  structure given by the  divisor 
 $D=\overline X_k \bigcup Z $ i. e.
 $$M_{\overline X}={\cal O}_{\overline X}\cap {\cal O}^*_{\overline X -D}$$
It is easy to see that  $\overline X$ is log smooth  over the log scheme $spec \, R$.

We would like to stress, 
that the class of schemes satisfying the property (\ref{C}) is stable under a finite extension of the base field.   

We remind a few general concepts which will be used in the construction of fundamental crystal.
 
{\bf 5.3. $\mathbb{Q}_p$-isocrystals}. 
Let $X$ be a fine log scheme over $spec \, k$. We consider the   $p$-formal crystalline 
site which consists of 
triples $(U,T,\delta )$, where $U$ is a scheme etale over $X$, $U\hookrightarrow T$ is
a {\it p-formal } logarithmic PD-thickining, such that 
           $$T_n = T\times _{Spf W(k)} spec\, W_n(k) $$ 
is flat over $spec \, W_n(k)$  and  $U\hookrightarrow T_n$ is  PD-nilpotent.
A $\mathbb{Q}_p$-isocrystal on $X$ associates to any object  $(U,T,\delta )$  
a sheaf $E_T$ of ${\cal O}_T \otimes \mathbb{Q}$-modules on $T$,
and for any morphism $g:(U,T,\delta ) \to  (U^{\prime},T^{\prime},\delta  ^{\prime})$
an isomorphism $g^*E_{T^{\prime}} \simeq E_T$ with the evident compatibility conditions. 
   
{\bf 5.4. Category of unipotent crystals.}
Given a $p$-formal logarithmic PD-thickening $U\hookrightarrow T$ of fine log smooth scheme $U$ and a  smooth proper
log  scheme $ Y \to U$, we
let ${\cal C}_r(Y/T)$  stand for  the category  of unipotent log crystals of 
level $r$ on $(Y/T)_{cris}$ . By definition, an object $E$ of the latter category is a crystal,
which possesses  a filtration 
$$ 0=E_0 \subset E_1 \subset E_2 \subset \cdot \cdot \cdot \subset E_r=E $$  
where  $E_i$ are subobjects, such that  each successive quotient  $E_i/E_{i-1}$ is isomorhpic to
the pullback of a coherent crystal on $(U/T)_{cris}$ ($=$ a coherent sheaf  on $T$).

It is easy to see, that ${\cal C}_r(Y/T) \otimes \mathbb{Q}$ is an abelian category.
(In fact,  ${\cal C}_r(Y/T)$  itself is abelian.)

Given section  $s: U\to Y $, we consider a functor 
      $$s^*: \; {\cal C}_r(Y/T) \otimes \mathbb{Q} \to  {\cal C}_r(U/T) \otimes \mathbb{Q}
 \simeq \{\textstyle{coherent\; sheaves \;
on} \; T \} \otimes \mathbb{Q} $$
\begin{lemma}\label{nevernaya}
i)The functor $s^*$ is exact, faithful and representable by an object $ B_r(Y/T,s)$ of
${\cal C}_r(Y/T) \otimes \mathbb{Q}$.

ii)Let $h: \, (U^{\prime}\hookrightarrow T^{\prime})
 \to (U\hookrightarrow T) $ be a morphism 
of  PD-thickenings, such that $U^{\prime}= U \times _T T^{\prime}$  . 
Then the canonical morphism
\begin{equation}
     B_r(Y\times _U U^{\prime} /T^{\prime}; h^*(s)) \longrightarrow
 \widetilde{h}^* B_r(Y/T,s)
\label{base}
\end{equation}
Here $ h^*(s): U^{\prime} \to Y\times _U U^{\prime}$ is the pullback of $s$ and
$\widetilde{h}$ stands for the obvious $Y\times _U U^{\prime} /T^{\prime} \to Y$. 
\end{lemma} 
{\bf Proof:} induction on $r$.  Assume that Lemma is proven for all $r^{\prime} < r$. It follows that log crystalline cohomology 
$I_r = H^1(Y/T;  B_{r-1}(Y/T,s)) $ is a locally free ${\cal O}_T \otimes \mathbb{Q}$-module on $T$ .
We  can define  $ B_r(Y/T,s)$ to be the canonical extension of $B_{r-1}(Y/T,s)$ by ${\cal H}om(I_r;{\cal O}_T \otimes \mathbb{Q})$.
Finally, the base change property follows from the corresponding property of the crystalline cohomology [Hyodo].

{\bf 5.5. Definition of $\Pi_r(\overline X_{R/p})$}.
The latter is a $\mathbb{Q}_p$-isocrystal on 
$(\overline X_{R/p}\times _{spec \, R/p}\overline X_{R/p}/ Spf \, \mathbb{Z}_p)_{cris}$. Here the
formal scheme $Spf \, \mathbb{Z}_p$ is endowed with {\it the trivial} log structure.
 
Let  $U$ be an etale scheme over  $\overline X_{R/p}\times _{spec \, R/p}\overline X_{R/p}$
  and  $U\hookrightarrow T$ be $p$-formal logarithmic PD-thickining. The projections 
$\overline X_{R/p}\times _{spec \, R/p}\overline X_{R/p}\stackrel{\to}{\to}\overline X_{R/p}$
define canonical sections $p_i:U \to \overline X_{R/p}\times _{spec \, R/p} U \;$, $(i=0,\, 1)$.
Define 
$$\Pi_r(\overline X_{R/p})(T)= p_2^* B_r(\overline X_{R/p}\times _{spec \, R/p} U/T,p_1)$$
Lemma (\ref{nevernaya}) implies that $\Pi_r(\overline X_{R/p})$ is  $\mathbb{Q}_p$-isocrystal.

{\bf 5.5. Another construction of $\Pi_r(\overline X_{R/p})$}. 
We start with a general categorical construction.

Let ${\cal H}$ and ${\cal G}$ be sheaves of additive categories on a site ${\cal A}$ and
${\cal F}_i: \,{\cal H}\to {\cal G} \; (i=0;\, 1)$ be functors.
We define a sheaf ${\cal M}or({\cal F}_0; {\cal F}_1)$  of abelian groups on  ${\cal A}$ 
in the following way. Given an object $U$ of  ${\cal A}$ we consider the category ${\cal A}_U$ 
of objects over $U$. Denote by $|_U $ the tautological restriction functor from sheaves
on ${\cal A}$ to sheaves on ${\cal A}_U$. Define
$${\cal M}or({\cal F}_0; {\cal F}_1)(U)=Mor({\cal F}_0|_U; {\cal F}_1|_U)$$
It is easy to check that  ${\cal M}or({\cal F}_0; {\cal F}_1)$ is a sheaf.

For a log scheme  $U$ over $spec \, R/p$ 
and  p-formal logarithmic PD-thickining $U\hookrightarrow T$ we consider the category
$${\cal C}_r(\overline X_{R/p}\times _{spec \, R/p} U/T)\otimes \mathbb{Q} $$
It defines a sheaf of categories ${\cal C}^X_r \otimes \mathbb{Q} $ 
on the big  $p$-formal crystalline site  $(spec \, R/p/Spf \, \mathbb{Z}_p)_{cris}$.
In particularly, it gives rise to  a sheaf of categories on 
${\cal C}^X_r \otimes \mathbb{Q} |_{\overline X_{R/p}\times _{spec \, R/p}\overline X_{R/p}}$ 
 on the $p$-formal  site
$(\overline X_{R/p}\times _{spec \, R/p}\overline X_{R/p}/Spf \, \mathbb{Z}_p )_{cris}$,

Next, we define another sheaf  $Mod_{{\cal O}\otimes \mathbb{Q}}$  on 
the big log cristalline site $(Spf\,\mathbb{Z}_p     /Spf\,\mathbb{Z}_p )_{cris}$ assigning 
to a  PD-thickening $U\hookrightarrow T$ the category of sheaves of ${\cal O}_T\otimes \mathbb{Q} $-modules.

For $i\in\{0,1\}$ we define a functor
  $${\cal F}_i:\; {\cal C}^X_r\otimes \mathbb{Q}|_{\overline X_{R/p}\times _{spec \, R/p}\overline X_{R/p}}
 \longrightarrow Mod_{{\cal O} \otimes \mathbb{Q}}|_{\overline X_{R/p}\times _{spec \, R/p}\overline X_{R/p}}$$
to be ${\cal F}_i=(p_i \times id)^*_{cris}$, where $p_i \times id :\; U \to        
\overline X_{R/p}\times _{spec \, R/p} U$ given by the projection 
$p_i:\overline X_{R/p}\times _{spec \, R/p}\overline X_{R/p}\to
\overline X_{R/p}$.

We claim that there is a canonical isomorphism:
$$\Pi_r(\overline X_{R/p}):={\cal M}or({\cal F}_0; {\cal F}_1)$$

{\bf 5.6.} As an application of the second construction, we show that  $\Pi_r(\overline X_{R/p})$
possesses a certain additional structure.

Denote by 
$$p_{ij}:\overline X_{R/p}\times _{spec \, R/p}\overline X_{R/p}
\times _{spec \, R/p}\overline X_{R/p} \longrightarrow 
 \overline X_{R/p}\times _{spec \, R/p}\overline X_{R/p}$$
the pojection given by the formula $p_{ij}(y_0,y_1,y_2)= (y_i;y_j)$.
It immediately  follows from the construction that there is a canonical morphism
\begin{equation}
p^*_{12}(\Pi_r(\overline X_{R/p}))\otimes p^*_{23}(\Pi_r(\overline X_{R/p})) \to
p^*_{13}(\Pi_r(\overline X_{R/p}))
\label{umnogenie}
\end{equation}
{\bf 5.7.} Next, we  compare $\Pi_r(\overline X_{R/p})$ with the vector bundle   
$ \Pi^{DR} _r( X/spec \, K)$ introduced at the beginning. 
 
Since the connection on $ \Pi^{DR} _r( X/spec \, K)$ is unipotent, the underlying vector bundle has a canonical (Deligne's) extension
 $\overline{ \Pi}^{DR} _r( X/spec \, K)$ to $\overline X_K$. The latter is uniquely characterized by saying that it is a vector bundle
 with unipotent log connection whose restriction to $X_K$ coincides with $ \Pi^{DR} _r( X/spec \, K)$.  

Denote by
$\widehat{\Pi }^{DR} _r( X/spec \, K)$ the corresponding sheaf of ${\cal O}_{\widehat
{\overline X}} \times  \mathbb{Q}_p$-modules  on the $p$-formal scheme $\widehat
{\overline X}$. 
\begin{lemma}\label{H}
There is a canonical horisontal isomorphism:
$$\Pi_r
(\overline X_{R/p})_{\widehat{\overline X}} \simeq \widehat{\Pi}^{DR} _r( X/spec \, K)$$ 
\end{lemma}
{\sf Proof} is similar to one of Lemma (\ref{nevernaya}): there is the evident map 
 $$\widehat{\Pi}^{DR} _r( X/spec \, K)\to \Pi_r
(\overline X_{R/p})_{\widehat{\overline X}}$$  and  using induction on $r$ one can easily prove that it is an isomorphism.
 
{\bf 5.8.} Let us give ourself  another  pair 
$(\overline X^{\prime}; Z^{\prime} \hookrightarrow \overline X^{\prime})$, 
satisfying the property (\ref{C}) 
and a morphism of the log schemes
$f:\overline X^{\prime}_{R/p} \to \overline X_{R/p}$.

The base change property (\ref{base}) implies, that for any  scheme $U^{\prime}$ etale over 
$\overline X^{\prime}_{R/p}$ and $p$-formal  PD-thickining   $U^{\prime}\hookrightarrow T^{\prime}$
we have canonical isomorphism:
$$ ((f\times id)p_2^{\prime})^*B_r(\overline X_{R/p} \times _{spec \, R/p} U^{\prime}/T^{\prime},
(f\times id)p_1^{\prime})\simeq (f^*\Pi_r(\overline X_{R/p}))(T^{\prime})$$
By the universal property of $B_r$ there is a natural map:
$$B(\overline X^{\prime}_{R/p} \times _{spec \, R/p} U^{\prime}/T^{\prime},  
p_1^{\prime}) \longrightarrow (f\times id)^*B_r
(\overline X_{R/p} \times _{spec \, R/p} U^{\prime}/T^{\prime},(f\times id)p_1^{\prime})$$
It defines a canonical morphism: 
 $$f_* : \; \Pi_r(\overline X^{\prime}_{R/p}) \longrightarrow f^*\Pi_r(\overline X_{R/p})$$
In particularly, we have the map:
 $$\phi : \; \Pi_r(\overline X_{R/p}) \longrightarrow Fr^*\Pi_r(\overline X_{R/p})$$
induced by the Frobenius.
\begin{proposition}
The morphism $\phi$  is an isomorphism.
\end{proposition}
 \begin{proof}
The crystal $\Pi_r(\overline X_{R/p})$ possesses a filtration by the kernels of the projections:
  $$\Pi_r(\overline X_{R/p}) \longrightarrow \Pi_{r^{\prime}}(\overline X_{R/p})$$
where $r^{\prime}\leq r$. The associated graded crystal $Gr_* \, \Pi_r(\overline X_{R/p})$ 
descends to the
logarithmic point $spec\, R/p$. Moreover, the morphism (\ref{umnogenie}) defines
{\it a ring structure} on   the latter crystal. It is easy to see that the ring is generated by
$Gr_1 \, \Pi_r(\overline X_{R/p})$. On the other hand, we have a canonical isomorphism:
    $$Gr_1 \, \Pi_r(\overline X_{R/p})\simeq (R^1p_* {\cal O}_{\overline X_{R/p} \, , cris} \otimes
\mathbb{Q})^*$$
Here $p$ stands for the map $\overline X_{R/p}\to spec\, R/p$. Hence, the 
Proposition follows from the
fact ([Hyodo],  Proposition 2.24) that the Frobenius
$$ \phi: \, R^ip_* {\cal O}_{\overline X_{R/p} \, , cris} \otimes
\mathbb{Q}\longrightarrow R^ip_* {\cal O}_{\overline X_{R/p} \, , cris} \otimes
\mathbb{Q}$$
 is an isomorphism.

\end{proof} 

{\bf 5.9. Remark}. One can show (see Section {\bf 9}) that $\Pi_r(\overline X_{R/p})$ is an object of
the category
${\cal C}_r(\overline X_{R/p}\times _ {spec \, R/p}\overline X_{R/p} /spec \, W(k))
 \otimes \mathbb{Q}$.

\section{Crystalline structure on  ${\cal P}^{DR}_{r}(X_{\overline K}) $}
{\bf 6.1.}
Let $X_K$ be a smooth variety over $K$. 
  Remind that  ${\cal P}^{DR}_{r}(X_{\overline K})$ is the category whose set of objects
is $X_K(\overline K)$ and whose group of morphisms between objects $x_0, x_1 \in X_K(\overline K) $ is
$\Pi^{DR} _r (x_0;x_1) $. 

We denote by $K_0^{ur}$ the maximal unramified extension of $K_0$ and by 
${\cal H}^{ur}_1$   the ring ${\cal H}_1 \otimes _{K_0} K_0^{ur} $. The latter is generated by  
$\phi,\, N \; a \in K_0^{ur}$, 
satisfying the  following relations 
              $$\phi N = p N \phi \; , \phi \cdot a = Fr(a) \cdot \phi \; , 
         N \cdot a = a  \cdot  N  $$

The main result of this section is the following theorem.
\begin{theorem}\label{E}
For any $X_K$ there are
 
1) a category  ${\cal P}_{r}(X_{K})$ whose set of objects is $X_K(\overline K)$ and  
whose groups of morphisms  $Mor _{{\cal P}_{r}(X_{ K})}(*,*)$ are endowed with
${\cal H}^{ur}_1$-module stucture
 
2) an isomorphism of categories
 \begin{equation}
 {\cal P}_{r}(X_{ K})\otimes_{K^{ur}_0}\overline K_{st} \simeq  {\cal P}^{DR}_{r}(X_{\overline K})
\otimes _{\overline K} \overline K_{st}
\label{AA}
\end{equation}
identical on objects.
        
For any morphism $f: X^{\prime} _{\overline K} \to X_{\overline K}$  there is a functor 
\begin{equation}      
f_*:{\cal P}_{r}(X_{K}) \longrightarrow {\cal P}_{r}(X^{\prime}_{ K})
\label{tak}
\end{equation}

such that the evident diagram is commutative.

Let $\Pi^{ur}_r (X_{K};x_0;x_1)$ stand for $Mor_{{\cal P}_{r}(X_{ K})}(x_0;x_1)$. 

3)The isomorphism  
 \begin{equation}
T: \; \Pi^{ur} _r (X_{ K};x_0;x_1)\otimes_{K^{ur}_0}\overline K_{st} \simeq \Pi^{DR} _r (X_{\overline K};x_0;x_1)
\otimes _{\overline K} \overline K_{st}
 \label{B}
\end{equation}
induced by (\ref{AA}) commutes with the action monodromy $N$. 

(The latter acts trivially on $\Pi^{DR} _r (X_{\overline K};x_0;x_1)$ and $N:\overline K_{st}\to \overline K_{st}$ is 
the a derivaton over $\overline K$, such that $N(l(p))=e$).  
\end{theorem}
{\bf 6.2. Remarks}. $1 \; .$ The theorem implies that $\Pi^{DR} _r (X_{\overline K};x_0;x_1)$ 
possesses a structure of log $F$-crystal on $spec \, R_L$ for some finite extension $L\supset K$. 
If the scheme $X_K$ has a semi-stable model over $spec \, R$ (or, more generally, a model satisfying
(\ref{C}) )  we can choose $L$ to be the field of definition of $x_0;x_1$. But in general,  it is not possible.

$2\; .$ The category ${\cal P}_{r}(X_{\overline K})$ depends on $X_{\overline K}$ and not on $X_K$. But the action of the 
monodromy operator
$N$ on morphisms in the latter category depends on the choice of a finite extension $K$. More precisely, there is a canonical 
isomorphism:  
$$ S_{L/K}: \,  \Pi _r (X_{K}\times spec \, L;x_0;x_1) \simeq \Pi _r (X_K;x_0;x_1) $$
compatible with (\ref{B}) and satisfying the following: 
 \begin{equation}
S_{L/K} \cdot N= e_{L/K} N \cdot S_{L/K}  \; ; \; \phi \cdot S_{L/K}  = S_{L/K} \cdot \phi
\label{HH}
\end{equation}
Here  $e_{L/K} $ is the ramification index of $L$ over $K$. 
 
{\bf 6.3. Notation}.
 The homomorphism $\overline K_{st} \to \overline K$, $l(p) \to 0$ defines a structure of $\overline K_{st}$
-module on $\overline K$. 
Tensoring (\ref{B}) with $\overline K $ over $\overline K_{st}$ we arrive to an isomorphism:
$$T_{l(p)=0}: \;   \Pi _r (X_K;x_0;x_1) \otimes _{K^{ur}_0} \overline K \simeq  \Pi^{DR} _r (X_{\overline K};x_0;x_1)$$
Note that the latter homomorphism determines a morphism $T$ of $\overline K_{st}[N]$-modules uniquely.

 The result is an almost immediate consequence  of the construction in the previous section,
de Jong's Alteration Theorem and the technique of descent. We start with recollecting some
aspects of the latter theory.

{\bf 6.4. Descent}. Consider the following diagram:
\begin{equation}
 \mathbf{X}_1 \;
{\begin{array}{l} 
p_1 \\
\longrightarrow \\
p_2 \\
\longrightarrow \\
\end{array}} 
\; \mathbf{X}_0
\stackrel{f}{\longrightarrow} X_{ K} 
\end{equation}
We assume 
that all schemes in the diagram are smooth over a field $K$ of characteristic $0$   and
the maps  $f$ and $(p_1 \times p_2): X_1 \to X_0 \times _{X_K} X_0$  are proper, surjective and generically etale.
In addition, we assume that $X_{ K}$ and $X_0$ are geometrically irreducible.

Let  ${\cal D}{\cal M}^{\theta}_{r} (X_0) $  be  the subcategory of  ${\cal D}{\cal M}_{r} (X_0)$, which consists of 
 a vector bundles $ E^{\prime}$  on $X_0$ together with an integrable unipotent connection such that
            $$p_1^* E^{\prime}\simeq p_2^* E^{\prime}$$
     \begin{proposition}\label{a}

We have an equivalence of categories:    
$$f^*: \, {\cal D}{\cal M}_{r} (X_K) \simeq {\cal D}{\cal M}^{\theta}_{r} (X_0)$$
\end{proposition}
\begin{proof}

First, we note that the map 
    $$H^1_{DR}(X_K) \longrightarrow H^1_{DR}(X_0)$$
is injective. It implies that the map 
              $$ \Pi^{DR}_r (X_0;x_0;x_0)  \longrightarrow         \Pi^{DR}_r (X_K;f(x_0);f(x_0))$$
is surjective. Hence, the functor $f^*$ is fully faithful. To complete the proof it suffices to show that 
$f^*$ induces a surjection  on $Ext^1$-groups.

Without loss of generality we may assume that $f$ is finite. Indeed, in general there exists a closed subscheme $Z \hookrightarrow X_K$ 
{\sf of codimension $\geq 2$ } such that $f$ is finite over $X_K-Z$. On the other hand, the restriction 
 $${\cal D}{\cal M}_{r} (X_K)\to  {\cal D}{\cal M}_{r} (X_K-Z)$$
is an equivalence of categories.

If, in addition, $f$ is etale, the claim immediately follows from the fact that 
for any vector bundle with a unipotent integrable connection (more generally, a connection with regular singularities)
   $$ H^i_{DR}(X_K; E) = ker(H^i_{DR}(X_0;f^*E) \stackrel{\longrightarrow}{\to} H^i_{DR}(X_0;p^*_1f^*E))$$
In general, we denote  by $D \hookrightarrow X_K$ the ramification divisor of $f$ and make use of the folowing isomorphism:
$$ H^1_{DR}(X_K; E)\simeq (H^1_{DR}(X_0;f^*E)\oplus H^1_{DR}(X_K- D ;E) \to H^1_{DR}(X_0 - f^{-1}(D);f^*E)$$

\end{proof}

{\bf 6.5. Alteration Theorem}.
 Let  $Z \hookrightarrow \overline X $ be a closed immersion  of regular schemes over $spec \, R$,
$\overline X_{k} = \bigcup _i \overline X_{k \, , i}$ and $Z=  \bigcup _j Z_j $ be the decompositions in their irreducible 
components.  

 Remind that $ (Z;\overline X)$ is  {\it a strict semi-stable pair } over $spec \, R$, if locally (for etale topology)     
 $\overline X$ is isomorphic to 
\begin{equation}
spec \, R[t_1, \cdot \cdot \cdot ,t_n, s_1,
\cdot \cdot \cdot , s_m]/(t_1\cdot \cdot \cdot t_n-\pi ) 
\label{D}
\end{equation}
 with  $Z_j$  given by 
the equations $ s_j=0$ for some $j \in \{1, \cdot \cdot \cdot , m\}$ and $\overline X_{k \, , i}\hookrightarrow \overline X_k$
 given by $t_i=0$. In particularly, the schemes $\overline X_{k \, , i}$ are smooth over  $spec \, R$. 

We are going to make use of the following fundamental result.

{\bf Theorem}.(de Jong.) Let $X$ be an integral, flat  and finite type over $spec \, R$ and   
$S \hookrightarrow  X $ be a proper closed subset. There exist a finite extension $L\supset K$ with ring of integers 
$R_L\supset R$, an integral scheme $ X_0$ over $spec \, R_L$, a proper generically etale morphism over  $spec \, R$:
        $$ f: X_0 \longrightarrow X $$
and an open immersion $j:  X_0 \hookrightarrow \overline X_0 $, with the following properties:

a)$ \overline X_0 $ is a projective scheme over $spec \, R_L$ , and

b)The pair $(\overline X_0, f^{-1}(S)_{red}\cup \overline X_0 \backslash j(X_0))$ is strict semi-stable.

{\bf 6.6. Proof of Theorem \ref{E} }. Without loss of generality we can assume, that $X_K$ is  geometrically irreducible. 

First, assume that we are given a smooth  scheme $X_K $ together with a compactification $X_K \hookrightarrow \overline X $
satisfying (\ref{C}) (with $Z_K =\overline X_K \backslash X_K $).  For a pair of points $x_i \in X_K(L)$ denote by
$i_{x_0,x_1}: spec \, L/p  \to \overline X_{R/p} \times _{spec\, R/p} \overline X_{R/p} $ the corresponding morphism.
We endow the scheme $spec \, L/p$ with the log structure associated to the prelogarithmic structure
 $\beta : \, \mathbb{N} \to L/p $, 
$ \, \beta (1)= \pi _K $ (here $\pi _K$ is a uniformizer of $K$) and  define
$$\Pi _r (X_K;x_0;x_1) :
= \Psi ^{un} (i_{x_0,x_1}^*\Pi_r(\overline X_{R/p}))$$
$$\Pi^{ur} _r (X_K;x_0;x_1) :
=  \Pi _r (X_K;x_0;x_1)  \otimes _{K_0} K_0^{ur} $$
 A construction   of the isomorphism (\ref{B})  immediately follows from the very definition of $ \Psi ^{un}$.
 In turn, Lemma \ref{I} implies that the latter 
commutes with the action of $N$. The direct image functor (\ref{tak})  is defined for any log morphism 
$X_{\overline K}^{\prime} \to X_{\overline K}$.
  
In general, we start with any open immersion $X_K \hookrightarrow \overline X$ into   
a proper, integral, flat over $spec \, R$ scheme.    
Next, using de Jong's theorem  we can construct    
a diagram   
\begin{equation}
 \mathbf{\overline X}_1 \;
{\begin{array}{l} 
\overline p_1 \\
\longrightarrow \\
\overline p_2 \\
\longrightarrow \\
\end{array}} 
\; \mathbf{\overline X}_0
\stackrel{\overline f}{\longrightarrow} \overline X \times spec\, R_L
\label{ba}
\end{equation} 
with the following properties:

1) Schemes  $\overline X_i$,  $(i=0, 1  )$ are proper over $spec \, R_L$.

2)The maps  $\overline f$ and $(p_1 \times p_2):\overline X_1 \to \overline X_0 \times _{\overline X \times spec\, R_L}\overline X_0$ 
are proper, surjective and generically etale.
 
3)  $\overline X_0$ is geometrically irreducible.

Let $X_i = \overline X_i \times _{\overline X \times spec\, R_L  } X_K\times spec\, L  \hookrightarrow \overline X_i $.

4)The pairs $(\overline X_i; X_i)$ satisfy (\ref{C}) (with $K$ replaced by $L$).
      
Let ${\cal A}_r(\overline X_i)$ be the category of $\mathbb{Q}_p$-linear functors from 
 ${\cal P}_{r}(X_{i})$ to the category $Vect_{K_0^{ur}}$ of finite-dimensional vector spaces over  $K_0^{nr}$.

Thanks to (\ref{AA}), we have a functor
  \begin{equation}
T_{l(p)=0}: \; {\cal A}_r(\overline X_i)   \longrightarrow {\cal D}{\cal M}_{r} (X_i\times spec \,\overline K )
\label{J}
\end{equation}
given by tensor product with $\overline K $.
 
Next, we define  a category  ${\cal A}_r(X_K)$ be the subcategory of ${\cal A}_r(\overline X_0)$, which consists of objects 
 $ E^{\prime}$  satisfying
            $$p_1^* E^{\prime}\simeq p_2^* E^{\prime}$$

The descent property (\ref{a}) together with (\ref{J}) imply the existence of a functor
 \begin{equation}
 {\cal A}_r(X_K)   \longrightarrow {\cal D}{\cal M}_{r} (X_{\overline K})
\label{ca}
\end{equation} 
By the very definition for any point $x \in X_0(\overline K)$ we have  a fiber functor:
$${\cal F}_x: \; {\cal A}_r(X_K)   \longrightarrow Vect_{K_0^{ur}}$$
Further,  
for a pair $\widetilde{x}_0,\widetilde{x}_1 \in X_0(\overline K)$ define   
$$ \Pi^{ur} _r (X_K;\widetilde{x}_0;\widetilde{x}_1)= Mor({\cal F}_{\widetilde{x}_0};{\cal F}_{\widetilde{x}_1})$$
The functor (\ref{ca}) induces a map
 \begin{equation}
T^{-1}_{l(p)=0}: \Pi^{DR} _r (X_{\overline K};x_0;x_1) \longrightarrow  
\Pi^{ur} _r (X_K;\widetilde{x}_0;\widetilde{x}_1)\otimes \overline K
\end{equation}
where $x_i=f(\widetilde{x}_i)$.
\begin{lemma}
The morphism $T^{-1}_{l(p)=0}$ is an isomorphism.
\end{lemma}
This is a simple exercise.

 The proposition (\ref{a}),  together with the lemma imply that for any object $E^{\prime}$ of  ${\cal A}_r(X_K)$  there is a {\it canonical} isomorphism 
           $$\theta : p_1^* E^{\prime}\simeq p_2^* E^{\prime}$$ 
As a consequence, we obtain a canonical identification between 
$\Pi^{ur} _r $ groups for different liftings of $x_i$. We define $\Pi _r (X_K;x_0;x_1)$ to be 
"the common value" of  $\Pi^{ur} _r (X_K;\widetilde{x}_0;\widetilde{x}_1)$

As we mentioned  above the latter gives rise to an isomorphism    of $\overline K_{st}[N]$-modules
(\ref{B}).
As usual, one can proof that $ \Pi^{ur} _r (X_K;x_0;x_1)$ does not depend on the choice of resolution (\ref{ba}) we made:
given two such diagrams, we can map them to a third one, and the induced maps on 
$ \Pi^{ur} _r (X_K;x_0;x_1)$ are isomorphisms
(since they are isomorphisms on $ \Pi^{DR} _r (X_{\overline K};x_0;x_1)$).

{\bf 6.7} Let $x_0;x_1$ be $K$-points of $X_K$. Then the vector space $ \Pi^{DR} _r (X_{\overline K};x_0;x_1)$ carries a
canonical $K$-structure:
  $$ \Pi^{DR} _r (X_{\overline K};x_0;x_1) \simeq \Pi^{DR} _r (X_{K};x_0;x_1) \otimes _K \overline K$$
It gives rise to a simi-linear action of $Gal(\overline K/K)$ on the former space. 

It is easy to see that 
the subspace $\Pi^{ur} _r (X_K;x_0;x_1)\subset \Pi^{DR} _r (X_{\overline K};x_0;x_1)$ is
invariant under this action. Moreover, the restriction of the induced action of   $Gal(\overline K/K)$ on 
 $\Pi^{ur} _r (X_K;x_0;x_1)$ to the inertia subgroup $I\subset Gal(\overline K/K)  $ factors through 
a finite quotient. 

\section{The Monodromy Conjecture for $\Pi^{ur}_r (X_K;x_0;x_1)$, where $X_K$ is a   proper scheme }
{\bf 7.1. Purity}. Let $V$ be an finite-dimensional ${\cal H}^{ur}_1$-module.
Choose a finite unramified extension $\mathbb{Q}_p \subset K_0$ and a  ${\cal H}_1$-stable lattice 
            $$ V_0 \subset V \; , \; V_0 \otimes _{K_0} \mathbb{Q}_p^{ur} = V $$
 Fix an integer $c$, such that $Fr^c$ acts trivially  on $K_0$.
 
The nilpotent operator $N$ gives rise to an increasing filtration $M_{\cdot} V_0\subset V$ characterazed by the property 
  $$N(M_iV_0)\subset M_{i-2}V_0 \; \textstyle{and} \; N^r: gr^M_r V_0 \simeq gr^M_{-r} V_0$$
The filtration is stable under $\phi ^c$; hence the latter acts on the associated graded space:
        $$ \phi ^c: \; gr^{M}_r V_0 \longrightarrow gr^{M}_r V_0 $$
The module $V$ is called {\it pure of weight} $i$  if any eigenvalue $\alpha $ of the latter operator is an algebraic 
number and for any embedding  $v:\overline{\mathbb{Q}} \hookrightarrow \mathbb{C}$:    
                  $$ |v(\alpha)|= p^{c\frac{r+i}{2}} $$
Clearly, that the the definition is independent of the choice of $V_0$.
           
One can easily check, that for any two  pure ${\cal H}^{ur}_1$-modules $V$ and $V^{\prime}$ of  weights
$i$ and $i^{\prime}$ respectively, we have 
\begin{equation}
Hom_{{\cal H}^{ur}_1}(V;V^{\prime})=0
\label{bei}             
\end{equation}
{\sf provided that $i> i^{\prime}$ }.

   It is conjectured (see, for example, [Illusie]) that the  log crystalline cohomology group $H^i_{st}(X) \otimes 
\mathbb{Q}_p^{ur}$ of a proper semi-stable scheme over $spec \, R$ is a pure   ${\cal H}^{ur}_1$-module of weight $i$.

{\bf 7.2. Mixed  ${\cal H}^{ur}_1$-modules.} The concept of mixed   ${\cal H}^{ur}_1$-module is an invention of 
Beilinson, Schneider and Illusie (unpuplished). 

A mixed module $V$ is a pair $(V; W_{\cdot}V)$, where $V$ is  a ${\cal H}^{ur}_1$-module and $ W_{\cdot}V$ is an
increasing filtration by  ${\cal H}^{ur}_1$-submodules ( $ W_{\cdot}V$ is called {\it the weight filtration}),
with the following condition:
  
for any integer $i$  the   ${\cal H}^{ur}_1$-module $gr^W_i V$ is pure of weight $i$.

A morphism of mixed modules is a homomorphism of underlying  ${\cal H}^{ur}_1$-modules preserving  the filtration. 

The property (\ref{bei}) implies the following result:
\begin{proposition} \label{beil}
Any morphism between mixed modules is strictly compatible with the  weight filtration. In particularly, the category of mixed modules
is abelian.
\end{proposition}
In addition, there is the evident tensor structure on the latter category.
\begin{lemma}\label{extension}
 Let $0\to V_1 \to V_2 \to V_3 \to 0 $ be an exact sequence of filtered ${\cal H}^{ur}_1$-modules  $(V_i; W_{\cdot}V_i)$
We  suppose that all morphisms are  strictly compatible with the filtration.  Assume that $V_1$ and $V_3$ are mixed. Then $V_2$ is also 
mixed. 
\end{lemma}
Proof is obvious.

{\bf 7.3. Mixed structure on $\Pi^{ur} _r (X_K;x_0;x_1)$}. Let $X_K$ be a smooth scheme over $spec \, K $, $x_0, x_1 \in 
X_K(K)$ .
 We define a filtration
on $\Pi^{ur}_r (X_K;x_0;x_1)$ to be:
 $$W_{-i}\Pi^{ur} _r (X_K;x_0;x_1): = ker(\Pi _r (X_K;x_0;x_1) \to \Pi _{i+1} (X_K;x_0;x_1))$$
for $0\leq i \leq r-1$,  and $W_{-i}\Pi _r (X_K;x_0;x_1)=0$ (resp. $=\Pi _r (X_K;x_0;x_1)$) for $i\geq r$ (resp. $i\leq 0$).
\begin{theorem}\label{M}
Suppose that  $X_K$ is proper. The
 pair $(\Pi _r (X_K;x_0;x_1);W_{\cdot })$ is a mixed   ${\cal H}^{ur}_1$-module.      
\end{theorem}  
\begin{proof}
By the very definition $\Pi _r (X_K;x_0;x_1)$ carries a structure of a  $\Pi _r (X_K;x_0;x_0)$-module. 
The $gr ^W_{\cdot} \Pi _r (X_K;x_0;x_0)$-module   $gr ^W_{\cdot} \Pi _r (X_K;x_0;x_1)$  has a canonical generator 
$\mathbf{1} \in  \Pi _1 (X_K;x_0;x_1)=K_0^{ur}$, which defines an isomorphism of  ${\cal H}^{ur}_1$-modules:
$$gr ^W_{\cdot} \Pi _r (X_K;x_0;x_0) \simeq  gr ^W_{\cdot} \Pi _r (X_K;x_0;x_1)$$
We note, that the ${\cal H}^{ur}_1$-module  $gr ^W_{\cdot} \Pi _r (X_K;x_0;x_0)$ does not depend on $x_0$. It justifies the notation 
 $gr ^W_{\cdot} \Pi _r (X_K)$ we use bellow.
   
Hence, it suffices to check, that $gr ^W_{\cdot} \Pi _{un} (X_K)$ is mixed. (We let $\Pi _{un} (X_K;x_0;x_0)$ stand for the projective limit of 
 $\Pi _r (X_K;x_0;x_0)$.)
Assume, first, that $X_K $ has a semi-stable compactification $X_K \hookrightarrow X $.
 
We make use of the following fact [DGMS]:

  $gr ^W_{\cdot} \Pi_{un} (X_K)$  is  a quadratic algebra generated by
 $gr ^W_{-1} \Pi _{un} (X_K)= (H^1_{st}(X))^*$  with the relations 
$$ Im((H^2_{st}(X))^* \to (H^1_{st}(X))^*\otimes (H^1_{st}(X))^*)$$
It is an immediate consequence of the well known result on the de Rham fundamental group.

Therefore, it is enough to check that
 
1) $H^1_{st}(X)$ is pure of weight $1$.

2) $ker(H^1_{st}(X)\otimes H^1_{st}(X) \to H^2_{st}(X))$ 

The latter follows from a result of Mokrane [Mo], 
who proved that $ H^i_{st}(X)$ are pure for a semi-stable $X$ of dimension $\leq 2$.

The  general case can be reduced to semi-stable using the following trick. First, assuming that $X_K$ is geometrically connected and using  de Jong's theorem 
  we can find  a finite
  extension $L\supset K$,
 a semi-stable geometrically connected scheme $Y$ over $spec \, R_L$  and  proper generically etale morphism $f: \, Y_L \to X_K \times spec \, L$ . 
Clearly, that the map 
$$f_* : \, gr ^W_{\cdot} \Pi _r (Y_L) \to gr ^W_{\cdot} \Pi _r (X_K \times spec \, L)$$ 
is surjective. Moreover, one can easily check that its kernel is a direct summand in a pure  ${\cal H}^{ur}_1$-module 
$gr ^W_{\cdot} \Pi _r (Y_L)$.
Hence, the image is also pure. It completes the proof. 
\end{proof}
 
{\bf 7.4. The Monodromy Conjecture for the crystalline cohomology.} In this subsection we remind the statement of the conjecture. It
is known to specialists, although we could not find it in the literature.

For  a smooth scheme $X_K$ finite type over $K$, one can construct a finite-dimensional 
 ${\cal H}^{ur}_1$-module $H^{*}_{pst}(X_K)$
together with an increasing filtration $W_iH^{*}_{pst}(X_K) \subset H^{*}_{pst}(X_K)$ by submodules and a 
 canonical isomorphism 
    $$H^{*}_{pst}(X_K) \otimes _ {\mathbb{Q}_p^{ur}} \overline K \simeq H_{DR}^*(X_K) \otimes _K \overline K $$ 
compatible with the weight filtration $W_i H_{DR}^*(X_K) \subset H_{DR}^*(X_K)$  on the de Rham cohomology.
 
In the case when $X_K$ has a compactification  
$X_K \hookrightarrow  \overline X $ as in {\bf 5.2} we define  
$$H^{*}_{pst}(X_K): = H^{*}_{st}(X) \otimes  \mathbb{Q}_p^{ur}$$

In general $H^{*}_{pst}(X_K)$ is defined using de Jong's Alteration theorem.
\begin{conjecture} 
The pair $(H^{n}_{pst}(X_K); W_{\bullet })$ is mixed  ${\cal H}^{ur}_1$-module.
\end{conjecture}
For $n \leq 2$, the conjecture follows form the already cited result of Mokrane (see Appendix).

  \section{Proof of Theorem B}
We let $X_K$ be a smooth scheme over $K$ and $x_0, x_1 \in X_K(K)$.
\begin{theorem}\label{main}
There exists a unique element 
$$C_{\cdot}(X_K,x_0;x_1) \in  \lim _{\leftarrow} \Pi ^{ur} _{r} (X_{K};x_0;x_1)$$    
characterized by the following properties:

1)$C_{1}(X_K,x_0;x_1)= \mathbf{1} \in \Pi ^{ur} _{1} (X_{K};x_0;x_1) = \mathbb{Q}_p^{ur}$

2)$C_{\cdot}(X_K,x_0;x_1)^{\phi} = C_{\cdot}(X_K,x_0;x_1)$

Choose an open immersion  $j: \, X_K \hookrightarrow \overline X_K$ into a proper smooth variety over $spec \, K$  such that the complement 
$D=\overline X_K - X_K$ is a divisor with simple normal crossings. 

3) $ N^{r-1} (j_*C_{r}(X_K,x_0;x_1))=0$
  \end{theorem}
{\bf Remark}. The element $C_{\cdot}(X_K,x_0;x_1)$ does not depend on the choice of $j$. Proof:  given two such immersion 
$j_i: \, X_K \hookrightarrow \overline X _i$
we can construct a third one $j: \, X_K \hookrightarrow \overline X_K$  and maps
$$
\def\normalbaselines{\baselineskip20pt
\lineskip3pt  \lineskiplimit3pt}
\def\mapright#1{\smash{
\mathop{\longrightarrow}\limits^{#1}}}
\def\mapdown#1{\Big\downarrow\rlap
{$\vcenter{\hbox{$\scriptstyle#1$}}$}}
\matrix{ X_K
 & \mapright{j_{i}} & \overline X _i  \cr
\mapdown{Id} &  &\mapdown{f_i } \cr
X_K & \mapright{j} &  \overline X _K \cr}
$$
The uniqueness implies that the element $C_{\cdot}(X_K,x_0;x_1)$ constructed using $ \overline X _K$ coincides with the one constructed using $\overline X _i$. 

{\bf Proof}.
Given  a finite-dimensional ${\cal H}^{ur}_1$-module $V$ and an integer $i$, we define a certain subspace $V^i \subset V$.
For this we choose a finite unramified extension $\mathbb{Q}_p \subset K_0$ and a  ${\cal H}_1$-stable lattice 
$V_0 \subset V$, $ V_0 \otimes _{K_0} \mathbb{Q}_p^{ur} = V $ , an integer $c$, such that $Fr^ca =a$ for  $a \in K_0$ and define 
$V_0^i \subset V_0 $ to be the maximal $\phi$-invariant subspace, where all eigenvalues of $\phi ^c$ are Weil numbers of weight $i$
(meaning that they are algebraic and  $ |v(\alpha)|= p^{ \frac{ci}{2}} $  for any embedding $v:\overline{\mathbb{Q}} 
\hookrightarrow \mathbb{C}$). Put $V^i=V_0^i\otimes K_0^{ur}$.

Next, we apply the above construction to $\Pi ^{ur} _{r} (X_{K};x_0;x_1)$. There is a decomposition:
 $$\Pi ^{ur} _{r} (X_{K};x_0;x_1)  = \bigoplus _{i\leq 0}\Pi ^{ur} _{r} (X_{K};x_0;x_1)^i$$
Indeed, it suffices to proof that all eigenvalues  of $\phi ^c$ on $gr ^W_{\cdot} \Pi^{un}_r (X_K)_0$ are Weil numbers of weight 
$i\leq 0$. The latter is generated by   $gr ^W_{-1} \Pi^{un}_r (X_K)_0$ for which the result is known. 
(As before, it suffices to treat 
the  semi-stable case. In the latter case we have $gr ^W_{-1} \Pi^{un}_r (X_K)= (H^1_{st}(\overline X))^*\otimes  K_0^{ur} $).     

Further, we define by induction on $r$  certain subspaces $L_r \subset \Pi ^{ur} _{r} (\overline X_{K};x_0;x_1)^0$:\\
$L_1=\Pi ^{ur} _{1} (\overline X_{K};x_0;x_1)^0 =K_0^{ur}$ and $L_{i+1}=p_{i+1}^{-1}(L_i)\cap ker\, N^i$,
 where $p_{i+1}$ is the projection
    $$  \Pi ^{ur} _{i+1} (\overline X_{K};x_0;x_1)^0 
\to \Pi ^{ur} _{i} (\overline X_{K};x_0;x_1)^0$$
\begin{lemma} 
For any $r$ $dim\, L_r= 1$. 
\end{lemma}
{\bf Proof} : Induction on $r$. Consider the following diagram:
$$
\def\normalbaselines{\baselineskip20pt
\lineskip3pt  \lineskiplimit3pt}
\def\mapright#1{\smash{
\mathop{\longrightarrow}\limits^{#1}}}
\def\mapdown#1{\Big\downarrow\rlap
{$\vcenter{\hbox{$\scriptstyle#1$}}$}}
\matrix{0 &    & 0 \cr
 \mapdown{} &  & \mapdown{ } \cr
gr ^W_{-i} \Pi^{ur}_{i+1} (\overline X_K)^{-2i} 
 & \stackrel{N^i}{\simeq} & gr ^W_{-i} \Pi^{ur}_{i+1} (\overline X_K)^{-2i}  \cr
\mapdown{} &  & \mapdown{Id } \cr
L_{i+1} & \mapright{N^i} &  gr ^W_{-i} \Pi^{ur}_{i+1} (\overline X_K)^{-2i} \cr 
\mapdown{p_{i+1}} &  & \mapdown{ } \cr 
L_i & \mapright{} &  0 \cr}
$$
Theorem \ref{M} implies that the first  horizontal arrow is an isomorphism.
It proves the lemma.
 
We have already proven the main theorem in the case when  $X_K=\overline X$. The general case follows from the following lemma.
\begin{lemma}
$j_*: \, \Pi^{ur}_{un} (X_K)^{0} \simeq \Pi^{ur}_{un} (\overline X_K)^{0}$
\end{lemma}
{\bf Proof}. It suffices to check that
\begin{equation} 
gr ^W_{\bullet}\Pi^{ur}_{un} (X_K)^{0} \simeq gr ^W_{\bullet}\Pi^{ur}_{un} (\overline X_K)^{0}
\label{oh}
\end{equation}
It is easy to see that 
the decomposition  
$$gr ^W_{\bullet}\Pi^{ur}_{un} (\overline X_K)=\bigoplus _{i\leq 0} (gr ^W_{\bullet}\Pi^{ur}_{un})^i$$ 
is compatible with the algebra structure. In particularly, (\ref{oh}) is a homomorphism of algebras.
Moreover, the algebra  $gr ^W_{\bullet}\Pi^{ur}_{un} (\overline X_K)^{0} $ is quadratic. Hence, it suffices  to check that (\ref{oh})
is an isomorphism in degrees $1$ and $2$. Assume that $\overline X_K$ has a strict semi-stable model. That is an open immersion 
 $\overline X_K \hookrightarrow   \overline X$ together with a closed subscheme $Z\hookrightarrow \overline X$ ,
$Z_K= X_K$ such that $(Z;\overline X)$ is  a strict semi-stable pair. Let $Z=\bigcup _{i\in I} Z_i $ be the decomposition 
in the irreducible components. We have the following exact sequence
\begin{equation}\label{somnenie}
0 \longrightarrow H^1_{st}(\overline X_K)\longrightarrow H^1_{st}( X_K)\longrightarrow \bigoplus _{i\in I} K_0(-1) \longrightarrow
 H^2_{st}(\overline X_K)\longrightarrow  H^2_{st}( X_K)
\end{equation}
Of course,  $ H^i_{st}(\overline X_K)$,  $H^i_{st}( X_K) $ are just notations for log crystalline cohomology of $\overline X$ with 
the evident log structures. It follows that $ H^1_{st}( X_K)^0 \simeq  H^1_{st}(\overline  X_K)^0 $ and the map
$ H^2_{st}(\overline X_K)^0\longrightarrow  H^2_{st}( X_K)^0$ is an injection. 

It proves  the lemma in the semi-stable case.   As usual, in general one  can make use of de Jong's theorem.

It completes the proof 
of Theorem \ref{main} along with Theorem B - the element $C_{\cdot}(X_K,x_0;x_1)$ gives rise to  a canonical
parallel translation.
\section{Mixed unipotent $F$-crystals}
In this section, for a smooth variety $X_{\overline K}$ over $spec \, \overline K$,
 we construct a certain category of  `` mixed unipotent $F$-crystals ``. 
This is a relative version of the concept of a mixed ${\cal H}^{ur}_1$-module (see {\bf 7.2.} ). 
As an application we construct a weight filtration  on $\Pi^{ur}_r (X_K;x_0;x_1)$, for any smooth scheme, and prove  
 prove the latter module together with the filtration is mixed.

{\bf 9.1. Mixed log $F$-crystals }. Throughout  this subsection we keep  the notations of {\bf 5.2.}
 Denote by ${\cal C}^{\phi}(\overline X  )$ the category of coherent log $F$-crystals on 
 $\overline X / spec \, W(k)$  ($spec \, W(k)$ is endowed with the  trivial log structure). 

An object $E$ of  ${\cal C}^{\phi}(\overline X ) \otimes \mathbb{Q}$ is called {\it pure of weight} $i$ 
if for any finite
extension $L\supset K$ and a point 
$i: spec \, R_L \to \overline X  $, with $i(spec\, L) \subset \overline X - Z $, 
the ${\cal H}^{ur}_1$-module 
$\Psi^{un}(i^*_{cris}E)\otimes \mathbb{Q}_p^{ur}$ is pure of weight $i$ (see {\bf 7.1}).

Next, we define the category ${\cal C}^{\phi}_{W_{\bullet}} (\overline X )$
 of mixed $F$-crystals: an object of the latter category  is a pair $(E, W_{\bullet}E)$, where
 $E$ is in ${\cal C}^{\phi}(\overline X ) \otimes \mathbb{Q}$ and $W_{\bullet}E $ 
is a filtration satisfying the following condition:

for any integer $i$ the crystal $gr^W_iE$ is pure of weight $i$.

It follows from Proposition \ref{bei} that  ${\cal C}^{\phi}_{W_{\bullet}} (\overline X )$ is an abelian 
category.

{\bf Remark.} Thanks to the results of {\bf 3.8.} we have: 
$$i^*: \{ \textstyle{coherent \, F-crystals \; on \; \overline X} \} \otimes \mathbb {Q}\simeq
\{ \textstyle{coherent \, F-crystals \; on \; \overline X_k} \} \otimes \mathbb {Q}$$
It shows that the category  ${\cal C}^{\phi}_{W_{\bullet}} (\overline X )$ depends on the special fibre $X_k$ only.
In particularly, for any finite totally ramified extension $L \supset K$ with the ring of integers $R_L$, 
  there is a canonical equivalence:
$$ {\cal C}^{\phi}_{W_{\bullet}} (\overline X ) \simeq {\cal C}^{\phi}_{W_{\bullet}} (\overline X \times spec \, R_L) $$
 
{\bf 9.2. Mixed unipotent $F$-crystals}.
 We define  a category ${\cal C}^{\phi}_{W_{\bullet}, \, un} (\overline X)$ to be
the full subcategory of ${\cal C}^{\phi}_{W_{\bullet}} (\overline X  )$ whose objects are mixed crystals
 $(E, W_{\bullet}E)$ with the following properties:

i) For any $i$ the $F$-crystals $gr^W_iE$
 are lifted from the logarithmic 
point $spec \,  R$.

ii)For any $i$, the vector bundle $W_iE/W_{i-2}E$ over $\overline X_K$ together with the logarithmic connection
is  {\sf a local system on $\overline X_K$} (i.e. the logarithmic connection is a connection in the usual sense).

We denote by   ${\cal C}^{\phi}_{W_{\bullet}, \, un} (\overline X \times spec\,  R_{K^{ur}})$ the 
injective limit of categories ${\cal C}^{\phi}_{W_{\bullet}, \, un} (\overline X \times spec\,  R_{K^{\prime }})$,
 where $K ^{\prime}$ runs over all 
finite unramified extensions of $K$. 

Let $X_{\overline K}$ be a smooth variety  over ${\overline K}$.
 We are going to  make use of de Jong's Theorem and the technique of descent 
(see {\bf 6.4} and {\bf 6.5}) to define a certain category  
${\cal C}^{\phi}_{W_{\bullet}, \, un} (X_{\overline K } )$ of
``mixed  unipotent $F$-crystals on $X_{\overline K }$ ''.

For this purpose, we choose a finite extension $L \supset \mathbb{Q}_p$,
an integral  proper, flat over $spec\,  R_L$  scheme  $\overline X$ together with an open immersion
$X_{\overline K} \hookrightarrow \overline X \times spec \, R_{\overline K}$ and a resolution (\ref{ba})
\begin{equation}
 \mathbf{\overline X}_1 \;
{\begin{array}{l} 
\overline p_1 \\
\longrightarrow \\
\overline p_2 \\
\longrightarrow \\
\end{array}} 
\; \mathbf{\overline X}_0
\stackrel{\overline f}{\longrightarrow} \overline X 
\label{resolution}
\end{equation} 
 We define  ${\cal C}^{\phi}_{W_{\bullet}, \, un} (X_{\overline K} )$ to be the full subcategory
of  $ {\cal C}^{\phi}_{W_{\bullet}, \, un} (\overline X_1 \times spec \,  R_{L ^{ur}} )$,
 which consists of crystals $E$ such that 
$$ \overline p_1^*E \simeq \overline p_2 ^*E$$
Using Proposition \ref{a} one can easily check that there is a canonical equivalence between 
the latter categories for different choices of the resolution.

Using Proposition \ref{a} one can easily check that there is a canonical equivalence between 
the latter categories for different choices of the resolution.

It also follows that there is a faithful functor
\begin{equation}\label{forget}
   {\cal C}^{\phi}_{W_{\bullet}, \, un} 
(X_{\overline K} ) \longrightarrow {\cal D}{\cal M}_{un}(X_{\overline K})
\end{equation}

{\bf Example.} {\it The category    ${\cal C}^{\phi}_{W_{\bullet}, \, un} (spec \, \overline K )$ 
is equivalent to the category 
of mixed  $ {\cal H}_1^{ur}$-modules.}

For a  mixed  $ {\cal H}_1^{ur}$-module $V$, we denote by $\underline V$ its pull-back on $X_{\overline K}$.

Given another smooth variety $ Y_{\overline K}$ and  a morphism  $f: Y_{\overline K} \to X_{\overline K} $ 
we can construct the restriction functor
$$f^*: {\cal C}^{\phi}_{W_{\bullet}, \, un} (X_{\overline K} )\longrightarrow 
 {\cal C}^{\phi}_{W_{\bullet}, \, un} (Y_{\overline K} )$$
In particularly, a point $i: spec \, \overline K \to X_{\overline K}$ defines a functor
$$i^*: {\cal C}^{\phi}_{W_{\bullet}, \, un} (X_{\overline K} )\longrightarrow
  {\cal C}^{\phi}_{W_{\bullet}, \, un} (spec \, \overline K )$$
$$ \longrightarrow  \{ \textstyle{Mixed} \; {\cal H}_1^{ur}-modules \}$$

For an object $E$ of ${\cal C}^{\phi}_{W_{\bullet}, \, un} (X_{\overline K} )$ one can define the log crystalline 
cohomology $H^{i}_{pst} (X_{\overline K}; E )$. This is a finite-dimensional  $ {\cal H}_1^{ur}$-module. 
There is a canonical isomorphism:
  $$H^{i}_{pst} (X_{\overline K}; E ) \otimes _{\mathbb{Q}_p^{ur}} \overline K 
\simeq  H^i_{DR}(X_{\overline K}; E)$$       

{\bf 9.3.} Given two objects $E$, $G$ of  the category ${\cal C}^{\phi}_{W_{\bullet}, \, un} (X_{\overline K} )$ we    
denote by ${\cal E}xt^i(E;G)$ the $ {\cal H}_1^{ur}$-module $H_{pst}^i(X_{\overline K}; E^*\otimes G))$.
 
For $i = 0,\,  1 $ the latter possesses a canonical  filtration  
$W_{\bullet}\, {\cal E}xt^i(E;G) \subset {\cal E}xt^i(E;G)$. The definition of the filtration on 
${\cal H}om(E;G)$ is obvious.
We  define the  filtration on  ${\cal E}xt^1(E;G)$  in the following way: choose a smooth compactification
 $X_K \hookrightarrow \overline X_{\overline K}$ and  let $M= E^*\otimes G $ and 
$W_{i+1}H^{1}_{pst}( X_{\overline K}; W_i M )$ be the preimage of
$$ H^{1}_{pst}( \overline X_{\overline K}; W_iM /W_{i-1}M) \subset H^{1}_{pst}( X_{\overline K}; W_iM/W_{i-1} M)$$
under the canonical map  
$$H^{1}_{pst}( X_{\overline K}; W_i M) \to  H^{1}_{pst}( X_{\overline K}; W_iM/W_{i-1} M)$$ 
Finally, we put  
$$W_{i+1}H^{1}_{pst}( X_{\overline K};  M) = Im \, W_{i+1}H^{1}_{pst}( X_{\overline K};  W_i M)$$
  Clearly,
the filtration does not depend on the choice of a compactification.

 \begin{proposition}\label{fundamental} 
For $i=0, \, 1$,  the pair $({\cal E}xt^i(E;G);W_{\bullet}\, {\cal E}xt^i(E;G))$ is a mixed   $ {\cal H}_1^{ur}$-module.
\end{proposition}
\begin{proof}

Without loss of generality we can assume that $X_{\overline K}$ is quasi-projective. Indeed, in general, we can
choose a smooth compactification $X_{\overline K} \hookrightarrow \overline X_{\overline K}$,
 a smooth projective variety $Y_{\overline K}$ and a morphism
$f: \, Y_{\overline K} \to  \overline X_{\overline K}$ which is a birational isomorphism. 
The existence of the pair $(Y_{\overline K}; f)$ follows from 
the Chow lemma. It is easy to see that the morphism $f:\,  f^{-1}(X_{\overline K}) \to X_{\overline K}$  induces 
an isomorphism on the fundamental groups. 

Next, we make use of the following lemma.
\begin{lemma}
Let $U \subset \mathbb{P}^n_{\overline K}$ be a quasi-projective variety of dimension $\geq 3$.
Then, for a general hyperplane section $U\cap H \hookrightarrow U$ , the restriction functor 
  $${\cal D}{\cal M}_{un}(U) \longrightarrow {\cal D}{\cal M}_{un}(U\cap H)$$
is an equivalence of the categories.
\end{lemma}
Proof is left to the reader.

It follows that there exists  a smooth surface $Z_{\overline K}$ and a morphism 
$g: \, Z_{\overline K} \to X_{\overline K}$ such that 
$$g^* \, {\cal D}{\cal M}_{un}(X_{\overline K}) \longrightarrow {\cal D}{\cal M}_{un}(Z_{\overline K}) $$
is an equivalence. In turn, it implies that the pull-back functor
 $$ g^* \, {\cal C}^{\phi}_{W_{\bullet}, \, un} (X_{\overline K } ) \to 
{\cal C}^{\phi}_{W_{\bullet}, \, un} (Z_{\overline K } )$$
 is fully faithful. Hence, is suffices to prove the proposition for surfaces. 

Of course, we can assume that  $E$ is the trivial $F$-crystal (i.e $E= \underline{\mathbb{Q}_p^{ur}}$).

 We prove the proposition by induction on length of the weight filtration on $G$.
  Let $r$ be the integer such that $W_r G \ne  0$ and $W_{r-1}G=0$.
Consider the following exact sequence of $ {\cal H}_1^{ur}$-modules:
$$0\longrightarrow  W_r G 
\longrightarrow H^{0}_{pst} (Z_{\overline K} ;G)  \longrightarrow H^{0}_{pst} (Z_{\overline K};  G/ W_r G)\longrightarrow
 W_r G \otimes H^{1}_{pst}(Z_{\overline K}) \longrightarrow         $$
$$H^{1}_{pst} (Z_{\overline K} ;G) \longrightarrow H^{1}_{pst} (Z_{\overline K};  G/ W_r G) \longrightarrow
 W_r G \otimes H^{2}_{pst}(Z_{\overline K})$$  
By the induction hypothesis $ H^{i}_{pst} (Z_{\overline K};  G/ W_r G)   $ are  mixed.

On the other hand, the result of Mokrane (see {\bf 7.4.}) implies that  
the same is true for $ H^{i}_{pst}(Z_{\overline K}) $. 
Hence,   $W_r G   \otimes H^{i}_{pst}(Z_{\overline K}) $  are  mixed.
 
Finally, it is easy to check that all the morphisms 
in the above sequence are strictly compatible with the weight filtration. Now the proposition follows from lemma 
\ref{extension}.
\end{proof}  
\begin{cor}\label{Leray} 
We have the following exact sequence:
$$0 \longrightarrow Ext^1_{{\cal C}^{\phi}_{W_{\bullet}, \, un} ( spec \, \overline K)}
( \mathbb{Q}^{ur}_p; {\cal H}om(E;G)) \longrightarrow $$
$$ Ext^1_{{\cal C}^{\phi}_{W_{\bullet}, \, un} (X_{\overline K} )}(E;G) 
\longrightarrow  Hom_{{\cal C}^{\phi}_{W_{\bullet}, \, un} ( spec \,\overline K )} (\mathbb{Q}^{ur}_p;
{\cal E}xt^1(E; G)) $$
\end{cor}

{\bf 9.4} For an integer $r$ we denote by ${\cal C}^{\phi}_{W_{\bullet}, \, r} (X_{\overline K} )$
 the full subcategory of 
 ${\cal C}^{\phi}_{W_{\bullet}, \, un} (X_{\overline K} )$ whose objects are unipotent $F$-crystals of level $r$.
\begin{theorem} 
There is a unique object  $\Pi_r^{ur}(X_{\overline K})$ of 
${\cal C}^{\phi}_{W_{\bullet}, \, un} (X_{\overline K}\times X_{\overline K} )$
characterised by the following property:
 
for any $E$ in  ${\cal C}^{\phi}_{W_{\bullet}, \, r} (X_{\overline K}\times X_{\overline K} )$
 there is a canonical isomorphism of  mixed $ {\cal H}_1^{ur}$-modules
$$ Hom(\Pi_r^{ur}(X_{\overline K}); E) \simeq H^{0}_{pst}( X_{\overline K}; \Delta ^* E)$$
Here $\Delta: X_{\overline K} \to X_{\overline K} \times X_{\overline K}$ stands for the diagonal embedding.

Moreover, the de Rham realization functor (\ref{forget}) sends  $\Pi_r^{ur}(X_{\overline K})$ to  
 $\Pi^{DR}_r(X_{\overline K} )$.
\end{theorem}

{\bf Proof.}

 We prove the theorem  by induction on $r$. 
If $r=1$ there is nothing to prove. Assume that  $\Pi_{r-1}^{ur}(X_{\overline K})$ is already constructed.
(And the assertion about the de Rham realization holds.)

The identity morphism $Id: \Pi_{r-1}^{ur}(X_{\overline K}) 
\to \Pi_{r-1}^{ur}(X_{\overline K})$ gives rise to a  morphism of mixed modules:
\begin{equation}
        \mathbf{1}:\underline{\mathbb{Q}^{ur}_p} \to   \Delta ^* \Pi_{r-1}^{ur}(X_{\overline K})
\end{equation}
Denote by $I$ the kernel of the map:
  $$ {\cal E}xt^1_{{\cal C}^{\phi}_{W_{\bullet}, \, un} (X_{\overline K}\times X_{\overline K} )}
(\Pi_{r-1}^{ur}(X_{\overline K});\mathbb{Q}^{ur}_p) \stackrel{\mathbf{1}}{\longrightarrow} 
{\cal E}xt^1_{{\cal C}^{\phi}_{W_{\bullet}, \, un}} (X_{\overline K} ) (\mathbb{Q}^{ur}_p;\mathbb{Q}^{ur}_p)$$ 
Proposition \ref{fundamental} implies that $I$ is a mixed $ {\cal H}_1^{ur}$-module.
 Note that the de Rham realization of $\underline I^*$ is canonically isomorphic to 
$ker\,  (\Pi^{DR}_r(X_{\overline K}  ) \to \Pi^{DR}_{r-1}(X_{\overline K } )$ 
\begin{lemma}
There is a unique element 
$[\Pi_{r}^{ur}(X_{\overline K})] 
\in Ext^1_{{\cal C}^{\phi}_{W_{\bullet}, \, un} (X_{\overline K}\times X_{\overline K} )}
(\Pi_{r-1}^{ur}(X_{\overline K});\underline I^*)$ satisfying the following properties:

i) The de Rham realization of $\Pi_{r}^{ur}(X_{\overline K})$ is $\Pi^{DR}_r(X_{\overline K}  )$.

ii)The image of $\Pi_{r}^{ur}(X_{\overline K})$ under the canonical morphism 
 $$ Ext^1_{{\cal C}^{\phi}_{W_{\bullet}, \, un} (X_{\overline K}\times X_{\overline K}  )}
(\Pi_{r-1}^{ur}(X_{\overline K});\underline I^*) \stackrel{\mathbf{1}}{\longrightarrow}  
Ext^1_{{\cal C}^{\phi}_{W_{\bullet}, \, un} (X_{\overline K} )}
(\underline{\mathbb{Q}^{ur}_p};\underline I^*)$$
is equal to $0$.
\end{lemma}
\begin{proof} 
First, we make use of Corollary \ref{Leray} and consider  the following commutative diagram:
 $$
\def\normalbaselines{\baselineskip20pt
\lineskip3pt  \lineskiplimit3pt}
\def\mapright#1{\smash{
\mathop{\to}\limits^{#1}}}
\def\mapdown#1{\Big\downarrow\rlap
{$\vcenter{\hbox{$\scriptstyle#1$}}$}}
\begin{matrix}{  Ext^1_{{\cal C}^{\phi}_{W_{\bullet}, \, un} ( spec \, \overline K )}
(\mathbb{Q}^{ur}_p       ; I^*)
&  \mapright{} & Ext^1 
 (\Pi_{r-1}^{ur}(X_{\overline K});\underline I^*) & \mapright {}
 {\cal E}xt^1(\Pi_{r-1}^{ur}(X_{\overline K});\underline I^*)  \cr
   \mapdown{Id} &    & \mapdown{\mathbf{1 \circ \Delta ^*}} &   &         \cr
  Ext^1_{{\cal C}^{\phi}_{W_{\bullet}, \, un} ( spec \, \overline K )}
(\mathbb{Q}^{ur}_p ; I^*)
 & \mapright{Id} & 
Ext^1_{{\cal C}^{\phi}_{W_{\bullet}, \, un} ( spec \, \overline K )}
(\mathbb{Q}^{ur}_p ; I^*)    &        & \cr }    
\end{matrix}
 $$
The first row is exact by {\it loc. cit.} It implies the uniqueness of $\Pi_{r}^{ur}(X_{\overline K})$.

It remains to prove the existence.
Choose a resolution (\ref{resolution}).  By the definition,  the category 
${\cal C}^{\phi}_{W_{\bullet}, \, un} (X_{\overline K}\times X_{\overline K} )$ is identified with a 
certain subcategory of ${\cal C}^{\phi}_{W_{\bullet}, \, un} (\overline X_{1}\times \overline X_1 \times 
spec \, R_{L^{ur}})$. 

Lemma \ref{obman} implies that there exists a unique extension 
\begin{equation}\label{vseo}
 0 \longrightarrow \underline I^* \longrightarrow    \Pi_{r}^{ur}(X_{\overline K}) 
\longrightarrow  \Pi_{r-1}^{ur}(X_{\overline K}) \longrightarrow 0
\end{equation}
in the category ${\cal C}(\overline X_0 \times \overline X_0 \times 
spec \, R_{L^{ur}}  )
\otimes \mathbb{Q}$ satisfying the properties i) and ii).

Fix a lifting:
\begin{equation}\label{vseodno}
        \mathbf{1}: \underline{\mathbb{Q}^{ur}_p} \to   \Delta ^* \Pi_{r}^{ur}(X_{\overline K})
\end{equation}  
 It is easy to see that there is a unique $F$-structure on the crystal $\Pi_{r}^{ur}(X_{\overline K})$ such that
all the morphisms in \ref{vseo} and \ref{vseodno} commute with $\phi$.

Finally, we define a weight filtration on $\Pi_{r}^{ur}(X_{\overline K})$. 

{\bf Sublemma 1.} The class of extension (\ref{vseo}) is in 
$W_0{\cal E}xt^1( \Pi_{r-1}^{ur}(X_ {\overline K});\underline I^*)$.

Proof is omitted. 

Consider the extension 
\begin{equation}\label{uraa}
 0 \longrightarrow \underline{I^*/W_i I^*} \longrightarrow \mathbf{?} \longrightarrow W_i \Pi_{r-1}^{ur}(X_{\overline K}) 
\longrightarrow 0
\end{equation}
induced by  (\ref{vseo}).

{\sf Notation:} for an extension of crystals $0\to E \to G \to K \to 0$, we let  
$Spl(K;G)$ be the preimage of $Id \in {\cal H}om(K;K)$ under a canonical map ${\cal H}om(K;K)\to {\cal H}om(K;G)$.

We have a canonical map:
$$ Spl(W_{i} \Pi_{r-1}^{ur}(X_{\overline K});\mathbf{?}) \longrightarrow
{\cal E}xt^1(\Pi_{r-1}^{ur}
(X_{\overline K})/W_{i} \Pi_{r-1}^{ur}(X_{\overline K});\underline{I^*/W_i I^*})$$  

Let $P$ be the preimage of $W_0{\cal E}xt^1(\Pi_{r-1}^{ur}
(X_{\overline K})/W_{i} \Pi_{r-1}^{ur}(X_{\overline K});\underline{I^*/W_i I^*})$ under the latter map.

{\bf Sublemma 2.} For any $i<0$ , the set $P$ consists of one element.

{\bf Proof}. Consider the following exact sequence
$${\cal H}om(\Pi_{r-1}^{ur}(X_{\overline K});\underline{I^*/W_i I^*})\longrightarrow
{\cal H}om(W_i\Pi_{r-1}^{ur}(X_{\overline K});\underline{I^*/W_i I^*})\longrightarrow $$
$${\cal E}xt^1(\Pi_{r-1}^{ur}
(X_{\overline K})/W_{i} \Pi_{r-1}^{ur}(X_{\overline K});\underline{I^*/W_i I^*})\stackrel{\alpha}{\longrightarrow}
{\cal E}xt^1(\Pi_{r-1}^{ur}(X_{\overline K});\underline{I^*/W_i I^*})$$
If $i<0$ the first map is $0$. Moreover, $ Spl(W_{i} \Pi_{r-1}^{ur}(X_{\overline K});\mathbf{?})=
\alpha^{-1}([\Pi_{r}^{ur}(X_{\overline K})]) $. Hence, Sublemma 1 implies that $P$ is not empty.
On the other hand, $W_0{\cal H}om(W_i\Pi_{r-1}^{ur}(X_{\overline K});\underline{I^*/W_i I^*})=0$. Therefore,
 $P$ consists of one element.

It is clear, that the element of $P$ defines a splitting of (\ref{uraa}).
In turn, the latter gives rise to a subcrystal $W_i\Pi_{r}^{ur}(X_{\overline K}) \subset \Pi_{r}^{ur}(X_{\overline K})$.
It completes the proof of the lemma.

\end{proof}
We leave to the reader to check that the mixed crystal $\Pi_{r}^{ur}(X_{\overline K})$ 
represents the functor  defined in the theorem.

{\bf 9.5.} Given a point $i_{x_0; x_1}: spec \, \overline K \to X_{\overline K} \times X_{\overline K}$, 
we define a mixed ${\cal H}_1^{ur}$-module
$$ \Pi^{ur}_r (X_{\overline K};x_0;x_1) := i^*_{x_0; x_1}\Pi_r(X_{\overline K})$$
It is easy to see that $ \Pi^{ur}_r (X_{\overline K};x_0;x_0)$ is {\sf an algebra} in the tensor category
${\cal C}^{\phi}_{W_{\bullet}, \, un} ( spec \,{\overline K} )$.  
 
\begin{cor} The restriction to a point $x_0 \in X_{\overline K}(\overline K)$ 
defines an equivalence between the category   
 ${\cal C}^{\phi}_{W_{\bullet}, \, un} (X_{\overline K} )$ and  the category of 
modules over $ \Pi^{ur}_r (X_{\overline K};x_0;x_0)$ in ${\cal C}^{\phi}_{W_{\bullet}, \, un} ( spec \,\overline K )$.
\end{cor}

{\bf 9.6. Remark:}   One
can show that   for any smooth variety $X_K$ over $spec \, K$ 
the underlying $ {\cal H}_1^{ur}$-module  $\Pi^{ur}_r (X_K \times spec \, \overline K ;x_0;x_1)$ 
is isomorphic to one defined in Section 6.

{\bf 9.7.  Proposition.} 
{\it Let ${\cal V}_r\subset \Pi^{ur}_r (X_K;x_0;x_1)$ be the subspace  which consists of elements 
$v\in  \Pi^{ur}_r (X_K;x_0;x_1)$ satisfying the following properties:

i)$\phi v= v$

ii)$ N^av \in W_{-a-1} \, \Pi^{ur}_r (X_K;x_0;x_1)$ for any $0<a<r$.

We have $dim {\cal V}_r = 1$. Moreover, the canonical morphism:
    $$ {\cal V}_r \longrightarrow   \mathbb{Q}_p$$
is an isomorphism.} 

Proof is left to the reader.  (Use induction on $r$.)

In particularly, it defines a distinguished 
 element $C_r (X_K;x_0;x_1) \in \Pi^{ur}_r (X_K;x_0;x_1)$ (the preimage of $1\in \mathbb{Q}_p^{ur}$).
We have just reproved Theorem B ! 
 
\section{Variations of p-adic Hodge structures}
{\bf 10.1. Notations.} Let  $\overline X$ be  a smooth proper scheme over $spec \, W(k)$,
$Z\subset \overline X$ be a normal crossings divisor relative to $spec \, W(k) $,  $X= \overline X - Z$, $X_k = X \times _{spec \, W(k)} spec\, k$,
 $X_K = X \times _{spec \, W(k)} spec \, K$.

We endow $\overline X$  with the log structure given by the divisor $Z$
and $spec \, W(k) $  with the trivial log structure (see {\bf 2.2.}).

 {\bf 10.2.} 
Faltings defined a certain category ${\cal M}{\cal F}_{[a;a+b]}( \overline X)$ of filtered $F$-crystals on $\overline X$, analogous to the category of variations of Hodge 
structures over the field of complex numbers. 
We copy  the definition from [Fa]  with merely decorative innovations.
Fix integers $a$, $b$, with $ 0 \leq b \leq p-1 $.

{\bf Definition.} An object of  ${\cal M}{\cal F}_{[a;a+b]}( \overline X)$  consists of the following data:

a)A coherent log  $F$-crystal $E$ on $\overline X$.

 c)A filtration on the corresponding coherent sheaf  $E_{\overline X}$ on $\overline X$
$$ \cdot \cdot \cdot\subset F^i E_{\overline X}\subset F^{i-1}E_{\overline X} \subset \cdot \cdot \cdot\subset E_{\overline X} $$
 by coherent subsheaves satisfying Griffiths-transversality:  
$$\nabla ( F^iE_{\overline X})\subset F^{i-1}E_{\overline X}\otimes \Omega^1_{\overline X}(d log\infty )$$ 
 and the  following stabilisation property: $ F^a E_{\overline X}= E_{\overline X}$,  $F^{a+b+1} E_{\overline X}=0$.
    
These are  subject to the following conditions:

i) The sheaves $F^iE_{\overline X}/F^{i+1}E_{\overline X}$ are vector bundles on $\overline X$.

To formulate next  condition , we first note  that the coherent sheaf $(Fr^*E)_{\overline X}$ possesses a canonical filtration:
                $$ \widetilde {F}^i (Fr^*E)_{\overline X}  \subset (Fr^*E)_{\overline X} $$
It can be constructed using local liftings of the Frobenius: if $U \subset \overline X$ together with a logarithmic Frobenius-lift  
$\widetilde{Fr}: \widehat{U} \to \widehat{U}$ ($\widehat{U}$ stands for $p$-adic completion), we define
      $$ \widetilde {F}^i (Fr^*E)_{\widehat{U} }    =\widetilde{Fr}^*(\sum _{k+m= i; \, m \geq 0} p^m F^{k}E_{\widehat {U} } )$$
If $b \leq p-1 $  the latter does not depend on the choice of filtration.

Or, equivalently, we can consider $E$ as a filtered crystal on $(\overline X_k/spec\, W(k))$ 
(see {\bf 2.5.}) and take the corresponding filtration on
$(Fr^*E)_{\overline X}$.

ii) The restriction of 
$$\phi  : Fr^*E_{\overline X } \to E_{\overline X } $$
to  $ \widetilde {F}^i (Fr^*E)_{\overline X}$ is divisible by $p^i$. 

iii) $\sum _i  p^{-i} \phi (\widetilde {F}^i (Fr^*E)_{\overline X})=  E_{\overline X } $ .

 Denote 
 $${\cal M}{\cal F}_{[a;a+b]}(\overline X ) \otimes \mathbb{Q} ={\cal M}{\cal F}_{[a;a+b]}^{\mathbb{Q}} (\overline X ) $$
We use the name variation of p-adic Hodge structure for an object of the latter category. 

Let $\underline{K_0}$ stand for the trivial Hodge structure on $spec \, W(k)$.

{\bf 10.3.}
The following result is proven by Faltings.
\begin{theorem} 
a) Let $E$ and $G$ be objects of    ${\cal M}{\cal F}_{[a;a+b]}^{\mathbb{Q}} (\overline X ) $  and  $f:\, E\to G $ is  morphism. 
Then $f$ is strict for the filtrations. 

b)The category  ${\cal M}{\cal F}_{[a;a+b]}^{\mathbb{Q}} (\overline X ) $ is abelian.

\end{theorem}
Because of  the restriction on $b$ the category of variations of Hodge structures does not possess a tensor structure.
Nevertheless, if $ 0\leq b_1, b_2$, $b_1+b_2 \leq p-1 $ we can define a tensor product:
$$ {\cal M}{\cal F}_{[a_1;a_1+b_1]}^{\mathbb{Q}  }(\overline X ) \times {\cal M}{\cal F}_{[a_2;a_2+b_2]}^{\mathbb{Q} }(\overline X ) 
\to {\cal M}{\cal F}_{[a_1+a_2;a_1+a_2+b_1+b_2]}^{\mathbb{Q} }(\overline X )$$
Under a similar restriction, we can define ${\cal H}om $  between two variations.

{\bf 10.4.}  Let  $\overline Y$ be another  smooth proper scheme over $spec \, W(k)$ 
together with a normal crossings divisor         $D\subset \overline Y$  relative to $spec \, W(k) $ and $f: \overline X \to \overline Y$ 
be a proper log smooth morphism of relative dimension $d$.
Assume that $f^*(D) \leq Z $ as divisors (i.e. $f$ is smooth in codimension $1$).
\begin{theorem}\label{push} ({\sf Faltings.})
 Let  $E$ be an object  of  ${\cal M}{\cal F}_{[a;a+b]}^{\mathbb{Q}} (\overline X ) $. For any $i$  with  $b+min(d,i)\leq p-2$  
 the filtered $F$-crystal  $R^if_*E$ is in   ${\cal M}{\cal F}_{[a;a+b+i]}^{\mathbb{Q}} (\overline Y ) $.
\end{theorem}
{\bf 10.5.}We have the evident  faithful functor 
      $$ {\cal M}{\cal F}_{[a;a+b]} ^{\mathbb{Q}}(\overline X )  \longrightarrow  {\cal D}{\cal M}(X_{K_0})$$
 to the category of vector bundles on $ X_{K_0}$ together with an integrable connection.

Let  $E$ and $G$ be objects of ${\cal M}{\cal F}_{[a;a+b]} ^{\mathbb{Q}}(\overline X ) $ and  $2b \leq p- 3$. 
Thanks to the previous result the vector spaces
  $$ Hom _{ {\cal D}{\cal M}(X_{K_0})}(E;G) \; , \;  Ext^1 _{ {\cal D}{\cal M}(X_{K_0})}(E;G)$$
carry canonical Hodge structures.
\begin{proposition}\label{ri} 
 We have the following exact sequence
 $$0 \longrightarrow  Ext^1(\underline{K_{0}}; Hom_{ {\cal D}{\cal M}(X_{K_0})}(E;G))
\stackrel{\alpha}{\longrightarrow}
 Ext^1_{ {\cal M}{\cal F}_{[a;a+b]} ^{\mathbb{Q}}(\overline X )}(E;G)$$ 
$$\stackrel{\beta}{\longrightarrow} ( Ext^1 _{ {\cal D}{\cal M}(X_{K_0})}(E;G)^{\phi=1}\cap F^0 \longrightarrow 0 $$
\end{proposition}
\begin{proof}
Given an extension 
       $$ 0\longrightarrow {\cal H}om(E;G) \longrightarrow A \longrightarrow \underline{K_{0}} \longrightarrow 0 $$
whose class $[A]$ in 
$Ext^1_{ {\cal M}{\cal F}_{[a;a+b]} ^{\mathbb{Q}}(\overline X )}(\underline{K_{0}};{\cal H}om(E;G))=
 Ext^1_{ {\cal M}{\cal F}_{[a;a+b]} ^{\mathbb{Q}}(\overline X )}(E;G)$ satisfies the property:
$\beta ([A])= 0$ we define $\gamma ([A]) \in
Ext^1(\underline{K_{0}}; Hom_{ {\cal D}{\cal M}(X_{K_0})}(E;G)) $ to be the class of the extension
   $$ 0\longrightarrow  Hom_{ {\cal D}{\cal M}(X_{K_0})}(E;G)    
\longrightarrow Hom _{ {\cal D}{\cal M}(X_{K_0})}( \underline{K_{0}}; A) \longrightarrow \underline{K_{0}} \longrightarrow 0 $$
It gives rise to an isomorphism $\gamma: ker \, \beta \simeq Ext^1(\underline{K_{0}}; Hom_{ {\cal D}{\cal M}(X_{K_0})}(E;G))$ with
$\gamma \alpha = 1$. 

It remains to show that the map $\beta$ is surjective. Let 
   $$ 0\longrightarrow G \longrightarrow B \longrightarrow E \longrightarrow 0 $$
be an extension with $[B] \in  ( Ext^1 _{ {\cal D}{\cal M}(X_{K_0})}(E;G)^{\phi=1}\cap F^0 $. Since
  $$   Ext^1 _{ cris}(E;G) \otimes \mathbb{Q} \simeq  Ext^1 _{ {\cal D}{\cal M}(X_{K_0})}(E;G) $$
we may assume that $B$ is a log crystal on $\overline X$. Moreover, we can choose $F$-structure 
$\phi : Fr^* B \to  B $ and a filtration on $F^iB \subset B$ on the underlying coherent sheaf compatible with those on 
$E$ and $G$. Then, for sufficiently large $n$ we have:
         $$p^n [B] \in Ext^1_{ {\cal M}{\cal F}_{[a;a+b]}(\overline X )}(E;G)$$  
  It completes the proof.
\end{proof}

 \section{Construction of etale local systems}
{\bf 11.1.} Faltings  has constructed  a fully faithful functor 
$$ \mathbf{D}: {\cal M}{\cal F}_{[a;a+p-2]}(\overline X ) \longrightarrow Sh^{et}_{\mathbb{Z}_p}(X_{K_0})$$ 
to the category of etale locally constant $\mathbb{Z}_p$-sheaves on $X_{K_0}$ [Fa].

We reproduce his construction here with some minor variations.

{\bf 11.2.} Given a log crystal $E$  on $\overline X$  we put $E_n= E/p^n$.

 An etale open affine subset $U = spec \,{\cal R} \subset \overline X $ is called ${\sf small}$ if there exists an
etale map of log schemes 
$$ U \to spec W(k)[T_1, T_2, \cdot \cdot \cdot, T_m]$$
 where  the log structure on $spec \, W(k)[T_1, T_2, \cdot \cdot \cdot, T_m]$ is given by the divisor  $\prod _{1\leq i \leq m}T_i = 0$
We note the  scheme $\overline X$ can be covered by small open subschemes. 

Choose a small subscheme $U$ and denote by $\widehat{U} =Spf\, \widehat{{\cal R}} $ the corresponding $p$-adic formal scheme.
Let $\bar {\cal R} $ be the union of all finite extensions $\widehat{{\cal R}_k}$ of $\widehat{{\cal R}}$ such that
$ (p {\cal I}_Z)^N \Omega ^1 _{\widehat{{\cal R}_k}/ \widehat{{\cal R}}}=0 $ for sufficiently  large  $N$ (${\cal I}_Z$
is the ideal corresponding to the immersion $Z \hookrightarrow \overline X$). We denote  by $\widehat{\bar U}=Spf\widehat{\bar {\cal R}}$  
the spectrum of the $p$-adic completion of  $\bar {\cal R} $.
 
We endow $\widehat{\bar U}$ with a log structure $$M_{\widehat{\bar U}}= \{f \in {\cal O}_
{\widehat{\bar U}} |f^N \in p^*M_{\widehat{U}} \textstyle{\, for\; some\; integer\; N \,}\}$$
So  the map $p:\widehat{\bar U} \to \widehat{U} $ is a morhism of log schemes.

{\bf 11.3.} Consider the log crystalline cohomology 
 $$    H^0_{cris}(\widehat{\bar U}\times _{spec \, W(k)} speck / spec \, W(k), p^*E_{n} |_{ \widehat{U}})$$
The Galois group of $\bar {\cal R}$ over $\widehat{{\cal R}}$ acts on the latter space.
On the other hand,the cohomology group  is also equipped
 with an action of the operator $\phi^*$ induced by $\phi: Fr^*E \to E $.
 Finally it possesses  a  canonical filtration.  The latter  comes from the identification 
$$H^0_{cris}(\widehat{\bar U}/spec\, W(k), p^*E_{n}|_{\widehat{U}})
 =H^0_{cris}(\widehat{\bar U}\times _{spec \, W(k)} speck /spec \, W(k), p^*E_{n}|_{\widehat{U}})$$       
and the filtration on $E$.
 
One can prove that the higher cohomology groups vanish (this depends on the assumption that $U$ is small).

   We define   
$$ \mathbf {D} (E_{n}|_{\widehat{U}}) =  (H^0_{cris}(\widehat{\bar U}\times _{spec \, W(k)} speck / spec \, W(k), p^*E_{n}|_{\widehat{U}}))^{\phi=1}\cap F^0 $$

Faltings has proven  that  $\mathbf {D} (E_{n}|_{\widehat{U}})$
 is a {\sf  finite abelian group of the same type as $E_n$}
(meaning that the coherent sheaf $E_U$ is locally isomorphic to 
$\mathbf {D} (E_{n}|_{\widehat{U}}) \otimes {\cal O}_U$)  
 equipped with an $Gal(\bar {\cal R}/\widehat{{\cal R}})$-action.      
  
Choose a covering of $\overline X$ by small open subschemes $U_i$. It easily follows from the result cited  above that we can glue an etale
 locally constant sheaf  on $\overline X_{K_0}$ from $\mathbf {D} (E_{n}|_{\widehat{U_i}})$.
 We denote the latter by  $\mathbf {D} (E_n)$.

 The sheaves $\mathbf {D} (E_i)$ of 
 $\mathbb{Z}/p^i\mathbb{Z}$-modules define a  locally constant
$p$-adic sheaf $\mathbf {D} (E)$.

{\bf 11.4.} The following  result is proven by Faltings.
\begin{theorem} 
  The functor
$$\mathbf {D}: \; {\cal M}{\cal F}_{[a;a+p-2]}(\overline X )\longrightarrow Sh^{et}_{\mathbb{Z}_p}(X_{K_0})$$
is  exact and fully faithful. Its image is closed under subobjects and quotients.

The functor $\mathbf{D}$ is compatible with the (partly defined) tensor product and ${\cal H}om$. 
\end{theorem}
   
Clearly it implies that the same is true  for the functor 
$$  {\cal M}{\cal F}_{[a;a+p-2]}^{\mathbb{Q}}(\overline X ) 
\longrightarrow  Sh^{et}_{\mathbb{Q}_p,}(X_{K_0})$$

\section{$B_{DR}$}
{\bf 12.1.} We keep the notations of the previous section.

The main result  of this section is the following theorem .
\begin{theorem}\label{Bst}
Let $i_x: spec \,  K  \to X_{K_0}$ be a point and $E$ be an object of the category
${\cal M}{\cal F}_{[a;a+ \frac{p-3}{2}]}^{\mathbb{Q}}(\overline X ) $   
There is  a canonical isomorphism 
 $$B_{DR} \otimes _{K}
 i_x^* E_{\overline X} \simeq B_{DR} \otimes _{\mathbb{Q}_p} i_{x}^* \circ  \mathbf {D}(E) $$     
of vector spaces over $B_{DR}$.

The isomorhism is compatible with the tensor structure.
\end{theorem}
{\bf 12.2.} The proof is based on the following 
result which is due to Kato [Ka2] (in fact Kato formulated the result in a slightly less general form but the proof 
works in our case as well).

Let $G$ be a $F$-crystal on the  scheme $ spec \, R / spec \, W(k) $ equipped with the canonical  log structure: $M_{ spec \, R} = R-0 $. 
Remind (Theorem \ref{psi}) that $G$ gives rise to  a ${\cal H}_1$ -module $\Psi^{un}(E)$ together with an isomorphism
$$\Psi^{un}(G)\otimes_{K_0} K\simeq G_R \otimes _R K$$ 
We denote by $\overline R$ the integral closure of $R$ in $\overline K$. Endow  $spec \, \overline R$ with a log structure $M_{ spec \, \overline R} =\overline R -0 $.
Consider the natural morphism of log schemes 
             $$q:\; spec \, \overline R  \longrightarrow spec \, R $$           
 \begin{theorem}
The  kernel of   
$$N \otimes 1 +1 \otimes N : \;\Psi^{un} (G )\otimes _{K_0} B^+_{st} \longrightarrow\Psi^{un}( G) \otimes _{K_0} B^+_{st}$$
is canonically isomorphic to $H^0_{cris}( spec \, \overline R / spec \, W(k) ; \,  q^*G )$. This isomorphism commutes with the action of $\phi$.
 
\end{theorem} 

{\bf Proof of theorem  \ref{Bst}}. We keep the notation of the previous section. The morphism   $i_x$ extends to $\overline i _x: spec \, R \to \overline X$.
Choose a small open subset $ x \in U \subset \overline X$ and a lifting
$ \overline i_{\widetilde {x}}: Spf \, \overline R \to \widehat{\bar U}$.
 
We have canonical morphisms
$$i^*_x \mathbf {D} (E) \to  H^0_{cris}(\widehat{\bar U} /spec \, W(k) , p^*E|_{\widehat{U}}) \to 
H^0_{cris}( spec \, \overline R / spec \, W(k) ; \,  q^* \overline i^*_x E) \to $$
$$\Psi^{un}( \overline i^*_x E) \otimes _{K_0} B_{st} \to \overline i^*_x E_{\overline X}\otimes _{K} B_{DR} $$
 It is easy to see that the composition $i^*_x \mathbf {D} (E) \to i^*_x E_{\overline X}\otimes _{K} B_{DR}$ does not depend on the choice of the lifting 
$ \overline i_{\widetilde {x}}$. Hence it defines a morphism of functors: 
    $$i^*_x \mathbf {D} (*) \otimes _{\mathbf {Q}_p} B_{DR} \to i^*_x E_{\overline X}\otimes _{K} B_{DR} $$ 
 The morphism is compatible with the tensor structure and the trace morphism . Thus it is an isomorphism.

\section{ Unipotent variations of Hodge structure: a $p$-adic analog of the theorem by Hain and Zucker}
{\bf 13.1. Definition.}
{\it A variation of Hodge structure $E$ is called unipotent of level $r $ if 
there exists a filtration $0=E_0\subset E_1 \subset
\cdot \cdot \cdot \subset E_r=E $ such that the quotients 
$E_i/E_{i-1}$ are constant variations.}

 We denote by  ${\cal M}{\cal F}^{\mathbb{Q}}_{[a;a+b],r}(\overline X )$  the full subcategory of
the category  ${\cal M}{\cal F}^{\mathbb{Q}}_{[a;a+b]}(\overline X )$ consiting of unipotent variations of level $r$.
 It is easy to see that ${\cal M}{\cal F}^{\mathbb{Q}}_{[a;a+b],r}(\overline X )$  is an abelian category.

There is  the  functor to the category of vector bundles with an 
unipotent integrable connection:
$$  {\cal M}{\cal F}_{[a;a+b], r}(\overline X ) \longrightarrow {\cal D}{\cal M}_{r} (X_{K_0}) $$
For a point $x \in  X(W(k)) $ we let 
 $\Pi _{r} ^{DR}(X_{K_0},x,x) $ stand for the algebra of 
endomorphisms of the corresponding fiber functor 
 $${\cal D}{\cal M}_{r} (X_{K_0}) \longrightarrow Vect_{K_0}$$

\begin{lemma}\label{glupaya} 
a) $\Pi^{DR} _{r}      (X_{K_0},x,x) $ is a finite-dimensional algebra.

b) The category  ${\cal D}{\cal M}_{r} (X_{K_0}) $
is equivalent to the category of finite-dimensional 
  $\Pi _{r} ^{DR}(X_{K_0},x,x) $-modules.
\end{lemma}
Proof is omitted.
 
\begin{theorem} \label{Hain}
Assume that $r \leq \frac{p-1}{2}$
There exists a unique object $\Pi^{DR}_r (X_{K_0},x)$ 
of  ${\cal M}{\cal F}^{\mathbb{Q}}_{[-r+1;0],{r}}(\overline X )$, called the fundamental variation,
characterized by the property:

For any  unipotent variation $E \in {\cal M}{\cal F}^{\mathbb{Q}}_{[-r+1;0],{r}}(\overline X ) $ there is a canonical isomorphism
\begin{equation}
Hom_{{\cal M}{\cal F}^{\mathbb{Q}}_{[-r+1;0],{r}}(\overline X )}(\Pi^{DR}_r (X_{K_0},x);E)
\simeq Hom_{{\cal M}{\cal F}^{\mathbb{Q}}_{[-r+1;0],{r}}(spec \, W(k))}
(\underline{K_0}; E_x)
\label{krasota}
\end{equation}
( $\underline{K_0}$ stands for the trivial Hodge structure.) 
 
In addition, the variation $\Pi^{DR}_r (X_{K_0},x)$ enjoys the following properties:

i) The fiber $\Pi^{DR}_r (X_{K_0},x)_{x}$
is identified  with  $\Pi _{r} ^{DR}(X,x,x)$
Moreover,   the canonical action of  $\Pi _{r} ^{DR}(X,x,x)$  on the fiber (defined via the equivalence (\ref{glupaya})) is the multiplication from the left.

ii) $\Pi _{r} ^{DR}(X_{K_0},x,x)$ is an algebra in the category of Hodge structures.
That latter means that the multiplication  
$$\Pi _{r} ^{DR}(X_{K_0},x,x) \otimes _{K_0}\Pi _{r} ^{DR}(X_{K_0},x,x)  \to  \Pi _{r} ^{DR}(X_{K_0},x,x)$$
and the injection $ \underline{K_0} \hookrightarrow \Pi _{r} ^{DR}(X_{K_0},x,x) $  are morphisms of Hodge structures.
  
iii) For any object $E$ of   ${\cal M}{\cal F}_{[-r+1;0],{r}}^{\mathbb{Q}}(\overline X )$
the homomorhism
     $$ \Pi _{r} ^{DR}(X_{K_0},x, x) \longrightarrow {\cal H}om(E_x;E_x)$$
is amorphism of Hodge structures. (Note that the first assertion of ii) is a special case of iii) and, in fact, both of them immediately follow from the universal property
(\ref{krasota}) ).
 
Moreover, it defines   an equivalence between the category
${\cal M}{\cal F}_{[-r+1;0],{r}}
^{\mathbb{Q}}(\overline X )$ and the category of 
Hodge theoretic representations of 
 $\Pi _{r} ^{DR}(X_{K_0},x,x)$
(The latter consists of pairs $(V,f)$, where
$V$ is an object of ${\cal M}{\cal F}_{[-r+1;0],{r}}
^{\mathbb{Q}}(spec W(k) ) $ and a morphism
$f:\Pi_{r} ^{DR}(X_{K_0},x,x) \longrightarrow 
{\cal H}om(V;V)$  of the Hodge structures.)
\end{theorem}   
  
{\bf 13.2. Remark.}
{\it We would like to point out that the theorem  
is merely a $p$-adic analog  of the well known result about unipotent variations of mixed Hodge structure which is  due to  Hain 
and Zucker  [HZ].  }

A proof of the theorem occupies the next two sections.

\section{Construction of $p$-adic Hodge structure on $\Pi ^{DR}_r(X_K)$}
{\bf 14.1. } Remind, that $\Pi ^{DR}_r(X_{K_0})$ is a vector bundle on $X_{K_0} \times X_{K_0} $ 
together  with a unipotent connection of level $r$ uniquely
characterized by the following property:
 
given another object $E$ of the category ${\cal D}{\cal M}_{r} (X_{K_0} \times X_{K_0})$, we have a canonical isomorphism 
$\mathbf{D}$-modules on $X_{K_0}$:
\begin{equation}
Hom _{{\cal D}{\cal M} (X_{K_0} \times X_{K_0} / X_{K_0})}(\Pi ^{DR}_r(X_{K_0}); E) \simeq \Delta ^* E
\label{skoree}
\end{equation}   
Here $\Delta$ stands for the diagonal embedding $X_{K_0} \hookrightarrow X_{K_0} \times X_{K_0}$ and 
${\cal D}{\cal M} (X_{K_0} \times X_{K_0} / X_{K_0})$ is the category of vector bundles on $X_{K_0} \times X_{K_0}$ together
with an integrable connection along the fibers of the projection $p_2: \, X_{K_0} \times X_{K_0} \to X_{K_0} $  

In particularly, we have 
\begin{equation}
Hom _{{\cal D}{\cal M} (X_{K_0} \times X_{K_0}) }(\Pi ^{DR}_r(X_{K_0}); E) \simeq  H^0_{DR} (\Delta ; \Delta ^* E )
\label{domoy}
\end{equation}
The identity morphism $Id:\, \Pi ^{DR}_r(X_{K_0}) \to \Pi ^{DR}_r(X_{K_0}) $ gives rise to a parallel section 
\begin{equation}
 \mathbf{1}: \,
{\cal O}_{\Delta} \hookrightarrow \Delta ^* \Pi ^{DR}_r(X_{K_0}) 
\label{main}		
\end{equation}

{\bf 14.2.}The main result of this section is the following theorem. 
\begin{theorem}
There exists a unique variation  of $p$-adic Hodge structure on $\Pi ^{DR}_r(X_{K_0})$ with the following property:
 for any unipotent variation $E \in {\cal M}{\cal F}^{\mathbb{Q}}_{[-r+1;0],{r}}(\overline X \times \overline X) $
the map (\ref{skoree}) is a morphism of variations of Hodge structures.
\end{theorem}
{\bf 14.3. Proof.}
We construct the  Hodge structure on $\Pi ^{DR}_r(X_{K_0})$  by induction on $r$. First, we endow  $\Pi ^{DR}_1(X_{K_0})=
{\cal O}_{X_K\times X_K}$ with the trivial  Hodge structure.  Next, we assume that  the Hodge structure  on $\Pi ^{DR}_{r^{\prime}}(K)$
 is already constructed  for all $r^{\prime} < r$.

There is a  canonical surjective homomorphism (given by $ \mathbf{1}$)
$$  p_r : \Pi ^{DR}_r(X_{K_0}) \longrightarrow \Pi ^{DR}_{r-1}(X_{K_0})    $$
  
\begin{lemma}\label{krug}
a)  The connection on the vector bundle $ker \, p_r$ is trivial                                                                 .
 Moreover, there is a  canonical isomorphism:
  $$ (ker \, p_r)^* \simeq $$ 
$$ ker\, (Ext^1_{ {\cal D}{\cal M}(X_{K_0}\times X_{K_0})}(\Pi ^{DR}_{r-1}(X_{K_0}); 
{\cal O}_{X_K\times X_K})  \stackrel{\mathbf{1}}{\longrightarrow}
 Ext^1_{ {\cal D}{\cal M}(\Delta )}({\cal O}_{\Delta}  ; {\cal O}_{\Delta})) $$
b) The group of automorphisms of the sequence
\begin{equation} 
 0 \longrightarrow  ker \, p_r  \longrightarrow \Pi ^{DR}_r(X_{K_0}) \longrightarrow 
\Pi ^{DR}_{r-1}(X_{K_0}) \longrightarrow 0
\label{Zuk}
\end{equation}
identical on the boundary terms is isomorphic to $ker \, p_r $
\end{lemma} 
Proof is omitted.

By Proposition \ref{ri} and the induction hypothesis  $ ker \, p_r $ and $\Pi ^{DR}_{r-1}(X_{K_0})$
are endowed with Hodge structures.

{\bf Key Lemma.} {\it There exists a unique variation of Hodge structure on the middle term of (\ref{Zuk})
 satisfying the following properties:

i)All maps in the exact sequence (\ref{Zuk})
are compatible with Hodge structures.

ii) (\ref{main}) is a morphism of the variations of Hodge structures.}

\begin{proof}

We start with  the exact sequence from proposition \ref{ri}.
$$  0 \longrightarrow Ext^1_{{\cal M}{\cal F}^{\mathbb{Q}} (spec \, W(k) )} (\underline{K_0} ;
Hom_{ {\cal D}{\cal M}(X_{K_0}\times X_{K_0})}(\Pi ^{DR}_{r-1}(X_{K_0}) ; ker \, p_r ))$$
$$ \longrightarrow Ext^1_{ {\cal M}{\cal F}^{\mathbb{Q}}(\overline X \times \overline X)}
  (\Pi ^{DR}_{r-1}(X_{K_0}) ; ker \, p_r) 
\longrightarrow $$
\begin{equation}
( Ext^1 _{ {\cal D}{\cal M}(X_{K_0}\times X_{K_0})}
(\Pi ^{DR}_{r-1}(X_{K_0}) ; ker \, p_r))^{\phi=1}\cap F^0 
\longrightarrow 0 
\label{ugas}
\end{equation}
The vector space  
$$ Ext^1 _{ {\cal D}{\cal M}(X_{K_0}\times X_{K_0})}(\Pi ^{DR}_{r-1}(X_{K_0}) ; ker \, p_r)=$$
$$ Hom_{K_0} ( ker\, (Ext^1_{ {\cal D}{\cal M}(X_{K_0}\times X_{K_0})}(\Pi ^{DR}_{r-1}(X_{K_0}); 
{\cal O}_{X_K\times X_K}) \stackrel{\mathbf{1}}{\longrightarrow} $$
$$ Ext^1_{ {\cal D}{\cal M}(\Delta )}({\cal O}_{\Delta}  ;
 {\cal O}_{\Delta}) ); Ext^1_{ {\cal D}{\cal M}(X_{K_0}\times X_{K_0})}(\Pi ^{DR}_{r-1}(X_{K_0}); 
{\cal O}_{X_K\times X_K})                                      )$$ 
has a distinguished element: the identical morphism. It is easy to see that it coincides with the class of the extension (\ref{Zuk}).
 
In particularly, the latter class is invariant under $\phi$ and lies in $F^0$. 
Hence its preimage in $Ext^1_{ {\cal M}{\cal F}^{\mathbb{Q}}(\overline X \times \overline X)}
  (\Pi ^{DR}_{r-1}(X_{K_0}) ; ker \, p_r)$  is 
a torser ${\cal T}$  over the first $Ext$ group in the sequence (\ref{ugas}).
We have to show that the condition ii) gives a trivialization of the torser.

The fact that the connection on $ ker \, p_r $ is trivial implies that
$$ Hom_{ {\cal D}{\cal M}(X_{K_0}\times X_{K_0})}(\Pi ^{DR}_{r-1}(X_{K_0}) ; ker \, p_r )$$
$$ = Hom_{ {\cal D}{\cal M}(X_{K_0}\times X_{K_0} )}({\cal O}_{X_K\times X_K} ; ker \, p_r ))  = ker \, p_r $$
 Consider  the following commutative diagram:
 $$
\def\normalbaselines{\baselineskip20pt
\lineskip3pt  \lineskiplimit3pt}
\def\mapright#1{\smash{
\mathop{\to}\limits^{#1}}}
\def\mapdown#1{\Big\downarrow\rlap
{$\vcenter{\hbox{$\scriptstyle#1$}}$}}
\begin{matrix}{ Ext^1_{{\cal M}{\cal F}^{\mathbb{Q}} (spec \, W(k) )} (\underline{K_0} ;  ker \, p_r ) \mapright{} & 
 Ext^1_{ {\cal M}{\cal F}^{\mathbb{Q}}(\overline X \times \overline X)}
  (\Pi ^{DR}_{r-1}(X_{K_0}) ; ker \, p_r)  \cr
 \mapdown{Id} &  \mapdown{} \cr
Ext^1_{{\cal M}{\cal F}^{\mathbb{Q}} (spec \, W(k) )} (\underline{K_0} ;  ker \, p_r )
  \mapright{} & 
Ext^1_{ {\cal M}{\cal F}^{\mathbb{Q}}(\Delta )}(\Delta ^*\Pi ^{DR}_{r-1}(X_{K_0});ker \, p_r)   \cr
 \mapdown{Id} & \mapdown{}\cr 
  Ext^1_{{\cal M}{\cal F}^{\mathbb{Q}} (spec \, W(k) )} (\underline{K_0} ;  ker \, p_r )         
\mapright{Id} &   Ext^1_{ {\cal M}{\cal F}^{\mathbb{Q}}(\Delta)}({\cal O}_{\Delta };ker \, p_r)  }   
\end{matrix}
 $$

Here the third  vertical map is given by embedding  $ {\cal O}_{\Delta } \hookrightarrow \Delta ^*\Pi ^{DR}_{r-1}(X_{K_0})  $ .

It implies that there exists a unique element in  ${\cal T}$ which maps to $0$ in 
$  Ext^1_{ {\cal M}{\cal F}^{\mathbb{Q}}(\Delta)}({\cal O}_{\Delta};ker \, p_r)           $ i.e satisfies property ii) . 
 
Let  $\Pi ^{\prime}$ be the corresponding extension of the variations of Hodge structures.
It remaines to show  that there exists a unique parallel isomorphism between  the underlying $\mathbf{D}$-module and
$\Pi ^{DR}_{r}(X_{K_0})  $ identical on $\Pi ^{DR}_{r-1}(X_{K_0}) $ and $ker \, p_r$ and satisfying ii).
The existence follows from part b) of lemma (\ref{krug}) and uniqueness  from
the following sublemma.

{\bf Sublemma}.
 Fix  a variation of Hodge structures on $\Pi ^{DR}_{r}(X_{K_0})$ 
satisfying i) and ii). Then for any variation $E \in  {\cal M}{\cal F}^{\mathbb{Q}}_{[-r+1;0],{r}}(\overline X \times \overline X) $
we have a canonical isomorphism:
$$Hom_{{\cal M}{\cal F}^{\mathbb{Q}}_{r}(\overline X \times \overline X)}( \Pi ^{DR}_{r}(X_{K_0});E)\simeq 
Hom_{{\cal M}{\cal F}^{\mathbb{Q}}_{r}(spec W(k))}(\underline{K_0};  H^0_{DR} (\Delta ; \Delta ^* E ) )$$

{\it Proof of Sublemma}. 
Indeed, by the very definition the morphism (\ref{domoy}) is compatible with the Hodge structures.   Sublemma follows.
 
 This completes the proof of                      Key-lemma along with Theorem.
\end{proof}

{\bf 14.4.} Denote by 
$$p_{ij}: X_K \times  X_K
\times  X_K \longrightarrow 
 X_K \times  X_K$$
the projection given by the formula $p_{ij}(x_0,x_1,x_2)= (x_i;x_j)$.

We note that there is  a canonical parallel morphism:
\begin{equation}
p_{12}^*\Pi ^{DR}_r(X_{K_0}) \times p_{23}^*\Pi ^{DR}_r(X_{K_0}) \longrightarrow  p_{13}^*\Pi ^{DR}_r(X_{K_0})
\label{prod}
\end{equation}

The reader can easily check that (\ref{prod}) is compatible with the Hodge structures.

\section{Proof  of Theorem \ref{Hain} } 

{\bf 15.1.} Given a point $x\in X(W(k))$ we denote by $\Pi ^{DR}_r(X_{K_0}, x )$ the restriction of the variation $\Pi ^{DR}_r(X_{K_0})$ to the subscheme $\overline X \times x
\hookrightarrow  \overline X \times \overline X $.

It is easy to see that  $\Pi ^{DR}_r(X_{K_0}, x )$ satisfies the universal property (\ref{krasota}).

Next, given a unipotent variation of Hodge structure $E \in  {\cal M}{\cal F}_{[-r+1;0], r} ^{\mathbb{Q}}(\overline X )$ we consider the morphism 
 $$ \Pi ^{DR}_r(X_{K_0}, x ) \longrightarrow E^*\otimes E_x $$
given by $Id \in E_x^*\otimes E_x $. The latter gives rise to a homomorphism of Hodge structures 
$$ \Pi ^{DR}_r(X_{K_0}, x ,x) \longrightarrow   {\cal H}om( E_x ; E_x )$$

It remains to proof the last assertion of the theorem on the equivalence between the category $ {\cal M}{\cal F}_{[-r+1;0],r} ^{\mathbb{Q}}(\overline X )$
  and the category of Hodge theoretic representations $ \Pi ^{DR}_r(X_{K_0}, x ,x)$. 
\begin{lemma}
(\sf Rigidity.)Let  $E$ and $G$ be objects of ${\cal M}{\cal F}_{[a;a+b]} ^{\mathbb{Q}}(\overline X ) $, $2b \leq p-3$,
$f: E \longrightarrow G $ be a homommorphism of the underlying
coherent sheaves preserving the connection. Suppose there exists a point 
$x_0 \in  X(W(R))$ such that the restriction of $f$ to the point $x$
$f_x: E_x \longrightarrow G_x $ is a morphism of the Hodge structures.
Then $f$ is a morphism of the  variations of Hodge structures.
\end{lemma}
\begin{proof}
Let $P={\cal H}om(E;G)= E^*\otimes G $. It is easy to see that 
$P$ is an object of  ${\cal M}{\cal F}_{[-b;b]}^{ \mathbb{Q} }(\overline X )$.

Consider the group log cristalline cohomology $H_{cris}^0(\overline X ; P)$. 
One can check that it possesses a canonical Hodge structure (it also follows
from Theorem \ref{push}). Moreover, the natural map 
$$ r_x: H_{cris}^0(\overline X ; P) \to {\cal H}om(E_x;G_x)$$
 is an injection of Hodge structures. The homomorphism $f$ defines an element $f \in H_{cris}^0(\overline X ; P)$. Since 
$$ r_x (f) \in ker(\phi_0 -Id: F^0({\cal H}om(E_x;G_x)) \to {\cal H}om(E_x;G_x))$$
 and morphisms of Hodge structures are strictly compatible with the filtration, we have
$f \in  ker(\phi_0 -Id:F^0(H_{cris}^0(\overline X ; P)) \to
H_{cris}^0(\overline X ; P))$. It completes the proof.
\end{proof}
 The lemma implies that that  the functor is fully faithful. It remains  show that its image contains all Hodge theoretic representations.
Given a representation $\rho: \Pi ^{DR}_r(X_{K_0}, x ,x) \otimes V \to V $,
we define a variation of Hodge structure to be 
       $$ coker \,( \Pi ^{DR}_r(X_{K_0}, x )   \otimes  \Pi ^{DR}_r(X_{K_0}, x ,x)        \otimes V 
\stackrel{m\otimes 1 - 1\otimes \rho}{\longrightarrow} \,( \Pi ^{DR}_r(X_{K_0}, x ) \otimes V ))$$
Here $m$ stands for the canonical morphism: 
$$   \Pi ^{DR}_r(X_{K_0}, x ) \otimes  \Pi ^{DR}_r(X_{K_0}, x ,x) \to      \Pi ^{DR}_r(X_{K_0}, x )$$
It defines the inverse functor.

\section{Proof of Theorems A}
{\bf 16.1.} We denote by $Sh^{et}_{\mathbb{Q}_p, r}(X_{K_0}) $ the full subcategory of
 $Sh^{et}_{\mathbb{Q}_p}(X_{K_0})$ which consists of unipotent etale local systems of level $ r$
( A local system ${\cal E}$ is called  unipotent  of level $r$ if there exists a filtration 
 $0= {\cal E}_0 \to {\cal E}_1 \to \cdot \cdot \cdot \to {\cal E}_r= {\cal E}$ 
such that ${\cal E}_i/ {\cal E}_{i-1}$ is isomorphic to the pullback of a sheaf on $spec\, K_0$)
 
Let $\Pi^{et}_r(X_{K_0})$ be the fundamental local system on $X_{K_0}\times X_{K_0}$. It is characterized by the property that for any 
${\cal E} \in Sh^{et}_{\mathbb{Q}_p, r}(X_{K_0} \times X_{K_0} )$ there is a canonical isomorphism:
$$Hom_{ Sh^{et}_{\mathbb{Q}_p, r}(X_{K_0} \times X_{K_0} )}(\Pi^{et}_r(X_{K_0}); {\cal E}) \simeq H^0_{et}(\Delta ; \Delta ^* E)$$
 \begin{theorem}
Assume that $p\leq \frac{p-1}{2}$.
There is a canonical isomorphism:
 $$\mathbf{D}(\Pi ^{DR}_r(X_{K_0})) \simeq \Pi^{et}_r(X_{K_0})$$
\end{theorem}
\begin{proof}
The canonical morphism 
$$ \mathbf{1}: \, {\cal O}_{\Delta} \hookrightarrow \Delta ^* \Pi ^{DR}_r(X_{K_0})$$ 
induces  a map 
$$ \underline {\mathbb{Q}_p} \longrightarrow \Delta ^* \mathbf{D}(\Pi ^{DR}_r(X_{K_0}))$$
which, in turn, gives rise to a morphism:
 $$f: \,  \Pi^{et}_r(X_{K_0}) \longrightarrow  \mathbf{D}(\Pi ^{DR}_r(X_{K_0})) $$
First, we show that $f$ is surjective. Since $\mathbf{D}$ is fully faithfull, it suffices to check that for any 
$E \in {\cal M}{\cal F}_{[-r+1;0], r} ^{\mathbb{Q}}(\overline X \times \overline X) $ the canonical morphism:
$$A:\, Hom_{ Sh^{et}_{\mathbb{Q}_p, r}(X_{K_0} \times X_{K_0} )}(\mathbf{D}(\Pi^{DR}_r(X_{K_0})); \mathbf{D}(E)) \longrightarrow $$
$$ Hom_{ Sh^{et}_{\mathbb{Q}_p, r}(X_{K_0} \times X_{K_0} )}( Im \, f ; \mathbf{D}(E))$$
is isomorphism.
It follows from the  universal property of $Pi^{DR}_r(X_{K_0})$ that
$$Hom_{ Sh^{et}_{\mathbb{Q}_p, r}(X_{K_0} \times X_{K_0} )}(\mathbf{D}(\Pi^{DR}_r(X_{K_0})); \mathbf{D}(E)) \simeq $$
$$Hom_{{\cal M}{\cal F}_{[-r+1;0]} ^{\mathbb{Q}}( spec \, W(k))}(\underline{K_0} ; H^0_{DR} (\Delta ; \Delta ^* E ))$$
On the other hand, we have
$$Hom_{ Sh^{et}_{\mathbb{Q}_p, r}(X_{K_0} \times X_{K_0} )}(\Pi^{et}_r(X_{K_0}); \mathbf{D}(E)) \simeq
 H^0_{et}(\Delta ; \Delta ^* \mathbf{D}E) $$
Finally, combining the latter with the isomorphism 
$$Hom_{{\cal M}{\cal F}_{[-r+1;0]} ^{\mathbb{Q}}( spec \, W(k))}(\underline{K_0} ; H^0_{DR} (\Delta ; \Delta ^* E ))\simeq
 H^0_{et}(\Delta ; \Delta ^* \mathbf{D}( E)) $$ we obtain the inverse map $A^{-1}$.
It proves subjectivity of $f$. Since $$dim_{\mathbb{Q}_p}\Pi^{et}_r(X_{K_0})= dim_{K_0}\Pi ^{DR}_r(X_{K_0})$$, $f$ is an isomorphism. 

\end{proof}
We can make use of Theorem  \ref{Bst} to derive the following result.
\begin{cor} 
Given a pair of points $x_0, x_1 \in X_{K_0}(\overline K)$, we have a canonical isomorphism
$$\Pi^{et}_r(X_{K_0},x_0, x_1) \times _{\mathbb{Q}_p} B_{DR}  \simeq \Pi^{DR}_r(X_{K_0},x_0, x_1)\times _{\overline K} B_{DR}$$
\end{cor}
It proves Theorem A.

{\bf 16.2. Remark.} It seems likely that one can make use of the $\mathbb{Q}_p$-theory [Fa] (associated convergent $F$-isocrystals ...) and generalize the above argument
to prove Theorem A for any smooth variety over a finite extension of $\mathbb{Q}_p$ (with no restrictions on $p$).     

\section{P-adic integration}
{\bf 17.1.}Let $X$ be a smooth variety over $\mathbb{C}$ and $E$ be a vector bundle together with integrable connection 
$$\nabla : \,   E \longrightarrow E\otimes \Omega^1_X $$
Parallel translation along a path
\begin{equation}
T_{\lambda}: E_{x_0}  \simeq  E_{x_1}
\label{translat}
\end{equation} 
is a generalization of integration of $1$-forms. Let us formulate it explicitly.    

Given a $1$-form $\omega$, we define  $E$ to  be the  trivial vector bundle on $X$ with a basis $e_{-1}, e_{0}, $
 and a connection $\nabla : E \to E\otimes \Omega^1_X $ is given  
in this basis by the formula:
                   $$\nabla(s)= A\, s + ds $$
where      
 $$
A=
\left(\begin{array}{cc}
0 & \omega  \\
0 & 0 
\end{array}\right)
$$ 
Then the parallel translation (\ref{translat}) is given by the following matrix 
\begin{equation}
T_{\lambda}=
\left(\begin{array}{cc}
1 & -\int\limits_{\lambda} \omega  \\
0 & 1 
\end{array}\right)
\label{formula} 
\end{equation}
{\bf 17.2.} Next, we let $X_K$ to be a smooth geometrically connected scheme over $K$ ,$\omega$ be a closed $1$-form and $x_0, x_1 \in X_K{K}$ 
 We  make use of the canonical parallel translation $C_{x_0, x_1, X_K}$ (\ref{canpar}) and  (\ref{formula}) {\sf to define } a $p$-adic  integral
$\int\limits^{x_1}_{x_0} \omega  \in K_{st}$.
\begin{proposition}
The integral $\int\limits^{x_1}_{x_0} \omega $ coincides with Colmez's integral.
\end{proposition}
It follows from Conjecture {\bf 1} (see ({\bf 1.18.}) proven in rank two case by Colmez.

\section{Appendix}
{\bf 18.1.}In this Appendix we prove  the Monodromy Conjecture for $H^2_{pst}(X_{\overline K})$.        
\begin{theorem}
Let $X_{\overline K}$ be a  proper smooth scheme finite type over $spec \, \overline K$.
Then, for $n \leq 2 $,  the ${\cal H}^{ur}_1$-module  $H^n_{pst}(X_{\overline K})$ is pure of weight $n$. 
(For the definition of purity,
see {\bf 7.1.}) 
\end{theorem}
\begin{cor}
Let $X_{\overline K}$ be a  smooth scheme finite type over $spec \, \overline K$. Then, for $n \leq 2 $,
$(H^n_{pst}(X_{\overline K}); W_{\cdot})$ is mixed.
\end{cor}
{\bf 18.2. Weakly admissible modules.} The definitions and results recollected bellow are due to Fontaine [FoIl].
   
A {\sf filtered $(\phi, N)$-module over $K$} is a ${\cal H}_1$ - module $V$ together with a decreasing filtration 
$F^iV$, $i\in \mathbb{Z}$ on $V_K=V\otimes_{K_0} K$ satisfying $\cup_i F^iV=V$ and $\cap_i F^iV=0$.
The filtered  $(\phi, N)$- modules form an additive category ${\mathbf MF}_K(\phi, N)$. 

For any object $V$ of  ${\mathbf MF}_K(\phi, N)$ of dimension $1$, we define $t_H(V)$ to be the biggest integer $i$ such that $F^iV_K=V_K$; 
if $\phi v = \lambda v$  define  $t_N(V)$ as the $p$-adic valuation of $\lambda$.  For any object $V$, whose underlying vector space is finite-dimensional
we put
                              $$t_H(V)=t_H(\wedge^{top}V) \; , \;   t_N(V)=t_N(\wedge^{top}V)$$ 
A subspace $V^{\prime}\subset V$  of an object of  ${\mathbf MF}_K(\phi, N)$ is called 
a subobject if  $V^{\prime}$ is  stable under the action of ${\cal H}_1$. We endow it with the filtration $F^iV_K^{\prime}:=F^iV\cap V_K^{\prime}$.

A finite-dimensional   object   $V$  is called weakly admissible if:\\
i)$t_H(V)= t_N(V)$\\
ii) For any subobject $V^{\prime}\subset V$, $t_H(V^{\prime} )\leq t_N(V^{\prime})$.

It is known that the category ${\mathbf MF}^a_K(\phi, N)$ of  weakly admissible modules is Tannakian 
(in particularly, abelian).

The latter result implies the weak admissibility  is preserved by the extension of
scalars i.e. we have a functor: 
        $${\mathbf MF}^a_K(\phi, N) \longrightarrow  {\mathbf MF}^a_L(\phi, N)$$
where $L\supset K$ is a finite extension. It allows one to define the category ${\mathbf MF}^a_{\overline K}(\phi, N)$.

 It follows from
results of [Tsuji]  that  for any smooth proper scheme finite type over $spec \, \overline K$
 the filtered $(\phi, N)$-module  $H^i_{pst}(X_{\overline K})$ (the filtration is the Hodge filtration on
$H^i_{pst}(X_{\overline K}) \otimes _{\mathbb{Q}_p^{ur}} \overline K \simeq H^i_{DR}(X_{\overline K}) $)
is  weakly admissible.

 {\bf 18.3. Lemma.} {\it Let $f: P\to Q$ be a morphism of finite-dimensional ${\cal H}^{ur}_1$-modules and
$V: = ker \, f$ . Assume that $P$ is pure of weight $n$ and $t_N(V)=\frac{n \cdot  dim \, V}{2}$.
Then $V$ is pure of weight $n$.}

{\bf Proof:} It suffices to prove the statement for a morphism of finite-dimensional ${\cal H}_1$-modules.
 Fix an integer $c$, such that $Fr^c$ acts trivially on $K_0$. We have the decomposition 
 $$V=\bigoplus_iV^i$$
 where  $V^i$ is the maximal $\phi$-invariant subspace such that all eigenvalues 
of $\phi^c|_{V^i}$ are Weil numbers of weight $ci$.
We have  to  show that for any $i$, the map
                         $$N^i : \, V^{n+i} \longrightarrow  V^{n-i}$$
is an isomorphism. Assume that this not the case. 
 It follows from the purity of  $P$ 
 that  the latter map is injective. Hence, the assumption implies that the eigenvalue of    $\phi^c$ on 
 the eigenvalue of the operator $\phi^c : \wedge^{top}V \to \wedge^{top}V $ is a  Weil number $\alpha $ of 
weight $\leq c(n \cdot dim \, V - 2)$.        
On the other hand,  the $p$-adic valuation of $\alpha$ is equal to $c\frac{n\cdot dim \, V}{2}$.
 It follows that $p^{-c\frac{n\cdot dim \, V}{2}} \alpha $ is an algebraic integer. 
On the other hand, it is a Weil number of negative weight. The contradiction 
completes the proof.

{\bf 18.4. Proof of the Theorem.} By the Chow Lemma we can find a smooth projective variety 
$Y_{\overline K}\subset \mathbb{P}^N $ and
a proper morphism $f:Y_{\overline K} \to   X_{\overline K}$, which is a birational equivalence.
Next, we  choose a plane $H\subset \mathbb{P}^N $ such that $S: = H\cap Y$ is a smooth connected surface.

The Weak Lefschetz Theorem implies that, for $n \leq 2$, the map
     $$ r: \;  H^n_{pst}(X_{\overline K})\longrightarrow H^n_{pst}(S)$$
is injective.
By the result of [Mokrane] the module $H^n_{pst}(S)$ is pure of weight $n$.
Since both modules  are weakly admissible, so is $coker \, r$. On the other hand, it is easy to see that
$t_H(coker \, r)=\frac{ n \cdot dim \, coker \, r}{2}$. 
(For example, we can embed $\overline K$ into $\mathbb{C}$ and use the fact that
$coker \, r$ possesses a Hodge structure of weight $2$). Hence, $t_N(coker \, r)= \frac{ n \cdot dim \, coker \, r}{2}$.
 It remains to 
apply the lemma (to the dual map).

\end{document}